\documentclass[11pt, twoside]{article}
\usepackage{amsmath, amsthm, amssymb, bbm, mathrsfs, graphicx, tikz-cd, calligra, enumerate, titlesec, fancyhdr}
\usepackage[hidelinks]{hyperref}
\usepackage[margin=1in,marginparwidth=0.8in]{geometry}

\numberwithin{equation}{section}

\newtheorem{theorem}[equation]{Theorem}
\newtheorem{lemma}[equation]{Lemma}
\newtheorem{proposition}[equation]{Proposition}
\newtheorem{corollary}[equation]{Corollary}

\theoremstyle{definition}
\newtheorem{definition}[equation]{Definition}
\newtheorem{example}[equation]{Example}
\newtheorem{construction}[equation]{Construction}

\newtheorem{remark}[equation]{Remark}
\newtheorem*{remark*}{Remark}

\titleformat{\section}
  {\normalfont\Large\centering}
  {\thesection}{1em}{}

\titleformat{\subsection}
  {\normalfont\large\bfseries}
  {\thesubsection}{1em}{}

\title{Motivic Borel-Moore Homology is Borel-Moore Motivic Homology}
\author{Spencer G. Taylor\footnote{University of California, Berkeley. }}
\date{\today}

\newcommand{\bb}[1]{\mathbb{#1}}
\newcommand{\p}{\mathfrak{p}}

\newcommand{\Hom}{\mathrm{Hom}}
\newcommand{\im}{\mathrm{im} \,}
\newcommand{\id}{\mathrm{id} }
\newcommand{\colim}{\varinjlim}
\newcommand{\spec}{\mathrm{Spec} \,}

\newcommand{\catlim}{\varprojlim}

\newcommand{\mf}[1]{\mathfrak{#1}}

\newcommand{\codim}{\mathrm{codim}}
\newcommand{\trd}{\mathrm{trdeg}}
\newcommand{\shom}{\mathcal{H} \kern - 3pt \text{\calligra \large om} \,}

\newcommand{\pr}{\mathrm{pr}}
\newcommand{\map}{\mathrm{map}}
\newcommand{\fil}{\mathrm{Fil}}
\newcommand{\gr}{\mathrm{gr}}
\newcommand{\Kgl}{\mathrm{KGL}}
\newcommand{\kgl}{\mathrm{kgl}}
\newcommand{\Ggl}{\mathrm{GGL}}
\newcommand{\sh}{\mathrm{SH}}
\newcommand{\Sp}{\mathrm{Sp}}

\newcommand{\bm}{\mathrm{BM}}
\newcommand{\psh}{\mathrm{PSh}}
\newcommand{\shv}{\mathrm{Sh}}
\newcommand{\lan}{\mathrm{Lan}} 
\newcommand{\fun}{\mathrm{Fun}}
\newcommand{\nis}{\mathrm{Nis}}
\newcommand{\sm}{\mathrm{Sm}}
\newcommand{\G}{\mathrm{G}}
\newcommand{\sli}{\mathrm{slice}}
\newcommand{\KH}{\mathrm{KH}}
\newcommand{\K}{\mathrm{K}}
\newcommand{\Map}{\mathrm{Map}}
\newcommand{\op}{\mathrm{op}}
\begin{document}
\let\placephi\phi
\let\phi\varphi
\let\varphi\placephi
\let\tilde\widetilde
\let\hat\widehat
\maketitle

\begin{abstract}
This paper studies the higher Chow groups of arbitrary finite type separated schemes using motivic homotopy theory, serving as a Borel-Moore counterpart to the \(\bb{A}^1\)-invariant motivic cohomology of Bachmann-Elmanto-Morrow. We construct a filtration on the \(\G\)-theory spectrum of a finite type separated scheme \(X\), inducing an Atiyah-Hirzebruch style spectral sequence whose associated graded pieces are a complex computing the higher Chow groups of \(X\). We prove that this complex is representable in the stable motivic homotopy category, and is a module over the \(\bb{A}^1\)-motivic complex, giving rise to an action of \(\bb{A}^1\)-invariant motivic cohomology on the Bloch-Levine cycle complex. This generalizes some results of Bloch and Levine relating higher Chow groups and motivic cohomology to the case of schemes over an arbitrary base, and which need not be smooth or equidimensional.
\end{abstract}

\tableofcontents

\section{Introduction:}
Over the last few decades, motivic cohomology has slowly been constructed in increasingly general situations. The first such construction, provided by Bloch, was that of the higher Chow group of a quasiprojective finite type separated scheme over a field  \cite{Blo}. Voevodsky's construction of motivic cohomology for finite type schemes over a field as the hypercohomology groups of certain complexes \cite{MVW} satisfied many of Beilinson and Lichtenbaum's proposed properties \cite{BeiMot, LichZet}. He also proved that this definition agreed with the higher Chow groups for smooth schemes over a field, and constructed so-called motivic Borel-Moore homology groups, which agreed with higher Chow groups for equidimensional schemes over a field which were not smooth.

At this point, the theory diverged in two directions. Levine generalized Bloch's cycle complex construction to any finite type scheme over a regular ring of Krull dimension one, defining the motivic Borel-Moore homology as such the hyperhomology of this cycle complex \cite{LevTechLoc}. He also constructed a spectral sequence with \(E^2\)-page the motivic Borel-Moore homology groups, and proved its convergence to the \(\G\)-groups of the scheme in question. On the other hand, Bachmann-Elmanto-Morrow extended the definition \(\bb{A}^1\)-invariant) motivic cohomology to any quasicompact quasiseparated scheme \cite{BEM}. Their motivic complex sidesteps many of the delicate cycle-theoretic methods of earlier constructions by working in the stable motivic homotopy category, obtaining it a mapping spectrum from the unit to the pullback of slices of the motivic spectrum representing homotopy \(\K\)-theory. This construction also provides a filtration of homotopy \(\K\)-theory inducing an algebraic version of the Atiyah-Hirzebruch spectral sequence, whose associated graded pieces are precisely the \(\bb{A}^1\)-invariant motivic cohomology groups.

This paper seeks to close the loop, obtaining the global sections of Levine's cycle complex as a complex defined in an analogous manner: a mapping spectrum in the stable motivic homotopy, this time from the unit to the \textit{exceptional} pullback of slices of the motivic spectrum representing homotopy \(K\)-theory. We arrive at a similar filtration on \(G\)-theory, with an associated spectral sequence whose associated graded pieces are motivic Borel-Moore homology. Additionally, we show that the action of homotopy \(\K\)-theory on \(\G\)-theory descends to the level of slices, so we may thus think of Levine's cycle complex as an object on which \(\bb{A}^1\)-motivic cohomology acts.

\subsection{Summary of results} 
\begin{theorem}\label{main-intro} Let \(X\) be a finite type separated scheme. There exists a complex 
\[\bb{Z}(n)^\bm(X) \in \mathrm{D}(\bb{Z}), \quad n \in \bb{Z}\]
called \emph{motivic Borel-Moore homology} satisfying the following properties: 
\begin{enumerate}
\item Borel-Moore functoriality: for smooth \(g \colon Y \to X\) of relative dimension \(d\), there is a morphism 
\[g^* \colon \bb{Z}(n)^\bm(X) \to \bb{Z}(n + d)^\bm(Y)\lbrack-2d \rbrack,\]
and for proper \(p \colon X \to C\), there is a morphism 
\[p_* \colon \bb{Z}(n)^\bm(X) \to \bb{Z}(n)^\bm(C)\]

\item Nisnevich descent: \(\bb{Z}(n)^\bm\) converts elementary distinguished squares of schemes to Cartesian squares of spectra.

\item Atiyah-Hirzebruch filtration: there exists a natural, exhaustive (decreasing) filtration \(\fil_\star^\bm \G(X)\) on \(G\)-theory with graded pieces are given by 
\[\gr_n^\bm \G(X) \simeq \bb{Z}(n)^\bm(X)\lbrack 2n \rbrack,\]
Writing \(H_p^\bm \big(X, \bb{Z}(q)\big) = H_p \big(\bb{Z}(q)^\bm(X)\big)\) for the \emph{motivic Borel-Moore homology groups}, then there exists an Atiyah-Hirzebruch-style spectral sequence 
\[E_{p, q}^2 = H_{p + q}^\bm\big(X, \bb{Z}(- q)\big) \Longrightarrow G_{p + q}(X)\]

\item Module structure over \(\bb{A}^1\)-motivic cohomology: \(\bb{Z}(- \star)^\bm(X)\) is a graded module over the graded \(\bb{E}_\infty\)-ring \(\bb{Z}(\star)^\bb{A}(X)\), in \(\mathrm{Gr} \, \mathrm{D}(\bb{Z})\), with action 
\[\bb{Z}(i)^\bb{A}(X) \otimes \bb{Z}(n)^\bm(X) \to \bb{Z}(n - i)^\bm(X)\]
descending from the action of spectra \(\KH(X) \wedge \G(X) \to \G(X)\).

\item Comparison to the Bloch-Levine cycle complex: for any \(n \in \bb{Z}\), there is an equivalence of complexes 
\[\bb{Z}(n)^\bm(X)\lbrack -2n \rbrack \simeq R\Gamma(X_\mathrm{Zar}, \mathcal{Z}_{n + 1}(X, *)\big),\]
Taking homology thus gives an isomorphism with higher Chow groups: 
\[H_j^\bm\big(X, \bb{Z}(n)\big) \cong \mathrm{CH}_{n + 1}(X, j - 2n)\]

\end{enumerate}
\end{theorem}

\subsection{Notation and conventions}
We write \(\Sp\) for the \(\infty\)-category of spectra. Given objects \(X, Y\) in an \(\infty\)-category \(\mathcal{C}\), we denote by \(\Map_\mathcal{C}(X, Y) \in \Sp\) the mapping space between the the objects. When, in addition, \(\mathcal{C}\) is stable, we write \(\map_\mathcal{C}(X, Y)\) for the mapping spectrum. We call a morphism of schemes \(f \colon X \to S\) admitting a factorization 
\[\begin{tikzcd} X \ar[r, "g"] & Z \ar[r, "i"] & S \end{tikzcd}\]
with \(g\) smooth and \(i\) a regular closed embedding a \emph{locally complete intersection (lci)} morphism.\footnote{Note that this is called a \emph{smoothable} lci morphism in other works}

Given a scheme \(X\), and point \(x \in X\), \(Z(x)\) denotes the closure \(\overline{\{x\}}\) with the reduced closed subscheme structure. Every morphism is finite type separated unless otherwise stated.

When \(D\) is some object in a category \(\mathscr{C}\), we write \(C_D\) to denote the corresponding constant simplicial object in the category \(\fun(\Delta^\op, \mathscr{C})\), and do the same for morphisms.

\subsection{Outline}
Here, we briefly outline the content of each section.

\begin{enumerate}
\item[\S\ref{prelim}] We discuss the stable motivic homotopy category, and the six-functor formalism it inherits. A key construction in this cagegory is the slice filtration, a motivic analogue of the Postnikov filtration of a spectrum. We review the motivic spectra representing homotopy \(K\)-theory and its role defining \(\bb{A}^1\)-invariant motivic cohomology and its equivalent \(\bb{A}^1\)-cdh counterpart. We then give a brief exposition of the representability of \(G\)-theory, bivariant theories, and their relation to homotopy K-theory. We finish with an overview of the Bloch-Levine cycle complex.

\item[\S\ref{G-fil}] In this section, we define a filtration on \(G\)-theory, and describe the convergence properties of the associated spectral sequence. We show that the action of \(\Kgl\) on \(\Ggl\) descends to the level of these filtrations.

\item[\S\ref{slice-comp-A}] This section beings the proof of the first half of main theorem, constructing a presheaf of spectra with the same derived global sections as the Bloch-Levine complex.

\item[\S\ref{slice-comp-B}] We continue the work of the previous section, building niveau filtrations on the spectral replacement and motivic Borel-Moore complex. We construct a filtered morphism between them, prove that it induces an equivalence on derived global sections by reducing to the case of smooth schemes over perfect fields, finishing the main theorem.

\end{enumerate}

\subsection{Acknowledgements}
The author would like to thank Elden Elmanto for suggesting this problem, for helpful conversations and advice, and for inviting him to give a talk on this paper. Without Elden's response to a cold email, literature recommendations, and patience with numerous questions, this paper never would have come to be. The author is also extremely grateful for his advisor David Nadler, who has provided constant support, mathematical and otherwise, despite numerous periods of the author falling off the map and continually running behind.

\vspace{1em}
\noindent
\textbf{LLM Statement:} An LLM was used to search literature and check for spelling errors and notational inconsistencies. An LLM also suggested an outline for the statement and proof of lemmas \ref{lemm:fun-fact} and \ref{lemm:supp-fun}, and construction \ref{con:ex-supp} that greatly shortened the much messier original statements and proofs of the author. Every word, symbol, and diagram of this paper were written by the author.

\section{Preliminaries}\label{prelim}
In this section, we introduce the underlying machinery needed to state and prove our results, and review the key results on which we build.

\subsection{The stable motivic homotopy category}
This construction follows \cite{Rob}, which uses modern \(\infty\)-categorical language. Let \(X\) be a quasicompact quasiseparated scheme. There exists a stable presentable symmetric monoidal \(\infty\)-category \(\sh(X)\), the \emph{stable motivic homotopy category} of \(X\), whose monoidal product we denote \(\otimes\), and refer to as the \emph{smash product} in analogy to algebraic topology. It is equipped with a symmetric monoidal functor \(M_X : \sm_X \to \sh(X)\), such that 
\begin{enumerate}
\item \(\mathbbm{1}_X = M_X(X)\) is the unit in the symmetric monoidal category \(\sh(X)\).

\item The morphism \(M_X(\bb{A}_X^1) \to M_X(X)\) induced by \(\bb{A}_X^1 \to X\) is an equivalence.

\item Smashing by \(\bb{T}_X = \mathrm{fib} \! \left(M_X\big(\bb{P}_X^1\big) \to \mathbbm{1}_X\right)\) is invertible in \(\sh(X)\). This object, whose subscript we often omit, is called the \emph{Tate motive}.

\item \(M_X\) sends a Nisnevich square in \(\sm_X\) into a Cartesian square in \(\sh(X)\).
\end{enumerate}
The pair \((\sh(X), M_X)\) is initial with respect to these properties.

\begin{remark}[Suspension] For a motivic spectrum \(E \in \sh(X)\), we denote by  
\[E \mapsto E[1] = \mathrm{cofib}(E \to 0), \quad E \mapsto E[- 1] = \mathrm{fib}(0 \to E)\]
the \emph{suspension} and \emph{loop} functors, respectively. Stability of \(\sh(X)\) means these are inverse equivalences. Motivic homotopy theory oftentimes uses the notation \(S^{a, b} = \bb{T}^{\otimes b}[a - 2b]\) for bigraded spheres, and \(\Sigma^{a, b}E = S^{a, b} \otimes E\) for bigraded suspensions. Following \cite{BEM}, we will avoid this and write everything in terms of \(\bb{T}\) and suspensions.

\end{remark}
\begin{remark}[Six functor formalism, following \cite{AY1, AY2}]
For every morphism \(f \colon X \to S\) of QcQs schemes (not necessarily finite type separated), there is an adjunction
\[f^* \colon \sh(S) \rightleftarrows \sh(X) \colon f_*\]
called \emph{pullback} and \emph{pushforward}, respectively. If \(f\) is smooth, \(f^*\) has a left adjoint: 
\[f_\sharp \colon \sh(X) \rightleftarrows \sh(S) \colon f^*\]
There is a second adjoint pair 
\[f_! \colon \sh(X) \rightleftarrows \sh(S) \colon f^!,\]
called \emph{exceptional pushforward} and \emph{pullback}, respectively. Finally, we have, for any \(E \in \sh(X)\), there is a third adjoint pair 
\[\underline{\Hom}(E, -) \colon \sh(X) \rightleftarrows ( - )\otimes E,\]
giving \(\sh(X)\) a closed symmetric monoidal structure. These functors satisfy the following compatibilities:
\begin{enumerate}
\item \(f^*\) is symmetric monoidal, and given \(W \in \sm_S\), there is an equivalence  
\[f^* M_S(W) \simeq M_X\left(X \times_S W\right)\footnote{Note that taking \(W = S\), this implies that \(f^*\mathbbm{1}_S \simeq \mathbbm{1}_X\)}.\]

\item For any open embedding \(j\), \(j_\sharp \simeq j_!\) and \(j^! \simeq j^*\).

\item For \(p\) proper, \(p_! \simeq p_*\).

\item For any closed embedding \(i\),  \(i_*\) is fully faithful.

\item Given a Cartesian square 
\[\begin{tikzcd} \tilde{X} \ar[d, "\tilde{g}", swap] \ar[r, "\tilde{f}", ""{name=h1}] & \tilde{Y} \ar[d, "g"]\\
X \ar[r, "f", swap, ""{name=h2}] & Y 
\ar[phantom, from=h1, to=h2, "\Delta"]
\end{tikzcd}\]
there are canonical equivalences 
\[g^* f_! \Longrightarrow \tilde{f}_! \tilde{g}^*, \quad \tilde{g}_* \tilde{f}^! \Longrightarrow f^! g_*,\]
called \emph{base change}. 

\item Any pair of motivic spectra \(E, F \in \sh(S)\) satisfy a \emph{projection formula}: there is a morphism 
\[f^!E \otimes f^*F \to f^!\big(E \otimes F\big)\]
inducing a lax \(\sh(S)\)-linear structure on \(f^!E\). When \(E\) has ring structure, this gives \(f^!E\) the structure of a \(f^*E\)-module.

\end{enumerate}
\end{remark}

\begin{remark}[Purity]
When \(f\) is a smooth morphism, there is a canonical equivalence of functors, 
\[\p_f \colon \mathrm{Th}(T_f) \otimes f^* \simeq f^!,\]
For \(T_f\) the relative tangent bundle. When \(i \colon X \to Y\) is a closed embedding, and \(f_X, f_Y \colon X, Y \to S\) are smooth, then there is an equivalence 
\[\mathrm{Th}(- N_X Y) \otimes f_X^* \simeq i^! f_Y^*\]
These equivalences are called \emph{smooth} and \emph{relative purity}, respectively. Very frequently (any time we work over slices of \(\Kgl_\bb{Z}\), for example), the tangent and normal bundles are trivializable, in which case their Thom spectra are some power of the Tate motive.
\end{remark}

\subsection{KGL and \(\bb{A}^1\)-motivic cohomology}

\begin{definition}[Filtered and graded objects]
To a stable \(\infty\)-category \(\mathcal{C}\), we may associate the stable \(\infty\)-categories of \emph{filtered objects} and \emph{graded objects}, 
\[\fil \, \mathcal{C} = \fun\big((\bb{Z}, \geq)^\op, \mathcal{C}\big) \quad \text{ and } \mathrm{Gr} \, \mathcal{C} = \fun(\bb{Z}^\delta, \mathcal{C}),\]
in which \((Z, \geq)\) is the totally ordered set of the integers, and \(\bb{Z}^\delta\) is the discrete category of the integers. A filtration on \(E \in \mathcal{C}\) is called \emph{exhaustive} when 
\[\fil^{- \infty} E = \colim_{j \to - \infty} \fil^j E \to E\]
is an equivalence, and \emph{complete} when 
\[\catlim_{j \to \infty} \fil^j E \simeq 0\]
\end{definition}

\begin{remark}[Increasing filtrations]
The slice filtration and filtration on \(\KH(X)\) are decreasing, but the filtration we construct on \(\G(X)\), as well as the niveau filtration, which plays a key role in proving the equivalence \ref{main-intro} (5), are increasing. As such, we write these filtrations and their associated gradeds with lower indices. We will abusively denote the categories of filtered and graded objects as above, but there will never be any ambiguity.
\end{remark}

\begin{definition}[The slice filtration]
The classical reference for this filtration, which are perfectly well-suited for our class of schemes, are \cite{VSli1} and \cite{VSli2}. The \(\infty\)-category of \textit{effective spectra}, \(\sh(S)^\mathrm{eff}\) is the smallest stable subcategory of \(\sh(S)\) closed under colimits and containing \(M_S(Y)[i]\) for all \(Y \in \sm_S\) and \(i \in \bb{Z}\). For \(j \in \bb{Z}\), the category of \textit{j-effective motivic spectra} is defined as 
\[\sh(S)^\mathrm{eff}(j) = \bb{T}^{\otimes j} \otimes \sh(S)^\mathrm{eff},\]
so that \(\sh(S)^\mathrm{eff}(j + 1) \subset \sh(S)^\mathrm{eff}(j)\). The inclusions \(\iota_j \colon \sh(S)^\mathrm{eff}(j) \hookrightarrow \sh(S)\) preserve colimits, thus possess right adjoints, denoted \(r_j\). Defining \(\fil_\sli^j = \iota_j \circ r_j\), gives a filtration on any \(E \in \sh(S)\): 
\[\begin{tikzcd} \cdots \ar[r] & \fil_\sli^{j + 1} E \ar[r] & \fil_\sli^j E \ar[r] & \fil_\sli^{j - 1} E \ar[r] & \cdots \ar[r] & E \end{tikzcd}\]
This is called the \textit{slice filtration}, and its cofibres, 
\[s^j E = \mathrm{cofib}\big(\fil_\sli^{j + 1} E \to \fil_\sli^j E\big),\]
are called \emph{slices} of the motivic spectrum \(E\). For any \(j \in \bb{Z}\) and \(F \in \sh(S)^\mathrm{eff}(j)\), a morphism of motivic spectra \(F \to E\) factors through \(\fil_\sli^j(E) \to E\). Additionally, for \(F \in \sh(S)^\mathrm{eff}(j + 1)\), \(E \in \sh(S)\), the mapping spectrum \(\map_{\sh(S)}(F, s^jE)\) vanishes.
\end{definition}

Just as topological \(\K\)-theory is representable in the category of spectra, homotopy \(\K\)-theory is representable in the motivic homotopy category \cite{Cis}: 

\begin{theorem}[Cisinski]
Let \(X\) be a QcQs scheme. There is a motivic spectrum \(\Kgl_X \in \sh(X)\) representing homotopy \(K\)-theory: for all \(Y \in \sm_X\), \(n \in \bb{Z}\), 
\[\pi_n \! \left(\map_{\sh(X)}\big(M_X(Y), \Kgl_X\big)\right) \cong \KH_n(Y)\]
This is \emph{Bott periodic}: there is a canonical equivalence \(\bb{T} \otimes \Kgl_X \simeq \Kgl_X\).

\end{theorem}

\begin{remark}[Bott periodicity and slices]
Bott periodicity will allow us to deduce a similar periodic description of slices of \(\Kgl_X\). By \cite[Lemma 2.1]{RonOst}, for any \(\mathcal{E} \in \sh(X)\), there is a functorial equivalence 
\[\fil_\sli^j \big(\bb{T} \otimes \mathcal{E}\big) \simeq \bb{T} \otimes \fil_\sli^{j - 1} \, \mathcal{E}, \quad s^j \big(\bb{T} \otimes \mathcal{E}\big) \simeq \bb{T} \otimes s^{j - 1} \mathcal{E}\]
As a result, we arrive at convenient description of slices : for all \(j\), \(s^j \Kgl_X \simeq \bb{T} \otimes s^{j - 1} \Kgl_X\).
\end{remark}

Bachmann-Elmanto-Morrow use slices of \(\Kgl\) to extend the definition of (\(\bb{A}^1\)-invariant) motivic cohomology

\begin{definition} Given QcQs scheme \(X\), its \emph{\(\bb{A}^1\)-motivic cohomology} is defined by 
\[\bb{Z}(j)^\bb{A}(X) = \map_{\sh(X)}\big(\mathbbm{1}_X, s^j \Kgl_X\big)[- 2j] \in \Sp,\]
for \(j \in \bb{Z}\), and we write \(H_\bb{A}^i\big(X, \bb{Z}(j)\big) = H^i\big(\bb{Z}(j)^\bb{A}(X)\big)\) for the corresponding \emph{\(\bb{A}^1\)-motivic cohomology groups}.
\end{definition}

We may also take the slices over the terminal scheme \(\spec \bb{Z}\), then pull them back along the structure morphism \(f \colon X \to \spec \bb{Z}\), giving an ostensibly different mapping spectrum:

\begin{definition} The \emph{\(\bb{A}^1\)-cdh-motivic cohomology} of \(X\) is defined by 
\[\bb{Z}(j)^{\bb{A}, \mathrm{cdh}}(X) = \map_{\sh(X)}\big(\mathbbm{1}_X, f^* s^j \Kgl_\bb{Z}\big)[- 2j] \in \Sp,\]
for \(j \in \bb{Z}\), and we write \(H_{\bb{A}, \mathrm{cdh}}^i\big(X, \bb{Z}(j)\big) = H^i\big(\bb{Z}(j)^{\bb{A}, \mathrm{cdh}}(X)\big)\) for the corresponding \emph{\(\bb{A}^1\)-cdh-motivic cohomology groups}.
\end{definition}

An extremely nontrivial result of Bachmann-Elmanto-Morrow \cite[Theorem 9.3 (2)]{BEM} is that these complexes are equivalent.

\begin{theorem}[Equivalence of theories] The canonical comparison map \(\bb{Z}(j)^{\bb{A}, \mathrm{cdh}}(X) \to \bb{Z}(j)^\bb{A}(X)\) is an equivalence of presheaves of spectra. 
\end{theorem}

As a result, we will always abusively write \(\bb{Z}(j)^\bb{A}(X)\) in place of \(\bb{Z}(j)^{\bb{A}, \mathrm{cdh}}(X)\).

\subsection{GGL and Bivariant theories}

Just as homotopy K-theory is representable in \(\sh(X)\), \(\G\)-theory is also representable. Jin \cite[Corollary 2.2.16]{Jin}, proves the following: 

\begin{theorem}[Representability of G-theory] Let \(X\) be some Noetherian scheme. There exists a motivic spectrum \(\Ggl_X \in \sh(X)\) representing G-theory: for all \(Y \in \sm_X\), 
\[\pi_n \! \left(\map_{\sh(X)}\big(M_X(Y), \Ggl_X\big)\right) \cong G_n(Y)\]
Additionally, for finite type separated \(f \colon X \to S\), there is a canonical equivalence
\begin{equation}\label{eq:G-shriek}
f^! \Ggl_S \simeq \Ggl_X
\end{equation}
\end{theorem}

Consequently, there is an action of motivic spectra, 
\[\Kgl_X \otimes \Ggl_X \to \Ggl_X\]
induced by the projection formula.

\begin{remark}[Regular schemes] When \(X\) is a regular Noetherian scheme, Thomason and Trobaugh \cite{TT} prove that there is a canonical equivalence of K-theory and G-theory spaces: \(\K(X) \simeq \G(X)\). Thus there is a canonical equivalence in \(\sh(X)\)
\begin{equation}\label{eq:G-reg}
\Ggl_X \simeq \Kgl_X
\end{equation}
\end{remark}

\begin{definition}[Bivariant theories and fundamental classes, following \cite{DJK}]
Let \(S\) be a scheme, and \(E \in \sh(S)\) a motivic spectrum. Given morphism \(f \colon X \to S\), and K-theory class \(v \in \K(X)\), we define the \emph{\(v\)-twisted bivariant spectrum of \(X\) over \(S\)} as 
\[E(X/S, v) = \map_{\sh(X)}\big(\mathrm{Th}(v), f^!E\big)\]
We write \(E_n(X/S, v)\) for the \(n^\mathrm{th}\) homotopy group of this spectrum. 
This satisfies the following axioms: 
\begin{enumerate}

\item \emph{Base change:} Given Cartesian square 
\begin{equation}\label{eq:cart-sq}
\begin{tikzcd}
\tilde{X} \ar[d, "\tilde{f}", swap] \ar[r, "\tilde{h}", ""{name=h1}] & X \ar[d, "f"]\\
\tilde{S} \ar[r, "h", swap, ""{name=h2}] & S 
\ar[phantom, from=h1, to=h2, "\Delta"]
\end{tikzcd}
\end{equation}
there is a canonical map 
\[\Delta^* \colon E(X/S, v) \to E\big(\tilde{X}/\tilde{S}, \tilde{h}^*v\big),\]
induced by base change.

\item \emph{Proper covariance:} Given proper morphism \(g \colon X \to Y\) over \(S\), there is a pushforward 
\[g_* \colon E(X/S, g^*v) \to E(Y/S, v),\]
induced by the unit \(g_! g^! \Longrightarrow \id_{\sh(X)}\) and the equivalence \(g_! \simeq g_*\) by propriety of \(g\).

\item \emph{Étale contravariance} Given étale \(g \colon X \to Y\), there is a pullback 
\[g^! \colon E(Y/S, v) \to E(X/S, g^*v),\]
induced by the purity transformation \(\mf{p}_g \colon g^! \simeq g^*\).

\item \emph{Product}
If \(E \in \sh(S)\) has multiplication \(\mu \colon E \otimes E \to E\), then for morphisms \(g \colon X \to Y\), \(q \colon Y \to S\) and classes \(v \in \K(X)\), \(w \in \K(Y)\), there is a product
\[E_n(X/Y, v) \otimes E_m(Y/S, w) \to E_{m + n}(X/S, v + g^*w),\]
in which we ``lift'' the bivariant theory to \(X\) by 
\[E_n(X/Y, v) = \map_{\sh(X)}\big(\mathrm{Th}_X(v)[n], g^! q^*E)\]
The bivariant product sends maps \(x \colon \mathrm{Th}_X(v)[n] \to g^!q^*E\) and \(y \colon \mathrm{Th}_Y(w)[m] \to q^!E\) to \(x \cdot y \colon \mathrm{Th}_X(v + g^*w)[m + n] \to (q \circ g)^! E\), given as the composite 
\[\begin{tikzcd} [column sep = large]
\mathrm{Th}_X(v)[n] \otimes \mathrm{Th}_X(g^*w)[m] \ar[r, "x \otimes \id"] & g^!q^*E \otimes g^*\mathrm{Th}_Y(w)[m] \ar[r, "\pr_*^!"] & g^!\big(g^*E \otimes \mathrm{Th}_Y(w)[m]\big) \ar[ddll, "g^!(\id \otimes y)", swap]\\
\\
g^!\big(q^*E \otimes q^!E\big) \ar[r, "q^! \pr_*^!", swap] & g^! q^! \big(E \otimes E\big) \ar[r, "g^! q^! \mu", swap] & (q \circ g)^! E
\end{tikzcd}\]

\item Localization sequence: given closed embedding \(i \colon Z \hookrightarrow X\), with quasicompact open complement \(j \colon U \hookrightarrow X\), there is a canonical fibre sequence in \(\Sp\) 
\[\begin{tikzcd}
E(Z/S, i^*v) \ar[r, "i_*"] & E(X/S, v) \ar[r, "j^*"] & E(U/S, j^*v)
\end{tikzcd}\]
natural in the étale contravariance and  proper covariance over \(S\).
\end{enumerate} 
\end{definition}

\begin{theorem}[Déglise-Jin-Khan]\label{thm:fund-class}
When \(f \colon X \to S\) is an lci morphism, the cotangent complex \(\bb{L}_f\) is an element of the K-theory space, and there is a canonical element, the \emph{fundamental class of \(f\)},
\[\eta_f \in E_0(X/S, \bb{L}_f)\]
When \(f \colon X \to S\) is smooth, \(\eta_f\) is induced by the purity equivalence 
\[\p_f \colon \mathrm{Th}(T_f) \otimes f^* \to f^!\]
Likewise, when \(i \colon Z \hookrightarrow X\) is a regular closed embedding, \(\eta_i\) is induced by the purity transformation 
\[\p_i \colon \mathrm{Th}(- N_Z X) \otimes i^* \to i^!\]
This satisfies \emph{transverse base change}: given a tor-independent Cartesian square \(\Delta\) as in \eqref{eq:cart-sq}, there is an equivalence 
\begin{equation}\label{eq:tranv-b-ch}
\Delta^*\eta_f \simeq \eta_{\tilde{f}}
\end{equation}
\end{theorem}

\begin{construction}[Gysin morphism]
Gysin map: any lci morphism \(p \colon X \to Y\) over \(S\) induces a \emph{Gysin map}, 
\[p^! \colon E(Y/S, v) \to E(X/S, p^*v + \bb{L}_f ), \quad x \mapsto \eta_p \cdot x\]
\end{construction}

We'll also use a similar construction with more flexible weights:

\begin{definition}[Borel-Moore homology groups]\label{def:BM}
Let \(S\) be a scheme, and \(E \in \sh(S)\) a motivic spectrum. Given \(f \colon X \to S\), and \(a, b \in \bb{Z}\), we have the \emph{Borel-Moore homology groups}, 
\[E_{a, b}^\bm(X/S) = \pi_0 \! \left(\map_{\sh(X)}\left(\bb{T}^{\otimes b}[a - 2b], f^!E\right)\right)\]
\end{definition}

\begin{remark}[Bivariant vs Borel-Moore homology groups]
Borel-Moore homology groups have all the same functorialities, products, and localization properties as the bivariant theories, and given lci \(f \colon X \to S\) of relative virtual rank \(d\) and trivializable virtual tangent bundle, then we have the identification
\[E_{2d, d}^\bm(X/S) \cong E(X/S, \bb{L}_f)\]
\end{remark}

As \(\G(X)\) is the homological counterpart of \(\KH(X)\), we would expect \(\Ggl_X\) to be the Borel-Moore counterpart of \(\Kgl_X\). Jin \cite[Corollarly 3.3.7]{Jin} proves that it is:

\begin{theorem}[Jin]
For any separated finite type \(X\) over \(S\), there is a natural isomorphism 
\[\Kgl_{m, n}^\bm(X/S) \cong G_{m - 2n}(X),\]
such that proper covariance for \(\Kgl^\bm\) is identified with proper functoriality of \(\G\)-theory, and the Gysin morphisms on \(\Kgl^\bm\) are identified with lci functoriality of \(\G\)-theory.
\end{theorem}

The main reason for introducing this closely related construction is the following computational tool: 

\begin{construction}[\(\delta\)-niveau spectral sequence]
Let \(X\) be finite type separated, with dimension function \(\delta\). Denote by \(X_{(p)} = \{x \in X \, | \, \delta(x) = p\}\), and \(Z(x)\) the reduced closure of \(x\) in \(X\). Given a motivic spectrum \(E \in \sh(\bb{Z})\), write 
\[\widehat{E}_{a + b, i}^\bm(x/\bb{Z}) = \colim_{V \subset Z(x)} E_{a + b, i}^\bm(V/\bb{Z}),\]
for \(V\) any nonempty open subscheme. Bondarko and Déglise \cite[Definition 3.1.5]{BonDeg} construct, for any \(i \in \bb{Z}\), a spectral sequence 
\begin{equation}\label{eq:dnss}
E_{a, b}^1 = \bigoplus_{x \in X_{(a)}} \widehat{E}_{a + b, i}^\bm(x/\bb{Z}) \Longrightarrow E_{a + b, i}^\bm(X/\bb{Z}),
\end{equation}
which converges for \(X\) Noetherian and of finite dimension.
\end{construction}

\subsection{The Bloch-Levine cycle complex}
In this section, we review Levine's generalization \cite{LevTechLoc} of Bloch's cycle complex.

\begin{definition}[Algebraic simplices] Let \(q \geq 0\). The \emph{algebraic \(q\)-simplex} is the scheme
\[\Delta^q = \spec \! \left(\bb{Z}[t_0, \dots , t_q] \big/ \left(\sum_i t_i - 1\right)\right),\]
with \emph{faces} the closed subschemes \(\Delta^{q - m} \cong V(t_{i_1}, \dots , t_{i_m})\). The closed embeddings \(\Delta^{q - 1} \hookrightarrow \Delta^q\) give coface maps, with analogous codegeneracy maps \(\Delta^q \to \Delta^{q - 1}\) induced by projection making \(\Delta^\bullet\) a cosimplicial scheme.
\end{definition}

\begin{construction}[Bloch-Levine cycle complex]\label{def:LCC}
Given finite type irreducible \(f \colon X \to \spec \bb{Z}\), let \(\eta\) be the generic point of \(X\), and \(\p\) its image under \(f\). The \emph{relative dimension of \(X\) over \(\bb{Z}\)} defined as 
\[\dim_\bb{Z} X = \trd \big(\kappa(\eta) : \kappa(\p)\big) - \mathrm{ht} \, \p + 1\]
Given \(q, n \in \bb{Z}\), and a scheme \(X\), define \(z_n(X, q)\) to be the free Abelian group on integral closed subschemes \(W\) of \(X \times \Delta^q\) of relative dimension \(n + q\) over \(\bb{Z}\) intersecting all faces properly, that is, 
\[\dim_\bb{Z} (W \cap \Delta^r) \leq n + r\]

This produces a simplicial Abelian group \(z_n(X, \bullet)\). Under the Dold-Kan correspondence, there is an associated non-negatively graded complex of Abelian groups, \(z_n(X, *)\). This upgrades to a complex of presheaves on \(X\), the \emph{Bloch-Levine cycle complex}, denoted by \(\mathcal{Z}_n(X, *)\), with rule \(\mathcal{Z}_n(X, *)(U) = z_n(U, *)\).
\end{construction}

\begin{definition} The \emph{higher Chow groups}, also called the \emph{motivic Borel-Moore homology groups of \(X\)}, are the Zariski hyperhomology groups 
\[\mathrm{CH}_n(X, j) = \bb{H}^{- j}\big(X_\mathrm{Zar}, \mathcal{Z}_n(X, *)\big)\]
\end{definition}

\begin{remark}[Functoriality]
Fixing \(q \geq 0\), a flat morphism \(f \colon X \to Y\) of relative dimension \(d\) induces a morphism \(f^* \colon z_n(Y, q) \to z_{n + d}(X, q)\), given by the usual flat pullback of Fulton \cite[1.7]{Ful}. Likewise, a proper \(p \colon C \to X\) induces a pushforward \(p_* \colon z_n(C, q) \to z_n(X, q)\), again defined as in Fulton \cite[1.4]{Ful}. These upgrade to the level of complexes, giving morphisms: \(f^* \colon \mathcal{Z}_n(Y, *) \to \mathcal{Z}_{n + d}(X, *)\),  \(p_* \colon \mathcal{Z}_n(C, *) \to \mathcal{Z}_n(X, *)\).
\end{remark}

\begin{remark}[Shifted dimension] We will use a different definition of dimension, which is simply one less than the relative dimension over \(\bb{Z}\) used to define \(\mathcal{Z}_n(X, *)\). Given \(X, \eta\), and \(\p\) as before, define 
\[\delta(X) = \trd\big(\kappa(\eta) : \kappa(\p)\big) - \mathrm{ht} \, \p\]
The reason for such a change is motivated by the following fact: if \(f \colon X \to \spec \bb{Z}\) factors as 
\[\begin{tikzcd} X \ar[rr, "f", swap, bend right = 19] \ar[r, "g"] & \spec \bb{F}_p \ar[r, "i", hookrightarrow] & \spec \bb{Z} \end{tikzcd}\]
with \(g\) smooth, then the virtual rank of the cotangent complex of \(f\) equals \(\delta(X)\). 
\end{remark}

\section{Filtration on G-Theory}\label{G-fil}

In this section, we construct a filtration on G-theory, and prove some convergence properties of the associated spectral sequence. We also demonstrate that the action of \(\Kgl\) on \(\Ggl\) is compatible with the respective associated graded pieces of the filtrations.

\subsection{The filtration and convergence properties}

The desire to obtain a Borel-Moore counterpart to the \(\bb{Z}(j)^\bb{A}(X)\) motivates an analogous construction, but with exceptional pullback instead of pullback.

\begin{definition}
Given finite type separated \(f \colon X \to \spec \bb{Z}\), its \emph{motivic Borel-Moore homology} is the spectrum (in fact, complex) 
\[\bb{Z}(n)^\bm(X) = \map_{\sh(X)}\big(\mathbbm{1}_X, f^! s^{- n} \Kgl_\bb{Z}\big)\lbrack 2n \rbrack, \quad n \in \bb{Z}\]
Write \(H_j^\bm\big(X, \bb{Z}(n)\big) = H_j\big(\bb{Z}(n)^\bm(X)\big)\) for the corresponding homology groups.\footnote{Hearkening back to \ref{def:BM}, this gives \(H_j^\bm\big(X, \bb{Z}(n)\big) \cong H\bb{Z}_{n, j}^\bm(X/\bb{Z})\), where \(H\bb{Z} = s^0 \Kgl_\bb{Z} \in \sh(\bb{Z})\). So we may think of \(\bb{Z}(n)^\bm(X)\) as the Borel-Moore groups of motivic cohomology, whence the title of this paper).}
\end{definition}

\begin{remark}[Borel-Moore functoriality of \(\bb{Z}(n)^\bm\)]
Given a smooth morphism \(g \colon Y \to X\) of relative dimension \(d\), applying the unit of \(g^* \dashv g_*\) to the target and smooth purity gives a morphism of spectra 
\[\bb{Z}(n)^\bm(X) \to \map_{\sh(X)}\big(\mathbbm{1}_X, g_*g^*f^! s^{- n} \Kgl_\bb{Z}\big)\lbrack 2n \rbrack \simeq \map_{\sh(Y)}\big(\mathbbm{1}_Y, \bb{T}^{\otimes - d} \otimes g^* f^!s^{- n}\Kgl_\bb{Z}\big) \lbrack 2n \rbrack,\]
thus giving smooth Borel-Moore functoriality,
\[g^* \colon \bb{Z}(n)^\bm(X) \to \bb{Z}(n + d)^\bm(Y) \lbrack -2d \rbrack\]
If \(p \colon X \to C\) is proper, then writing \(h \colon C \to \spec \bb{Z}\) as structure morphism applying the counit of \(p_! \dashv p^!\), along with \(p_! \simeq p_*\) by propriety of \(p\), gives a morphism of spectra 
\[\map_{\sh(X)}\big(p^*\mathbbm{1}_C, p^!h^! s^{- n} \Kgl_\bb{Z}\big)\lbrack 2n \rbrack \simeq \map_{\sh(C)}\big(\mathbbm{1}_C, p_*p^!h^! s^{- n} \Kgl_\bb{Z}\big)\lbrack 2n \rbrack \to \map_{\sh(C)}\big(\mathbbm{1}_C, h^! s^{- n}\Kgl_\bb{Z}\big)\lbrack 2n \rbrack,\]
thus giving proper Borel-Moore functoriality, 
\[p_* \colon \bb{Z}(n)^\bm(X) \to \bb{Z}(n)^\bm(C)\]
These are exactly the functorialities enjoyed by \(\mathcal{Z}_{n + 1}(X, *)\).
\end{remark}

\begin{lemma}[Dimension vanishing]\label{lemm:dim-van}
Write \(D = \dim X\) for the Krull dimension. Then \(\bb{Z}(n)^\bm(X) \simeq 0\) for \(n \geq D\).
\end{lemma}

\begin{proof}
In degree \(q\), \(\mathcal{Z}_{n + 1}(X, *)\) is the complex \(z_{n + 1}(U, q)\) of cycles of \(\delta = n + q + 1\) on \(U \times \Delta^q\), which is of Krull dimension at most \(D + q\), so \(n \geq D\) means there are no such cycles, and \ref{main-intro} (5), which we prove later, tells us that \(\bb{Z}(n)^\bm(X) \lbrack -2n \rbrack \simeq 0\), so its shifting by \(2n\) gives the first statement.\end{proof}

The following lemma helps us establish convergence properties of \(\fil_\star^\bm \G(X)\).

\begin{lemma}\label{lemm:comp:KGL}
The slice filtration on \(\Kgl_\bb{Z}\) is complete.
\end{lemma}

\begin{proof}
\(\sh(\bb{Z})\) is compactly generated by the objects \(M_\bb{Z}(S) \otimes \bb{T}^{\otimes a}\), for \(a \in \bb{Z}\) \cite[Proposition C.12 (1)]{Hoy}. Writing \(f_S \colon S \to \spec \bb{Z}\) for the structure morphism, we have that 
\[\map_{\sh(\bb{Z})}\big(M_\bb{Z}(S) \otimes \bb{T}^{\otimes a}, \fil_\sli^j\Kgl_\bb{Z}\big) \simeq \map_{\sh(\bb{Z})}\big((f_S)_\sharp \mathbbm{1}_S, \bb{T}^{\otimes - a} \otimes \fil_\sli^j\Kgl_\bb{Z}\big)\]
\[\simeq \map_{\sh(S)}\big(\mathbbm{1}_S, f_S^*\fil_\sli^{j - a} \Kgl_\bb{Z}\big) \simeq \map_{\sh(S)}\big(\mathbbm{1}_S, \fil_\sli^{j - a} \Kgl_S\big) = \fil_\bb{A}^{j - a} \KH(S) \]
\(S\) is Noetherian, and of finite type, so has finite Krull dimension which is equal to its valuative dimension by \cite[Proposition 2.3.2 (9)]{EHMK}. The filtration \(\fil_\bb{A}^\star \KH(S)\) is complete for finite valuative dimension schemes \cite[Corollary 9.5]{BEM}, thus 
\[\catlim_{j \to \infty} \map_{\sh(\bb{Z})}\big(M_\bb{Z}(S) \otimes \bb{T}^{\otimes a}, \fil_\sli^j \Kgl_\bb{Z}\big) \simeq \catlim_{j \to \infty} \KH_\bb{A}^j(S) \simeq 0\]
Compact generation thus gives that 
\[\catlim_{j \to \infty} \fil_\sli^j \Kgl_\bb{Z} \simeq 0,\]
so we are done.
\end{proof}

\begin{theorem}
The spectrum \(\bb{Z}(n)^\bm(X) \in \Sp\) has the following properties: 
\begin{enumerate}

\item There is a functorial decreasing filtration \(\fil_\star^\bm \G(X)\) on \(\G(X)\) such that the associated gradeds are given by 
\[\gr_n^\bm \G(X) \simeq \bb{Z}(n)^\bm(X)\lbrack -2n \rbrack\]

\item This filtration is exhaustive and complete, hence inducing a natural Atiyah-Hirzebruch-style spectral sequence 
\[E_{p, q}^2 = H_{p - q}^\bm\big(X, \bb{Z}(- q)\big) \Longrightarrow \G_{p + q}(X)\]

\item Letting \(D\) denote the Krull dimension of \(X\) again, \(\fil_n^\bm \G(X) \to \G(X)\) is an equivalence for \(n \geq D - 1\). Writing \(N(j) = - j - D - 3\), \(\pi_j \fil_n^\bm \G(X) \cong 0\) for \(n \leq N(j)\). Consequently, the spectral sequence in (2) converges strongly.

\item The presheaf \(\bb{Z}(n)^\bm\) converts elementary distinguished squares of schemes to Cartesian squares of spectra.

\item \(\bb{Z}(n)^\bm\) is \(\bb{A}^1\)-invariant up to a twist: 
\[\bb{Z}(n)^\bm(\bb{A}_X^1) \simeq \bb{Z}(n - 1)^\bm(X)[2]\]

\end{enumerate}
\end{theorem}

\begin{proof}
Let \(f \colon X \to \spec \bb{Z}\) be finite type separated. On \(\bb{Z}\), by \eqref{eq:G-reg} \(\Ggl_\bb{Z} \simeq \Kgl_\bb{Z}\). Equipping \(\Ggl_\bb{Z}\) with the slice filtration, and applying the functor \(f^!\), we get a filtration in \(\sh(X)\), 
\[\begin{tikzcd} \cdots \ar[r] & f^!\fil_\sli^{j + 1} \Kgl_\bb{Z} \ar[r] & f^!\fil_{\sli}^j \Kgl_\bb{Z} \ar[r] & \cdots \ar[r] & f^!\Kgl_\bb{Z} \stackrel{\eqref{eq:G-shriek}} \simeq \mathrm{GGL}_X \end{tikzcd}\]
Define 
\[\fil_\star^\bm \G(X) = \map_{\sh(X)}\big(\mathbbm{1}_X, f^! \fil_\sli^\star \Kgl_\bb{Z}\big)\]
Taking cofibres commutes with the exact functor \(f^!\), so the associated gradeds of the filtration are shifts of \(\bb{Z}(n)^\bm(X)\) by construction, whence (1).

To prove exhaustiveness, note that by \cite[Definition 1.3.6]{BonDeg}, \(f_! \mathbbm{1}_X\) is compact in \(\sh(\bb{Z})\).\footnote{This result actually concerns compactness in the homotopy category of \(\sh(\bb{Z})\), but this is equivalent to compactness in the underlying \(\infty\)-category as well by \cite[Remark 1.4.4.3]{HA}} Adjointness gives us that 
\[\fil_\infty^\bm \G(X) = \colim_{j \to - \infty} \map_{\sh(X)}\big(\mathbbm{1}_X, f^! \fil_\sli^j \Kgl_\bb{Z}\big) \simeq \colim_{j \to - \infty} \map_{\sh(\bb{Z})}\big(f_! \mathbbm{1}_X, \fil_\sli^j \Kgl_\bb{Z} \big),\]
so compactness of \(f_! \mathbbm{1}_X\) and exhaustiveness of the slice filtration \cite[Remark 3.38]{BEM} equate this to 
\[\map_{\sh(\bb{Z})} \! \left(f_! \mathbbm{1}_X, \colim_{j \to - \infty} \fil_\sli^j \Kgl_\bb{Z}\right) \simeq \map_{\sh(\bb{Z})}\big(f_! \mathbbm{1}_X, \Kgl_\bb{Z}\big) \simeq \G(X)\]
Similarly, 
\[\catlim_{j \to \infty} \fil_j^\bm \G(X) = \catlim_{j \to \infty} \map_{\sh(X)}\big(\mathbbm{1}_X, \fil_\sli^j \Kgl_\bb{Z}\big) \simeq \map_{\sh(X)} \!\left(\mathbbm{1}_X, \catlim_{j \to \infty} \fil_\sli^j \Kgl_\bb{Z}\right) \simeq 0,\]
giving completeness.

Dimension vanishing \ref{lemm:dim-van} tells us that \(\fil_{n - 1}^\bm \G(X) \to \fil_n^\bm \G(X)\) is an equivalence for \(n \geq D\), thus \(\fil_n^\bm \G(X)\) is constant for \(n \geq D - 1\). Exhaustiveness identifies the constant value with \(\G(X)\).

For the other vanishing statement, once again compare \(\bb{Z}(n)^\bm\) to \(\mathcal{Z}_{n + 1}(X, *)\) as in the proof of dimension vanishing. When \(n + 1 + q < 0\), then there are no cycles of this relative dimension in degree \(q\), giving a complex concentrated in homological degree \(q \geq - (n + 1)\), thus cohomological degree \(n + 1\). \(X\) is Noetherian, so its Zariski site has cohomological dimension at most \(D\) by Grothendieck vanishing. Hence \(R\Gamma_\mathrm{Zar}\) is of cohomological degree \(\leq n + 1 + D\), so that \(\pi_j = H^{- j}\) vanishes for \(- j > n + 1 + D\). Thus 
\[\pi_j \gr_n^\bm \G(X) \cong \pi_{j + 1}^\bm \gr_n^\bm \G(X) \cong 0 \quad \text{ for } n \leq N(j)\]
Combined with the cofibre sequence 
\[\begin{tikzcd}
\fil_{n - 1}^\bm \G(X) \ar[r] & \fil_n^\bm \G(X) \ar[r] & \gr_n^\bm \G(X),
\end{tikzcd}\]
this implies that \(\pi_j \fil_{n - 1}^\bm \G(X) \cong \pi_j \fil_n^\bm \G(X)\) when \(n \leq N(j)\), so that \(\big(\pi_j \fil_n^\bm \G(X)\big)_j\) is constant for \(j \leq N(j)\). Taking \(n \to - \infty\), we obtain the Milnor sequence 
\[\begin{tikzcd}
0 \ar[r] & \underset{n \to - \infty}{\catlim^1} \pi_{j + 1} \fil_n^\bm \G(X) \ar[r] & \displaystyle \pi_j \! \left(\catlim_{n \to - \infty} \fil_n^\bm \G(X)\right) \ar[r] & \displaystyle \catlim_{n \to - \infty} \pi_j \fil_n^\bm \G(X) \ar[r] & 0
\end{tikzcd}\]
The middle term vanishes by completeness, meaning the right term vanishes as well, so we see that \(\pi_j \fil_n^\bm \G(X) \cong 0\) for \(n \leq N(j)\).

For (4), let 
\[\begin{tikzcd} [row sep = large, column sep = large]
U \times_X V \ar[d] \ar[r, "\tilde{j}"] & V \ar[d, "p"]\\
U \ar[r, "j", swap, hookrightarrow] & X
\end{tikzcd}\]
be a Nisnevich distinguished square, in which étale \(p\) is an isomorphism on \(Z = X - U\). Writing \(i \colon Z \hookrightarrow X\) for the complementary closed embedding, application of the mapping spectrum \(\map_{\sh(X)}(\mathbbm{1}_X, -)\) to the localization sequence 
\begin{equation}\label{eq:nis-ex}
\begin{tikzcd} 
i_*i^! f^! s^{- n} \Kgl_\bb{Z} \ar[r] & f^! s^{- n} \Kgl_\bb{Z} \ar[r] & j_*j^! f^! s^{- n} \Kgl_\bb{Z},
\end{tikzcd}
\end{equation}
adjointness, and shifting by \(2n\) means that 
\[\begin{tikzcd}
\bb{Z}(n)^\bm(Z) \ar[r] & \bb{Z}(n)^\bm(X) \ar[r] & \bb{Z}(n)^\bm(U)
\end{tikzcd}\]
is a fibre sequence. As an elementary distinguished square, \(p|_Z \colon Z' = p|_Z^{- 1}(Z) \to Z\) is an isomorphism of schemes. This gives the pullback square 
\[\begin{tikzcd} [row sep = large, column sep = large]
Z' \ar[d, "p|_Z", swap] \ar[r, "\tilde{i}"] & V \ar[d, "p"]\\
Z \ar[r, "i", swap] & X
\end{tikzcd}\]

Applying \(p^! \simeq p^*\) to the sequence \ref{eq:nis-ex} gives 
\[\begin{tikzcd} p^! i_* i^! f^! s^{- n} \Kgl_\bb{Z} \ar[r] & p^! f^!\Kgl_\bb{Z} \ar[r] & p^*j_*j^! f^! s^{- n} \Kgl_\bb{Z}
\end{tikzcd}\]
By base change, \(p^! i_* \simeq \tilde{i}_* p|_Z^!\) and \(p^*j_* \simeq \tilde{j}_*\tilde{p}^*\), so applying the exact functor \(\map_{\sh(V)}(\mathbbm{1}_V, -)\), using adjointness, and shifting by \(2n\) gives a fibre sequence 
\[\begin{tikzcd}
\bb{Z}(n)^\bm(Z') \ar[r] & \bb{Z}(n)^\bm(V) \ar[r] & \bb{Z}(n)^\bm(U \times_X V)
\end{tikzcd}\]
So we have a commutative diagram of fibre sequences 
\[\begin{tikzcd} [row sep = large, column sep = large]
\bb{Z}(n)^\bm(Z) \ar[d, "p^*", swap] \ar[r] & \bb{Z}(n)^\bm(X) \ar[d, "p^*"] \ar[r] & \bb{Z}(n)^\bm(U) \ar[d, "p^*"]\\
\bb{Z}(n)^\bm(Z') \ar[r] & \bb{Z}(n)^\bm(V) \ar[r] & \bb{Z}(n)^\bm(U \times_X V)
\end{tikzcd}\]
Since \(p\) is an isomorphism from \(Z\) to \(Z'\), the left vertical leg is an equivalence, \(\bb{Z}(n)^\bm(Z) \simeq \bb{Z}(n)^\bm(Z')\), thus the rightmost square is Cartesian.

The final statement follows by factoring \(f_{\bb{A}_X^1} \colon \bb{A}_X^1 \to \spec \bb{Z}\) as \(f \circ \pr_1\), for \(\pr_1 \colon \bb{A}_X^1 \to X\), which is smooth of relative dimension 1. Smooth purity and adjointness thus give equivalences 
\[\bb{Z}(n)^\bm\big(\bb{A}_X^1\big) \simeq \map_{\sh(\bb{A}_X^1)} \! \left(\mathbbm{1}_{\bb{A}_X^1}, \bb{T} \otimes \pr_1^* f^! s^{- n} \Kgl_\bb{Z}\right)\!\lbrack 2n \rbrack \simeq \map_{\sh(X)} \! \left((\pr_1)_\sharp \mathbbm{1}_{\bb{A}_X^1}, f_X^! s^{- n + 1} \Kgl_\bb{Z}\right)\!\lbrack 2n \rbrack,\]
and \(\bb{A}^1\)-invariance of \(\pr_1\) means this is equivalent to 
\[\map_{\sh(X)}\big(\mathbbm{1}_X, f^! s^{- (n - 1)} \Kgl_\bb{Z}\big)\lbrack 2n \rbrack = \bb{Z}(n - 1)^\bm(X) \lbrack 2 \rbrack\]
\end{proof}

\subsection{Action by K-Theory}\label{K:action}
We first prove that the action by \(\Kgl\) on \(\Ggl\) described by Jin descends to the level of slices, giving motivic Borel-Moore homology the structure of a complex acted upon by \(\bb{A}^1\)-motivic cohomology.

\begin{theorem}[Slice actions and linearity]
Let \(X\) be a finite type separated scheme.
\begin{enumerate}

\item \(\fil_\star^\bm \G(X)\) is a filtered module over \(\fil_\bb{A}^\star \KH(X)\), with structure 
\[\fil_\bb{A}^j \KH(X) \otimes \fil_n^\bm \G(X) \to \fil_{n - j}^\bm \G(X)\]
Taking colimits over \(j, n\) recovers the action of \(\KH(X)\) on \(\G(X)\).

\item \(\bb{Z}(- \star)^\bm(X)\) is a graded module over the graded \(\bb{E}_\infty\)-ring \(\bb{Z}(\star)^\bb{A}(X)\), with action 
\[\bb{Z}(j)^\bb{A}(X) \otimes \bb{Z}(n)^\bm(X) \to \bb{Z}(n - j)^\bm(X)\]

\item The unique \(\bb{E}_\infty\) map \(\bb{Z} \to \bb{Z}(0)^\bb{A}(X)\) gives \(\bb{Z}(- \star)^\bm(X)\) the structure of a graded module in \(\mathrm{Gr} \, \mathrm{D}(\bb{Z})\), and the morphism above may be viewed in that category.

\end{enumerate}
\end{theorem}

\begin{proof}
Formation of the slice filtration is a lax monoidal functor \(\fil_\sli^\star \colon \sh(\bb{Z}) \to \fil \, \sh(\bb{Z})\), with the target endowed with Day convolution as its symmetric monoidal structure \cite[\S 13.4]{BH}. This gives \(\fil_\sli^\star \Kgl_\bb{Z}\) the structure of a filtered module over itself. Closure of \(\sh(\bb{Z})^\mathrm{eff}\) under smash product implies that \(\sh(\bb{Z})^\mathrm{eff}(p) \otimes \sh(\bb{Z})^\mathrm{eff}(q) \subset \sh(\bb{Z})^\mathrm{eff}(p + q)\).
\[\fil_\sli^p \Kgl_\bb{Z} \otimes \fil_\sli^q \Kgl_\bb{Z} \in \sh(\bb{Z})^\mathrm{eff}(p + q)\]
Application of the universal property of \(\fil_\sli^{p + q}(\Kgl_\bb{Z} \otimes \Kgl_\bb{Z})\) to the canonical morphism 
\[\fil_\sli^p \Kgl_\bb{Z} \otimes \fil_\sli^q \Kgl_\bb{Z} \to \Kgl_\bb{Z} \otimes \Kgl_\bb{Z}\]
gives a unique morphism 
\[\fil_\sli^p \Kgl_\bb{Z} \otimes \fil_\sli^q \Kgl_\bb{Z} \to \fil_\sli^{p + q}(\Kgl_\bb{Z} \otimes \Kgl_\bb{Z})\]
Application of the exact functor \(f^!\) gives a morphism 
\[f^!\big(\fil_\sli^p \Kgl_\bb{Z} \otimes \fil_\sli^q \Kgl_\bb{Z}\big) \to f^! \fil_\sli^{p + q}\big(\Kgl_\bb{Z} \otimes \Kgl_\bb{Z}\big),\]
and the projection formula gives a morphism 
\[f^* \fil_\sli^p \Kgl_\bb{Z} \otimes f^! \fil_\sli^q \Kgl_\bb{Z} \to f^! \big(\fil_\sli^p \Kgl_\bb{Z} \otimes \fil_\sli^q \Kgl_\bb{Z}\big),\]
so composing these, along with, 
\[f^! \fil_\sli^{p + q}\mu \colon f^! \fil_\sli^{p + q} \big(\Kgl_\bb{Z} \otimes \Kgl_\bb{Z}\big) \to f^! \fil_\sli^{p + q}\Kgl_\bb{Z},\]
for \(\mu\) the multiplication structure on \(\Kgl_\bb{Z}\), gives the desired morphism for \(p, - q\). The projection formula exhibits \(f^! \colon \sh(\bb{Z}) \to \sh(X)\) as a lax linear functor when we give \(\sh(X)\) the structure of an \(\sh(\bb{Z})\)-module via the symmetric monoidal functor \(f^*\), proving (1). Exhaustiveness of the two filtrations gives the claim regarding the colimit.

Taking associated gradeds is a symmetric monoidal functor \(\gr \colon \fil \, \Sp \to \gr \, \Sp\), so by the construction of the filtrations, we obtain the action (2).

Write \(R^\star = \bb{Z}(\star)^\bb{A}(X)\). By \cite[Corollary 9.12 (1)]{BEM}, \(R^0 = \bb{Z}(0)^\bb{A}(X) \cong R\Gamma_\mathrm{cdh}(X, \bb{Z})\), so is a coconnective spectrum. \(H\bb{Z}\) is connective, meaning there is a unique \(\bb{E}_\infty\) map \(H\bb{Z} \to R^0\), since it factors through \(\tau_{\geq 0} R^0\), which is discrete. Letting \(H\bb{Z}\langle 0 \rangle\) be the graded \(\bb{E}_\infty\)-ring with \(H\bb{Z}\) placed in degree 0, we have an adjoint equivalence 
\[\Map_{\mathrm{CAlg}(\mathrm{Gr} \, \Sp)} (H\bb{Z}\langle 0 \rangle, R^\star) \simeq \Map_{\mathrm{CAlg}(\Sp)}(H\bb{Z}, R^0),\]
The former is thus a contractible mapping space, giving \(\bb{Z}(\star)^\bb{A}(X)\) a unique \(\bb{E}_\infty\)-map from \(\mathrm{CAlg}\big(\mathrm{Gr} \, D(\bb{Z})\big)\), which in turn gives the graded module \(\bb{Z}(- \star)^\bm(X)\) a unique graded \(\bb{Z}(\star)^\bb{A}(X)\)-module structure in \(\mathrm{Gr} \, D(\bb{Z})\).
\end{proof}

\section{Slice Computation A}\label{slice-comp-A}
In this section, we begin the proof of the main result, which is that \(\bb{Z}(n)^\bm(X)\) computes the higher Chow groups. We construct an auxiliary object, the \emph{spectral replacement} of \(\mathcal{Z}_{n + 1}(X, *)\), and demonstrate that they have the same complex of derived global sections.

\subsection{Spectral replacement}

We first establish several lemmata needed to justify the formation of colimits of spectra.

\begin{definition}[Sections with support]
Let \(i \colon C \hookrightarrow S\) be a closed embedding, with complementary open embedding \(j \colon U \hookrightarrow S\). Write \(\sh_C(S)\) to be the full subcategory of \(\sh(S)\) on motivic spectra \(E\) with \(j^*E \simeq 0\). Taking \emph{sections with support} is the functor, 
\[\Gamma_C = i_*i^! \colon \sh(S) \to \sh(S),\]
and, \(\epsilon^i\), the counit of \(i_* \simeq i_! \dashv i^!\), gives a natural transformation from \(\Gamma_C\) to the identity, which we call \emph{forgetting support on \(C\).}
\end{definition}

\begin{lemma}
The functor \(i_* \colon \sh(C) \to \sh(S)\) is fully faithful with essential image \(\sh_C(S)\), and \(\Gamma_C\) is right adjoint to the inclusion \(\sh_C(S) \hookrightarrow \sh(S)\), with counit \(\epsilon^i\).
\end{lemma}

\begin{proof}
\(i_*\) being fully faithful is encoded in the six-functor formalism, as \(i\) is closed. Six functors also gives the localization sequence  
\[\begin{tikzcd}
j_!j^*E \ar[r] & E \ar[r] & i_*i^*E,
\end{tikzcd}\]
which proves that \(E \simeq i_*i^*E\), whence the essential image result, and adjointness gives, for any \(E \in \sh_C(S)\), \(F \in \sh(S)\), equivalences of mapping spectra 
\[\map_{\sh(S)}(E, F) \simeq \map_{\sh(S)}(i_*i^*E, F) \simeq \map_{\sh(C)}(i_*E, i^!F) \simeq \map_{\sh_C(S)}(E, \Gamma_C F),\]
the first equivalence by full faithfulness of \(i_*\), and the third equivalence by the fact that \(\sh_C(F)\) is a full subcategory of \(\sh(S)\) containing the object \(i^!F\). This identifies \(\Gamma_C\) as the right adjoint.
\end{proof}

We will use the following lemma several times to factor functors:

\begin{lemma}\label{lemm:fun-fact}
Let \(\mathcal{C} \subset \mathcal{D}\) be a full subcategory, with inclusion \(\iota \colon \mathcal{C} \hookrightarrow \mathcal{D}\) possessing a right adjoint, \(R\) and counit \(\epsilon \colon \iota R \Longrightarrow \id\). If \(F, G \colon \mathcal{E} \to \mathcal{D}\) are functors with values in \(\mathcal{C}\), then composition with \(\epsilon G\) gives equivalence 
\[\begin{tikzcd}
\mathrm{Map}_{\fun(\mathcal{E}, \mathcal{D})}(F, \iota RG) \ar[r, "\sim"] & \mathrm{Map}_{\fun(\mathcal{E}, \mathcal{D})}(F, G)
\end{tikzcd}\]
Thus each natural transformation \(F \Longrightarrow G\) factors through \(\epsilon G\), unique up to contractible choice, meaning \(G \Longrightarrow \iota RH\) are equivalent if they agree after composition by \(\epsilon H\).
\end{lemma}

\begin{proof}
\(R \circ - \) is right adjoint to \(\iota \circ -\) on functor categories, and the counit of this adjuction has value \(\epsilon G\) at \(G\). The adjoint equivalence is expressed in terms of the counit as 
\[\Map_{\fun(\mathcal{E}, \mathcal{C})}(F', RG) \simeq \Map_{\fun(\mathcal{E}, \mathcal{D})}(\iota F', G), \quad \alpha \mapsto \epsilon G \circ \iota \alpha,\]
for any functor \(F' \colon \mathcal{E} \to \mathcal{C}\). \(\iota\) is fully faithful, thus composition by \(\iota\) is an equivalence in the functor category, fitting into the commutative diagram 
\[\begin{tikzcd} [row sep = large, column sep = large]
\Map_{\fun(\mathcal{E}, \mathcal{C})}(F', RG) \ar[dr, "\iota", swap] \ar[r, "\epsilon G \circ \iota (- )"] & \Map_{\fun(\mathcal{E}, \mathcal{D})}(\iota F', \iota RG) \ar[d, "\epsilon G \circ - "]\\
& \Map_{\fun(\mathcal{E}, \mathcal{D})}(\iota F', G)
\end{tikzcd}\]
As \(\iota\) is fully faithful, \(\eta\), the unit of the adjunction \(\iota \dashv R\), is an equivalence, so that \(\iota R F \simeq F\), for any functor \(F \colon \mathcal{E} \to \mathcal{D}\) with values in \(\mathcal{C}\). Taking \(F'\) in the diagram above to be \(RF\), \(\iota F' \simeq F\), so we arrive at the desired result.
\end{proof}

\begin{lemma}[Formation of sections with support is functorial]\label{lemm:supp-fun}
Denoting by \(\mathcal{P}\) the poset category of closed subsets of \(S\), then 
\[\mathcal{P} \to \fun\big(\sh(S), \, \sh(S)\big)_{/\id}, \quad C \mapsto \epsilon^i \colon\Gamma_C \Longrightarrow \id,\]
is a functor, unique up to contractible choice. Given inclusion \(C_1 \subset C_2\), the natural transformation \(\beta^{12} \colon \Gamma_{C_1} \Longrightarrow \Gamma_{C_2}\) induced by counits is the unique natural transformation such that
\begin{equation}
\epsilon^{i_2} \circ \beta^{12} = \epsilon^{i_1}
\end{equation}
\end{lemma}

\begin{proof}
\(C_1 \subset C_2\) means \(U_2 \subset U_1\), so \(\sh_{C_1}(S) \subset \sh_{C_2}(S)\) as a full subcategory. Write \(\iota \colon \sh_{C_2}(S) \hookrightarrow \sh(S)\) as the inclusion functor, which we know has right adjoint \(R = \Gamma_{C_2}\). \(\iota R\) is \(\Gamma_{C_2}\), considered a functor with target \(\sh(S)\) instead of \(\sh_{C_2}(S)\). Then \ref{lemm:fun-fact} applied to \(\sh_{C_2}(S) = \mathcal{C}\), \(\sh(S) = \mathcal{D} = \mathcal{E}\), \(F = \Gamma_{C_1}\) (considered a functor with target \(\sh(S)\)), and \(G = \id_{\sh(S)}\) gives us that 
\[\begin{tikzcd} [column sep = large]
\Map_{\fun\big(\sh(S), \, \sh(S)\big)}(\Gamma_{C_1}, \Gamma_{C_2}) \ar[r, "\epsilon^{i_2} \circ -"] & \Map_{\fun\big(\sh(S), \, \sh(S)\big)}(\Gamma_{C_1}, \id)
\end{tikzcd}\]
is an equivalence. In 
\[\fun\big(\sh(S), \, \sh(S)\big)_{/\id}\]
the mapping space from \((\Gamma_{C_1}, \epsilon^{i_1})\) to \((\Gamma_{C_2}, \epsilon^{i_2})\) is the fibre of this map over \(\epsilon^{i_1}\). Equivalences have contractible fibres, so the unique point in this fibre up to contractible choice is \(\beta^{12} \colon \Gamma_{C_1} \Longrightarrow \Gamma_{C_2}\) such that \(\epsilon^{i_2} \circ \beta^{12} \simeq \epsilon^{i_1}\), giving the requisite factorization.

In the case that \(C_1\) is not contained in \(C_2\), the fibre is empty. So the full subcategory of this slice category on the pair \((\Gamma_C, \epsilon^i)\) has contractible or empty mapping spaces, and is thus a poset. This poset is equivalent to \(\mathcal{P}\), as \(\sh_C(S)\) determines \(C\), giving functoriality.
\end{proof}

\begin{lemma}
Take \(i, j\) as before, and a morphism \(p \colon \tilde{S} \to S\). Form the Cartesian squares
\begin{equation}\label{eq:pullbacks}
\begin{tikzcd} [row sep = large, column sep = large] \tilde{C} \ar[d, hookrightarrow, swap, "\tilde{i}"] \ar[r, "\tilde{p}", ""{name=h1}] & C \ar[d, hookrightarrow, "i"]
& & \tilde{U} \ar[d, hookrightarrow, swap, "\tilde{j}"] \ar[r, "\tilde{p}_j", ""{name=h3}] & U \ar[d, hookrightarrow, "j"]\\
\tilde{S} \ar[r, swap, "p", ""{name=h2}] & S & & \tilde{S} \ar[r, swap, "p", ""{name=h4}] & S
\ar[phantom, from=h1, to=h2, "\Delta_i"] \ar[phantom, from=h3, to=h4, "\Delta_j"] \end{tikzcd}
\end{equation}
Then \(p^!\) maps \(\sh_{C}(S)\) into \(\sh_{\tilde{C}}(\tilde{S})\).
\end{lemma}

\begin{proof}
Take \(E \in \sh_C(S)\). \(\tilde{j}\) is an open embedding, so exceptional pullback is equivalent to pullback, and commutativity of \(\Delta_j\) gives  
\[\tilde{j}^*p^!E \simeq \tilde{j}^! p^!E \simeq \tilde{p}_j^!j^!E \simeq \tilde{p}_j^! j^*E \simeq 0,\]
thus \(p^!E \in \sh_{\tilde{C}}(\tilde{S})\).
\end{proof}

\begin{construction}[Support exchange]\label{con:ex-supp}
This means that \(p^! \Gamma_C\) takes a motivic spectrum \(E \in \sh(S)\) to \(\sh_{\tilde{C}}(\tilde{S})\). Then applying \ref{lemm:fun-fact} to \(\mathcal{D} = \sh(\tilde{S})\), \(\mathcal{C} = \sh_{\tilde{C}}(\tilde{S})\), \(\mathcal{E} = \sh(S)\), and \(R = \Gamma_{\tilde{C}}\), the natural transformation \(p^!\epsilon^i \colon p^! \Gamma_C \Longrightarrow p^!\) from \(\sh(S)\) to \(\sh(\tilde{S})\), is the image under \(\epsilon^{\tilde{i}}p^!\) of a unique natural transformation of functors from \(\sh(S)\) to \(\sh(\tilde{S})\),
\[\tau^C \colon p^!\Gamma_C \Longrightarrow \Gamma_{\tilde{C}}p^!\]
The transformation \(\tau^C\) is called \emph{support exchange}, and by construction satisfies
\begin{equation}\label{eq:tau-epsilon}
\epsilon^{\tilde{i}}p^! \circ \tau^C \simeq p^! \epsilon^i
\end{equation}
\end{construction}

\begin{lemma}
\(\tau^C\) is an equivalence, induced by base change.
\end{lemma}

\begin{proof}
Use base change, and commutativity of the square \(\Delta_i\) from \eqref{eq:pullbacks}: 
\[p^! \Gamma_C = p^!i_*i^! \simeq \tilde{i}_* \tilde{p}^!i^! \simeq \tilde{i}_* (i \tilde{p})^! \simeq \tilde{i}_*\tilde{i}^!p^! = \Gamma_{\tilde{C}} p^!\]
Composing this sequence of equivalences with \(\epsilon^{\tilde{i}}p^!\) gives \(p^!\epsilon^i\), so by uniqueness of the factorization above, this natural equivalence must \(\tau^C\).
\end{proof}

\begin{remark}[Inverse of \(\tau^C\)]\label{rem:tau-inv}
As a natural equivalence, \(\tau^C\) has an inverse, \(\tilde{\tau}^C \colon \Gamma_{\tilde{C}} p^! \Longrightarrow p^! \Gamma_C\), which satisfies 
\begin{equation}\label{eq:twid-tau-epsilon}
p^! \epsilon^i \circ \tilde{\tau}^C = \epsilon^{\tilde{i}} p^!
\end{equation}
\end{remark}


\begin{remark}[Notation] Throughout the remainder of this section, \(n \in \bb{Z}\) will remain fixed, we will write \(d = n + q\) for a given \(q \geq 0\), and \(f \colon X \to \spec \bb{Z}\) will be the structure morphism. Given \(U \in X_\nis\), we write \(f_U \colon U \to \spec \bb{Z}\) for the structure morphism, and \(f_{U, q}\) for the structure morphism of \(U \times \Delta^q\), so that these morphisms satisfy \(f_U \circ \pr_1 = f_{U, q}\).
\end{remark}

\begin{construction}[Spectral Replacement]\label{con:spec-rep}
For \(q \geq 0\), let \(\mathcal{S}_n(U, q)\) be the poset category of closed subschemes \(i \colon W \hookrightarrow U \times \Delta^q\) satisfying \(\delta(W) \leq d\) and which meeting all faces properly. Note that since \(\delta(W_1 \cup W_2) = \max \{\delta(W_1), \delta(W_2)\}\), and \(W_1 \cup W_2\) meets faces properly, \(\mathcal{S}_n(U, q)\) is a filtered category. For \(W \in \mathcal{S}_n(U, q)\), define the spectrum 
\[F_U^W = \map_{\sh(U \times \Delta^q)}\big(\mathbbm{1}_{U \times \Delta^q}, \Gamma_W f_{U, q}^! s^{- d} \Kgl_\bb{Z}\big),\]
Given \(W_1, W_2 \in \mathcal{S}_n(U, q)\), with \(W_1 \subset W_2\), by \ref{lemm:supp-fun} there is a natural transformation \(\Gamma_{W_1} \Longrightarrow \Gamma_{W_2}\), and this induces a morphism of spectra, \(F_U^{W_1} \to F_U^{W_2}\), so the assignment \(\mathcal{S}_n(U, q) \to \Sp\), \(W \mapsto F_U^W\) is a functor. Denote the colimit of this functor as 
\[F_n(U, q) = \colim_{W \in \mathcal{S}_n(U, q)} F_U^W\]
\end{construction}

\begin{definition} We call both \(F_n( -, \bullet) \in \psh(X_\nis \times \Delta, \Sp)\) and its geometric realization \(|F_n( -, \bullet)| \in \psh(X_\nis, \Sp)\) the \emph{spectral replacement of \(\mathcal{Z}_{n + 1}(X, *)\)}.
\end{definition}

We table proving that this is an object of \(\fun(X_\nis \times \Delta, \Sp)\), as it will follow from computations of the homotopy groups of \(F_n(U, q)\).

\subsection{Comparison to \(\mathcal{Z}_{n + 1}(X, *)\)}
In this section, we prove that the hyperhomology groups of the Bloch-Levine cycle complex may be computed as the homotopy groups of the global sections of the spectral replacement. The following theorem, due to Frankland-Spitzweck \cite[Proposition 8.3]{FrSp}, will feature prominently:

\begin{theorem}[Arithmetic Purity]\label{thm:ar-pure}
Let \(i_p \colon \spec \bb{F}_p \hookrightarrow \spec \bb{Z}\) be the canonical closed embedding. Then the purity transformation \(\p_{i_p} \colon s^0 \Kgl_{\bb{F}_p} \otimes \bb{T}^{\otimes - 1} \to i_p^! s^0 \Kgl_\bb{Z}\) is an equivalence.
\end{theorem}

\begin{remark}[\(H\bb{Z}\)]
We use the notation \(H\bb{Z} = s^0 \Kgl_\bb{Z} \in \sh(\bb{Z})\) since this particular motivic spectrum will frequently appear in this section. This differs from the standard topological Eilenberg-Mac Lane spectrum, but that spectrum will not make any further appearances, so there should be no confusion.
\end{remark}

\begin{lemma}[Borel-Moore homology of \(H\bb{Z}\)]\label{lemm:HZBM}
Let \(f_T \colon T \to \spec \bb{Z}\) be finite type separated, with \(\delta(T) \leq d\). Then 
\begin{equation}\label{eq:HZBM}
H\bb{Z}_{*, d}^\bm(T/\bb{Z}) \cong \begin{cases} \displaystyle \bigoplus_{x \in T_{(d)}} \bb{Z} & * = 2d\\
0 & \mathrm{else} \end{cases}
\end{equation}
\end{lemma}

\begin{proof}
We use the spectral sequence \eqref{eq:dnss} with \(i = d\) and \(E = H\bb{Z}\) to compute \(H\bb{Z}_{*, d}^\bm(T/\bb{Z})\). Each \(x \in T_{(a)}\) contributes a direct summand of \(E_{a, b}^1\),  
\[\hat{H\bb{Z}}_{a + b, d}^\bm(x/\bb{Z}) = \colim_{V \subset Z(x)} \pi_0 \! \left(\map_{\sh(V)}\big(\mathbbm{1}_V \otimes \bb{T}^{\otimes d}[a + b - 2d], f_V^! H\bb{Z}\big)\right)\]
\[\cong \colim_{V \subset Z(x)} \pi_{a + b - 2d}\left(\map_{\sh(V)}\big(\mathbbm{1}_V, \bb{T}^{\otimes - d} \otimes f_V^! H\bb{Z}\big)\right) \Longrightarrow H\bb{Z}_{a + b, d}^\bm(T/\bb{Z})\]

Fix \(x \in T\). \(Z(x)\), an integral scheme, is regular at its generic point \(x\). It is also finite type over the excellent ring \(\bb{Z}\), meaning its regular locus is open. So we may assume without loss of generality that the \(V\) in the colimit are regular schemes. Further assume that \(V = \spec R\) is affine. Write \(s = f_V(x)\). We have two possibilities. 
\begin{enumerate}
\item \(s = 0\). Then \(\mathcal{O}_{\bb{Z}, s} \to \mathcal{O}_{V, x}\) is the field extension \(\bb{Q} \hookrightarrow \kappa(x)\), so \(f_V\) is flat at \(x\). The fibre of \(f_V\) over \(s\) is \(V_\bb{Q} = \spec (R \otimes_\bb{Z} \bb{Q}) = \spec R_\bb{Q}\). \(R_\bb{Q}\), the localization of \(R\) at \(\bb{Z} - \{0\}\), is an integral domain, and its spectrum has the same generic point, \(x\). Thus \(V_\bb{Q}\) is regular at \(x\), and as a finite type scheme over a perfect field, regularity at \(x\) implies smoothness of \(V_s \to \spec \kappa(s)\) at \(x\).

\cite[Tag 01V9]{Sta} then implies that \(f\) is smooth at \(x\) of relative dimension \(\dim_x V_\bb{Q} = \trd\big(\kappa(x) : \kappa(0)\big)\), which is precisely \(\delta(x)\) since \(0\) has height zero. As a finite type scheme over the excellent ring \(\bb{Z}\), the smooth locus of \(f_V\) is open. So by cofinality, we may further assume without loss of generality that each \(f_V \colon V \to \spec \bb{Z}\) is smooth of relative dimension \(\delta(x)\). Smooth purity allows us to make the conversion 
\[f_V^!H\bb{Z} \simeq \bb{T}^{\otimes \delta(x)} \otimes f_V^* s^0 \Kgl_\bb{Z} \simeq \bb{T}^{\otimes \delta(x)} \otimes s^0 \Kgl_V.\]

\item \(s = p\). As a closed point of \(\spec \bb{Z}\), we have that \(f_V(V) = f_V(\overline{\{x\}}) \subset \overline{f_V(x)} = \overline{\{p\}} = \spec \kappa(p)\), so we may factor \(f_V\) as 
\[\begin{tikzcd} V \ar[rr, "f_V", bend right = 18, swap] \ar[r, "f_{V, p}"] & \spec \bb{F}_p \ar[r, "i_p", hookrightarrow] & \spec \bb{Z} \end{tikzcd}\]
Applying the reasoning from the \(s = 0\) case, we may assume without loss of generality that each \(f_{V, p}\) is smooth at \(x\) of relative dimension \(\trd\big(\kappa(x) : \kappa(p)\big)\), which is precisely \(\delta(x) + 1\) since \(p\) has height one. Smooth purity combined with arithmetic purity \ref{thm:ar-pure}, we see that 
\[f_V^! H\bb{Z} = f_{V, p}^! i_p^!s^0 \Kgl_\bb{Z} \simeq \bb{T}^{\delta(x) + 1} \otimes f_{V, p}^* \big(\bb{T}^{- 1} \otimes s^0 \Kgl_{\bb{F}_p}\big) \simeq \bb{T}^{\delta(x)} \otimes f_{V, p}^* s^0 \Kgl_{\bb{F}_p} \simeq \bb{T}^{\otimes \delta(x)} \otimes s^0 \Kgl_V\footnote{At this point, it should be apparent why we modified the definition of relative dimension}.\]
\end{enumerate}
In either case, \(\bb{T}^{\otimes - d} \otimes f_V^! H\bb{Z} \simeq s^{\delta(x) - d} \Kgl_V = s^{a - d} \Kgl_V\), meaning 
\[\hat{H\bb{Z}}_{a + b, d}^\bm(x/\bb{Z}) = \colim_{V \subset Z(x)} \pi_{a + b - 2d} \! \left(\map_{\sh(V)}\big(\mathbbm{1}_V, s^{a - d} \Kgl_V\big)\right)\]
The cofinal system \(\{f_V \colon V \to \spec \bb{Z}\}_{V \subset Z(x)}\) has limit \(f_x \colon \spec \kappa(x) \to \spec \bb{Z}\). By continuity of \(\sh\) \cite[1.3.19]{BonDeg} and commutativity of homotopy groups with filtered colimits, 
\[\hat{H\bb{Z}}_{a + b, d}^\bm(x/\bb{Z}) = \pi_{a + b - 2d} \! \left(\map_{\sh(\kappa(x))}\big(\mathbbm{1}_{\kappa(x)}, s^{a - d} \Kgl_{\kappa(x)}\big)\right)\]
That mapping spectrum is just a shift of \(\bb{Z}(a - d)^\bb{A}\big(\kappa(x)\big)\), so this simplifies to
\[\pi_{a + b - 2d} \left(\bb{Z}(a - d)^\bb{A}\big(\kappa(x)\big)[2(a - d)]\right) \cong \pi_{b- a} \! \left(\bb{Z}(a - d)^\bb{A}\big(\kappa(x)\big)\right)\]
\(\bb{A}^1\)-motivic cohomology vanishes in negative weights \cite[Remark 4.20]{BEM}, so these groups are zero  when \(a < d\). And \(a\), the relative dimension of \(x \in T\), may not exceed \(d\), the maximum dimension of \(T\). This means that the page \(E_{a, b}^1\) has zero terms except in the \(d^\mathrm{th}\) column.

Using \cite[Corollary 9.1.2]{BEM} again, we have an equivalence \(\bb{Z}(0)^\bb{A}\big(\kappa(x)\big) \simeq R\Gamma_\mathrm{cdh}\big(\kappa(x) ,\bb{Z}\big)\). This complex is concentrated at degree zero, with \(H^0 R\Gamma_\mathrm{cdh}\big(\kappa(x), \bb{Z}\big) \cong \bb{Z}\). Thus the rows are all zero except at \(b = a = d\), in which there is one direct summand  \(\bb{Z}\) for each \(x \in T_{(d)}\). As a result, the first page degenerates, and so is isomorphic to the abutment, giving the desired result.
\end{proof}

\begin{definition}
Given \(x \in T_{(d)}\), we write \(\eta_x\) for the corresponding generator of the Borel-Moore homology group. Given an integral scheme \(C\), we avoid the cumbersome notation \(\eta_{\eta_C}\) by simply writing \(\eta_C\). This notation is reminiscent of the fundamental class: by \ref{thm:fund-class}, the purity isomorphism
\[H\bb{Z}_{2d, d}^\bm(T/\bb{Z}) \cong H^0 R\Gamma(X_\mathrm{cdh}, \bb{Z}) \cong \bb{Z}\]
in the calculation with the \(\delta\)-niveau spectral sequence is precisely given by the fundamental class of \(f_V \colon V \to \spec \bb{Z}\). As \ref{lemm:BM-excise} will show, this only depends on \(x\), further justifying the notation. We denote by 
\[\gamma_T \colon \begin{tikzcd} H\bb{Z}_{2d, d}^\bm(T/\bb{Z}) \ar[r, "\sim"] & \displaystyle \bigoplus_{x \in T_{(d)}} \bb{Z} \cdot \eta_x \end{tikzcd}\]
the isomorphism from the previous lemma, and call it the \emph{cycle map}. Note that in the case of \(\delta(T) < d\), this is simply the zero morphism.
\end{definition}

\begin{lemma}\label{lemm:EM}
For any \(U \in X_\nis\) and \(q \geq 0\), there is an equivalence of spectra 
\begin{equation}\label{eq:EM}
F_n(U, q) \simeq Hz_{n +1}(U, q),
\end{equation}
for \(z_{n + 1}(U, q)\) the free Abelian group defined in \ref{def:LCC}
\end{lemma}

\begin{proof}
Take \(W \in \mathcal{S}_n(U, q)\). Writing \(f_W \colon W \to \spec \bb{Z}\) for the structure morphism, adjointness and the factorization \(f_W = f_{U, q} \circ i\) gives isomorphisms 
\[\pi_k \! \left(F_U^W\right) = \pi_k \! \left(\map_{\sh(U \times \Delta^q)}\big(\mathbbm{1}_{U \times \Delta^q}, i_* i^! f_{U, q}^! s^{- d} \Kgl_\bb{Z}\big)\right) \cong \pi_k \! \left(\map_{\sh(W)}\big(\mathbbm{1}_W, f_W^! s^{- d}\Kgl_\bb{Z}\big)\right)\]
\[\cong \pi_k \! \left(\map_{\sh(W)}\big(\bb{T}^{\otimes d}, f_W^! s^0 \Kgl_\bb{Z}\big)\right) \cong \pi_0 \! \left(\map_{\sh(W)}\big(\bb{T}^{\otimes d}[k], f_W^!H\bb{Z}\big)\right) = H\bb{Z}_{2d + k, d}^\bm(W/\bb{Z})\]
So by \eqref{eq:HZBM}, we have 
\[\pi_k \! \left(F_U^W\right) \cong \begin{cases} \displaystyle \bigoplus_{x \in W_{(d)}} \bb{Z} \cdot \eta_x & k = 0\\
0 & \text{else}
\end{cases}\]
Write \(M_W = \pi_0(F_U^W)\). Then \(F_U^W \simeq HM_W\), an Eilenberg-Mac Lane spectrum. Given \(i_{12} \colon W_1 \to W_2\) in \(\mathcal{S}_n(U, q)\), the resulting map of spectra \(F_U^{W_1} \to F_U^{W_2}\) is induced by the counit \(\beta^{i_{12}} \) in the target. The isomorphism \(\gamma_{W_1}\) is precisely using adjointness to remove application of supports \(\Gamma_{W_1}\), and likewise with \(\gamma_{W_2}\) and \(\Gamma_{W_2}\).

The map of Borel-Moore homology groups \((i_{12})_* \colon H\bb{Z}_{2d, d}^\bm(W_1/\bb{Z}) \to H\bb{Z}_{2d, d}^\bm(W_2/\bb{Z})\) is by definition induced by the counit \(\epsilon^{i_{12}}\), which is \(\beta^{i_{12}} = (i_2)_* \epsilon^{i_{12}} i_2^!\) with supports removed. Thus by \ref{lemm:cyc-emb}, the morphism 
\[\pi_0 \! \left(F_U^{W_1}\right) \to \pi_0 \! \left(F_U^{W_2}\right)\]
is, under identification of \(\gamma_{W_1}\) and \(\gamma_{W_2}\), inclusion of direct summands on dimension \(d\) points \(MW_1 \hookrightarrow MW_2\). This gives a morphism of spectra \(HM_{W_1} \to HM_{W_2}\). Formation of Eilenberg-Mac Lane spectra commutes with filtered colimits, so 
\[F_n(U, q) = \colim_{W \in \mathcal{S}_n(U, q)} F_U^W \simeq \colim_{W \in \mathcal{S}_n(U, q)} HM_W \simeq H \! \left(\colim_{W \in \mathcal{S}_n(U, q)} M_W\right)\]
Since every \(W \in \mathcal{S}_n(U, q)\) meets all faces properly, so does \(Z(x)\), for each \(x \in W_{(d)}\). Each integral closed subscheme \(Z \subset U \times \Delta^q\) with \(\delta(Z) = d\) meeting all faces properly has a unique generic point \(\eta\) with \(\delta(\eta) = d\), so that \(Z = Z(\eta)\), meaning there is a bijection between \(x \in W_{(d)}\) for \(W \in \mathcal{S}_n(U, q)\) and integral closed subschemes of \(U \times \Delta^q\) meeting all faces properly. The transition maps in the colimit are inclusions, so 
\[\colim_{W \in \mathcal{S}_n(U, q)} M_W \cong \bigoplus_{x \in (U \times \Delta^q)_{(d)}} \bb{Z} \cdot \eta_x,\]
in which \(Z(x)\) meets all faces properly. So the colimit is the free Abelian group on all integral closed subschemes of \(U \times \Delta^q\) with \(\delta = d\), the ``Levine dimension'' of which is \(d + 1 = n + 1 + q\). This group is precisely \(z_{n + 1}(U, q)\), so \(F_n(U, q) \simeq Hz_{n + 1}(U, q)\).
\end{proof}

\begin{lemma}[Determined by \(\pi_0\)]\label{lemm:pi-nought}
\(F_n( -, \bullet)\) is a presheaf of simplicial spectra exactly if \(\pi_0 F_n( -, \bullet)\) is a presheaf of simplicial Abelian groups.
\end{lemma}

\begin{proof}
By \ref{lemm:EM}, the assignment \(F_n( -, \bullet) \colon X_\nis^\op \times \Delta^\op \to \Sp\) factors through \(\Sp^\heartsuit\). \(\pi_0 \colon \Sp^\heartsuit \to N(\mathrm{Ab})\) is an equivalence, so restriction along \(\pi_0\) yields an equivalence 
\[\fun \! \left( \!N(X_\nis^\op \times \Delta^\op), \Sp^\heartsuit\right) \simeq N\big(\psh(X_\nis \times \Delta, \mathrm{Ab})\big)\]
Thus the structure of a simplicial presheaf of Abelian groups on \(\pi_0\) lifts to a structure of a simplicial presheaf of spectra up to contractible choice.
\end{proof}

This allows us to construct structure morphisms on the level of spectra, and show that they satisfy the necessary compatibilities at the level of \(\pi_0\).

\begin{lemma}\label{lemm:senh-struc}
Fix \(U \in X_\nis\), \(q \in \Delta\). We have the following structure morphisms: 
\begin{enumerate}[i).]

\item Face maps: \(\iota \colon U \times \Delta^{q - 1} \hookrightarrow U \times \Delta^q\) induces a morphism of spectra, \(F_U^W \to F_U^{\iota^{- 1}(W)}\).

\item Degeneracy maps: \(p \colon U \times \Delta^{q + 1} \to U \times \Delta^q\) induces a morphism of spectra, \(F_U^W \to F_U^{p^{- 1}(W)}\).

\item Étale maps: A morphism \(V \to U\) in \(X_\nis\) induces a morphism of spectra \(F_U^W \to F_V^{g^{- 1}(W)}\), in which \(g \colon V \times \Delta^q \to U \times \Delta^q\) is the product.

\end{enumerate}
\end{lemma}

\begin{proof}
\item (i). Let \(\iota\) be any coface inclusion as above. Given \(W \in \mathcal{S}_n(U, q)\), form the pullback 
\begin{equation}\label{eq:face-square}
\begin{tikzcd} [row sep = large] \iota^{- 1}(W) \ar[d, "\tilde{i}", swap] \ar[r, "\tilde{\iota}"] & W \ar[d, "i"]\\
U \times \Delta^{q - 1} \ar[r, "\iota", swap] & U \times \Delta^q \end{tikzcd}
\end{equation}
Since \(W\) meets all faces properly, \(\iota^{- 1}(W)\), the intersection with a face, has \(\delta\) at most \(d - 1\). Likewise  \(\iota^{- 1}(W)\) meets all faces of \(U \times \Delta^{q - 1}\), themselves faces of \(U \times \Delta^q\), properly. Thus \(\iota^{- 1}(W) \in \mathcal{S}_n(U, q - 1)\), and compatibility of pullback with closed embeddings gives that \(\iota^{- 1} \colon \mathcal{S}_n(U, q) \to \mathcal{S}_n(U, q - 1)\) is a functor.  Applying the functor \(\iota^!\) to an element of \(F_U^W\) gives a morphism
\[\iota^!\mathbbm{1}_{U \times \Delta^q} \to \iota^!\Gamma_W f_{U, q}^! s^{- d} \Kgl_\bb{Z}\]
Projection from \(U \times \Delta^q\) and \(U \times \Delta^{q - 1}\) is smooth, so by relative purity, the source is equivalent to \(\mathbbm{1}_{U \times \Delta^{q - 1}} \otimes \bb{T}^{- 1}\), and application of support exchange \ref{con:ex-supp} to \(S = U \times \Delta^q\), \(\tilde{S} = U \times \Delta^{q - 1}\), \(p = \iota\) gives an equivalence
\[\tau^W \colon \iota^! \Gamma_W \simeq \Gamma_{\iota^{- 1}(W)}\iota^!\]
Since \(\iota^! \circ f_{U, q}^! = f_{U, q - 1}^!\), we have a composite morphism 
\[\begin{tikzcd} \mathbbm{1}_{U \times \Delta^{q - 1}} \otimes \bb{T}^{- 1} \ar[r, "\iota^! \phi"] & \iota^! \Gamma_W f_{U, q}^! s^{- d} \Kgl_\bb{Z} \ar[r, "\tau^W"] & \Gamma_{\iota^{- 1}(W)} f_{U, q - 1}^! s^{- d} \Kgl_\bb{Z} \end{tikzcd}\]
Smashing by \(\bb{T}\) is an equivalence which commutes with \(\Gamma_{\iota^{- 1}(W)} f_{U, q - 1}^!\) and converts \(s^{- d}\Kgl_\bb{Z}\) to \(s^{- (d - 1)} \Kgl_\bb{Z}\), giving an element of the spectrum \(F_U^{\iota^{- 1}(W)}\). Thus \(\iota\) induces a morphism of spectra, 
\[\iota^! \colon F_U^W \to F_U^{\iota^{- 1}(W)}\]

\item (ii). Let \(p\) be a codegeneracy morphism corresponding to \(t_i + t_{i + 1}\) in barycentric coordinates. Take the pullback 
\[\begin{tikzcd} [row sep = large]
p^{- 1}(W) \ar[d, "\tilde{i}", swap] \ar[r, "\tilde{p}"] & W \ar[d, "i"]\\
U \times \Delta^{q + 1} \ar[r, "p", swap] & U \times \Delta^q
\end{tikzcd}\]
Then by \cite[Tag 02R8]{Sta}, as a smooth morphism of relative dimension 1, 
\[\delta\big(p^{- 1}(W)\big) = \delta(W) + 1 \leq d + 1,\]
and an analogous result on all faces ensures proper intersection, meaning \(p^{- 1}(W) \in \mathcal{S}_n(U, q + 1)\), and \(p^{- 1} \colon \mathcal{S}_n(U, q) \to \mathcal{S}_n(U, q + 1)\) is a functor. Applying the functor \(p^!\) to an element of \(F_U^W\) gives a morphism 
\[p^! \psi \colon p^!\mathbbm{1}_{U \times \Delta^q} \to p^! \Gamma_W f_{U, q}^! s^{- d} \Kgl_\bb{Z}\]
Smooth purity gives us equivalences \(p^! \mathbbm{1}_{U \times \Delta^q} \simeq p^* \mathbbm{1}_{U \times \Delta^q} \otimes \mathrm{Th}(T_p) \simeq \mathbbm{1}_{U \times \Delta^{q + 1}} \otimes \bb{T}\). Support exchange \ref{con:ex-supp} applied to \(S = U \times \Delta^q\), \(\tilde{S} = U \times \Delta^{q + 1}\) gives an equivalence 
\[\tau^W \colon p^! \Gamma_W \simeq \Gamma_{p^{- 1}(W)} \, p^!\]
Since \(p^! f_{U, q}^! = f_{U, q + 1}^!\), we have a composite morphism 
\[\begin{tikzcd}
\mathbbm{1}_{U \times \Delta^{q + 1}} \otimes \bb{T} \ar[r, "p^! \psi"] & p^! \Gamma_W f_{U, q}^! s^{- d} \Kgl_\bb{Z} \ar[r, "\tau^W"] & \Gamma_{p^{- 1}(W)} f_{U, q + 1}^! s^{- d} \Kgl_\bb{Z}
\end{tikzcd}\]
Smashing by \(\bb{T}^{- 1}\) is an equivalence which commutes with \(\Gamma_{p^{- 1}(W)} f_{U, q + 1}^!\) and converts \(s^{- d} \Kgl_\bb{Z}\) to \(s^{- (d + 1)} \Kgl_\bb{Z}\), giving an element of the spectrum \(F_U^{p^{- 1}(W)}\). Thus \(p\) induces a morphism of spectra,
\[p^! \colon F_U^W \to F_U^{p^{- 1}(W)}\]

\item (iii). A morphism \(V \to U\) in \(X_\nis\) is étale, thus so is its product, \(g \colon V \times \Delta^q \to U \times \Delta^q\). Given \(W \in \mathcal{S}_n(U, q)\), we may flatly pull this back to \(g^{- 1}(W) \in \mathcal{S}_n(V, q)\): 
\[\begin{tikzcd} [row sep = large]
g^{- 1}(W) \ar[d, "\tilde{i}", hookrightarrow, swap] \ar[r, "\tilde{g}"] & W \ar[d, "i", hookrightarrow]\\
V \times \Delta^q \ar[r, "g", swap] & U \times \Delta^q
\end{tikzcd}\]
As an étale morphism, \(g\) is flat of relative dimension zero, so by  \cite[Tag 02R8]{Sta} again, \(\delta\big(g^{- 1}(W)\big) = \delta(W) \leq d\). Similar reasoning applied to the pullback of the faces of \(U \times \Delta^q\) implies that \(g^{- 1}(W)\) meets all faces of \(V \times \Delta^q\) properly. So \(g^{- 1} \colon \mathcal{S}_n(U, q) \to \mathcal{S}_n(V, q)\) is a functor. Applying the functor \(g^!\) to an element of \(F_U^W\)  gives a morphism 
\[g^!\xi \colon g^*\mathbbm{1}_{U \times \Delta^q} \to g^!\Gamma_W f_{U, q}^! s^{- d} \Kgl_\bb{Z}\]
Since \(g\) is étale, \(g^! \simeq g^*\) by purity, so the unit is pulled back to the unit of \(\sh(V \times \Delta^q)\). Support exchange \ref{con:ex-supp} applied to \(\tilde{S} = V \times \Delta^q\), \(S = U \times \Delta^q\), and \(p = g\) yields an equivalence 
\[\tau^W \colon g^! \Gamma_W \simeq \Gamma_{g^{- 1}(W)} g^!\]
\(g\) is the product of the morphism \(V \to U\) over \(X\), meaning \(g^! f_{U, q}^! = f_{V, q}^!\), thus we have a composite morphism 
\[\begin{tikzcd}
\mathbbm{1}_{V \times \Delta^q} \ar[r, "g^! \xi"] & g^! \Gamma_W f_{U, q}^! s^{- d} \Kgl_\bb{Z} \ar[r, "\tau^W"] & \Gamma_{g^{- 1}(W)} f_{V, q}^! s^{- d} \Kgl_\bb{Z},
\end{tikzcd}\]
which is an element of the spectrum \(F_V^{g^{- 1}(W)}\). Thus \(g\) induces a morphism of spectra,
\[g^! \colon F_U^W \to F_V^{g^{- 1}(W)},\]
so we have all three structure morphisms.
\end{proof}

\(\G\)-theory, \(H\bb{Z}^\bm\), and \(z_{n + 1}( -, \bullet)\) all have nice functorialities, and the following results show relations between them, which will eventually be used to prove the desired structure result on \(\pi_0 F_n( -, \bullet)\).

\begin{lemma}[Cycle map and embeddings]\label{lemm:cyc-emb}
Let \(j \colon U \hookrightarrow T\) be an open embedding. Then the diagram 
\[\begin{tikzcd} [column sep = large] H\bb{Z}_{2d, d}^\bm(T/\bb{Z}) \ar[d, "j^*", swap] \ar[r, "\gamma_T"] & \displaystyle \bigoplus_{x \in T_{(d)}} \bb{Z} \ar[d]\\
H\bb{Z}_{2d, d}^\bm(U/\bb{Z}) \ar[r, "\gamma_U", swap] & \displaystyle \bigoplus_{x \in U_{(d)}} \bb{Z} \end{tikzcd}\]
commutes, for \(j^*\) the étale pullback of Borel-Moore homology groups in \ref{def:BM}. Likewise, for a closed embedding \(i \colon C \hookrightarrow T\), the diagram 
\[\begin{tikzcd} [column sep = large] H\bb{Z}_{2d, d}^\bm(C/\bb{Z}) \ar[d, "i_*", swap] \ar[r, "\gamma_C"] & \displaystyle \bigoplus_{x \in C_{(d)}} \bb{Z} \ar[d]\\
H\bb{Z}_{2d, d}^\bm(T/\bb{Z}) \ar[r, "\gamma_T", swap] & \displaystyle \bigoplus_{x \in T_{(d)}} \bb{Z} \end{tikzcd}\]
commutes, for \(i_*\) the proper pushforward of Borel-Moore homology groups in \ref{def:BM}.
\end{lemma}

\begin{proof}
The \(\delta\)-niveau spectral sequence arises from the filtration of \(H\bb{Z}_{*, d}^\bm(T/\bb{Z})\) by supports: 
\[F_a  = \bigcup_{\delta(Z) \leq a} \im \! \Big ((i_Z)_* \colon H\bb{Z}_{*, d}^\bm(Z/\bb{Z}) \to H\bb{Z}_{*, d}^\bm(T/\bb{Z})\Big)\]
Both \(i_*\) and \(j^*\) preserve this filtration, so are compatible with it on the \(E^1\) page, which is isomorphic to the abutment via \(\gamma\).
\end{proof}

\begin{lemma}[Generic determination of Borel-Moore homology groups]\label{lemm:BM-excise}
Given open \(j \colon U \hookrightarrow T\), with \(i \colon T - U \hookrightarrow T\) the complementary closed embedding, then if \(\delta(T - U) < d\), \(j^* \colon H\bb{Z}_{2d, d}^\bm(T/\bb{Z}) \to H\bb{Z}_{2d, d}^\bm(U/\bb{Z})\) is an isomorphism.
\end{lemma}

\begin{proof}
As a finite type scheme over \(\spec \bb{Z}\), \(T\) is Noetherian, so \(j\) is a quasicompact open embedding, and we may apply the localization sequence to get exact sequence 
\[\begin{tikzcd} H\bb{Z}_{2d, d}^\bm(T - U/\bb{Z}) \ar[r, "i_*"] & H\bb{Z}_{2d, d}^\bm(T/\bb{Z}) \ar[r, "j^*"] & H\bb{Z}_{2d, d}^\bm(U/\bb{Z}) \ar[r] & H\bb{Z}_{2d - 1, d}^\bm(T - U/\bb{Z}), \end{tikzcd}\]
the outer terms of which vanish due to the dimension assumptions on \(T - U\).
\end{proof}

\begin{remark}[Generators and classes]
Denote by \([C] \in H\bb{Z}_{2d, d}^\bm(C/\bb{Z})\) the inverse image \(\gamma_C^{- 1}(\eta_C)\). By \ref{lemm:cyc-emb} and \ref{lemm:BM-excise}, this is equivalently the unique class restricting to 
\[1 \in \bb{Z} \cong \hat{H\bb{Z}}_{2d, d}^\bm(\eta_C/\bb{Z})\]
Thus for any scheme \(T\), each class in \(H\bb{Z}_{2d, d}^\bm(T/\bb{Z})\) may be uniquely written as 
\[\sum_{x \in T_{(d)}} n_x (i_x)_* [Z(x)],\]
for \(i_x \colon Z(x) \hookrightarrow T\) the closed embedding of the integral scheme \(Z(x)\), and \(n_x \in \bb{Z}\).
\end{remark}

\begin{definition}\label{def:stalky}
For \(\delta(T) \leq d\), \(x \in T_{(d)}\), let \(V_x \subset T\) be an integral open subscheme containing \(x\) as its generic point. 
By \ref{lemm:BM-excise}, there is an isomorphism 
\[H\bb{Z}_{2d, d}^\bm(V_x/\bb{Z}) \cong \hat{H\bb{Z}}_{2d, d}^\bm(x/\bb{Z}) \cong \bb{Z} \]
independent of the choice of such a \(V_x\). Writing \(j_x \colon V_x \hookrightarrow T\) as the open embedding, we have a map
\[\begin{tikzcd} H\bb{Z}_{2d, d}^\bm(T/\bb{Z}) \ar[r, "j_x^*"] & H\bb{Z}_{2d, d}^\bm(V_x/\bb{Z}) \ar[r, "\sim"] & \hat{H\bb{Z}}_{2d, d}^\bm(x/\bb{Z}) \ar[r, "\sim"] & \bb{Z} \end{tikzcd}\]
Call the image of a class \(\alpha\) under this map \(\alpha|_x\). We see that 
\begin{equation}\label{eq:functionals}
\gamma_T(\alpha) = \sum_{x \in T_{(d)}} \alpha|_x \cdot \eta_x,
\end{equation}
as \(\gamma\) commutes with pulling back by opens and is \(\bb{Z}\)-linear.
\end{definition}

\begin{definition}
Denote \(\kgl_\bb{Z} = \fil_\sli^0 \Kgl_\bb{Z}\) for the connective cover, and write \(\nu \colon \kgl_\bb{Z} \to H\bb{Z}\) for the canonical morphism \(\fil_\sli^0 \Kgl_\bb{Z} \to s^0 \Kgl_\bb{Z}\). 
Periodicity of \(\kgl_\bb{Z}\)  allows us to lift the orientation of \(\Kgl_\bb{Z}\) to \(\kgl_\bb{Z}\), so \(\nu\) is compatible with Thom isomorphisms, which is necessary for later purity results.
\end{definition}

\begin{lemma}
For finite type separated \(f_T \colon T \to \spec \bb{Z}\) with \(\delta(T) \leq d\), the canonical morphism \(\kgl_\bb{Z} \to \Kgl_\bb{Z}\) induces an isomorphism 
\[\kgl_{a, d}^\bm(T/\bb{Z}) \cong \Kgl_{a, d}^\bm(T/\bb{Z}),\]
for each \(a \in \bb{Z}\). 
\end{lemma}

\begin{proof}
Exhaustiveness of the slice filtration on \(\Kgl_\bb{Z}\) means that we may write the cofibre \(C\) of \(\kgl_\bb{Z} \to \Kgl_\bb{Z}\) as 
\[C \simeq \colim_{j \to - \infty} \fil_\sli^j \Kgl_\bb{Z} / \fil^0 \Kgl_\bb{Z}\]
Letting \(j < 0\), such a term is a finite extension of the slice \(s^i \Kgl_\bb{Z}\), for \(j < i < 0\). 
\[(s^i \Kgl_\bb{Z})_{a, d}^\bm(T/\bb{Z}) \cong H\bb{Z}_{a - 2i, d - i}^\bm(T/\bb{Z})\]
by Bott periodicity, so \ref{lemm:HZBM} tells us that this group is zero, as there are no points of relative dimension \(d - i > d \geq \delta(T)\). Compactness of \((f_T)_!\mathbbm{1}_T\) means that \(C_{a, d}^\bm(T/\bb{Z}) \cong 0\), so we have the desired result.
\end{proof}

\begin{definition}
For \(\delta(T) \leq d\) as before, write \(\nu_d\) for the composite 
\[\begin{tikzcd}
G_0(T) \ar[r, "\sim"] & \Kgl_{2d, d}^\bm(T/\bb{Z}) \ar[r, "\sim"] & \kgl_{2d, d}^\bm(T/\bb{Z}) \ar[r, "\nu"] & H\bb{Z}_{2d, d}^\bm(T/\bb{Z})
\end{tikzcd}\]
in which the first isomorphism comes from Jin's identification of \(\G\)-theory as the Borel-Moore theory associated to \(\Kgl\). Jin proves that \(\nu_d\) is compatible with proper pushforward \cite[Lemma 3.3.4]{Jin}, and restriction to open subschemes \cite[Lemma 3.3.5]{Jin}.
\end{definition}

\begin{proposition}[Cycle maps of coherent sheaves]\label{prop:cyc-coh}
Let \(\mathcal{F}\) be a coherent sheaf on \(T\). Then 
\begin{equation}
\gamma_T\big(\nu_d[\mathcal{F}]\big) = \sum_{x \in T_{(d)}} \ell_{\mathscr{O}_{T, x}}(\mathcal{F}_x) \cdot \eta_x
\end{equation}
\end{proposition}

\begin{proof}
By \ref{eq:functionals}, we must prove that for any \(x \in T_{(d)}\) and coherent \(\mathcal{F}\), the image of \(\nu_d[\mathcal{F}]\) under the isomorphism
\[\hat{H\bb{Z}}_{2d, d}^\bm(x/\bb{Z}) \cong \bb{Z}\]
is \(\ell(\mathcal{F}_x)\). Commutativity of \(\nu_d\) with open restriction means that the image of \(\nu_d[\mathcal{F}]\) is 
\[\begin{tikzcd}
\displaystyle \colim_{x \in V} G_0(V)\ar[r] & \displaystyle \colim_{x \in V} H\bb{Z}_{2d, d}^\bm(V/\bb{Z})
\end{tikzcd}\]
Continuity of G-theory \cite[Proposition 2.2]{Qu} exhibits the first colimit as \(G_0\big(\mathscr{O}_{T, x}\big)\). \(x\) is a top-dimensional point of \(T\), so the local ring \(\mathscr{O}_{T, x}\) is Artinian. As a result, the length function 
\[[\mathscr{G}] \mapsto \ell_{\mathscr{O}_{T, x}}(\mathscr{G})\]
provides an isomorphism \(G_0(\mathscr{O}_{T, x}) \cong \bb{Z}\) \cite[Tag 00JI]{Sta}. We must show that the generator \([\kappa(x)]\) is mapped to 1 under such an isomorphism.

Choose a regular open \(V \subset Z(x)\) as in the proof of \ref{lemm:HZBM}, so that \(V\) is smooth over \(\spec \bb{F}_p\) or \(\spec \bb{Z}[1/N]\), and \(f_V\) is lci. By \cite[Corollary 3.3.7]{Jin}, the Gysin morphism \(f_V^!\) on \(\Kgl_{2d, d}^\bm\) identifies with pulling back coherent sheaves on \(G_0\). Thus \(f_V^![\mathscr{O}_\bb{Z}] \cong [\mathscr{O}_V]\) identifies with the fundamental class \(\eta_{f_V}^\Kgl\), up to shifting \(\mathrm{Th}(\bb{L}_{f_V}) \simeq \bb{T}^{\otimes - d}\) to the source. As \(\nu\) is a morphism of motivic spectra, \(\nu(\eta_{f_V}^\Kgl) = \eta_{f_V}^{H\bb{Z}}\), and the isomorphism 
\[H\bb{Z}_{2d, d}^\bm(V/\bb{Z}) \cong H_\bb{A}^0\big(V, \bb{Z}(0)\big) \cong \bb{Z}\]
is exactly the purity equivalence \ref{thm:fund-class}, so capping by the fundamental class \(\eta_{f_V}^{H\bb{Z}}\). Thus this maps \(\eta_{f_V}^{H\bb{Z}}\) to 1, completing the proof. 
\end{proof}

\begin{corollary}
For an integral closed subscheme \(C \subset T\) with \(\delta(C) = d\), 
\begin{equation}\label{eq:nu-int}
\nu_d[\mathscr{O}_C] = [C]
\end{equation}
\end{corollary}

\begin{proof}
The generic point \(x = \eta_C\) is a top-dimensional point of \(T\), so \(\mathscr{O}_{C, x}\) has length one over \(\mathscr{O}_{T, x}\), and \([C] \in H\bb{Z}_{2d, d}^\bm(T/\bb{Z})\) is by definition the class mapping to \(\eta_C\) under the isomorphism \(\gamma_T\).
\end{proof}

\begin{lemma}\label{lemm:gys-struct}
Fix \(U \in X_\nis\). Let \(\phi\) be one of the three structure morphisms in \ref{lemm:senh-struc}. Given \(W \in \mathcal{S}_n(U, q)\), let \(\tilde{\phi} \colon \tilde{W} \to W\) be the pullback of \(\phi\) along \(W \hookrightarrow U \times \Delta^q\). Under the isomorphism \(\pi_0 F_U^W \cong H\bb{Z}_{2d, d}^\bm(W/\bb{Z})\) of \ref{lemm:HZBM}, the morphism on \(F_U^W\) induced by \(\phi\) is the Gysin morphism 
\[\tilde{\phi}^! \colon H\bb{Z}_{2d, d}^\bm(W/\bb{Z}) \to H\bb{Z}_{2(d + r), d + r}^\bm(\tilde{W}/\bb{Z}),\]
with \(r = 0, 1\), and \(- 1\), in the case of \(\phi\) étale, a codegeneracy, and a coface, respectively.
\end{lemma}

\begin{proof}
The étale and codegeneracy cases immediately give smooth \(\phi\) of relative dimension \(0, 1\), respectively. Letting the coface inclusion be given by \(V(t)\), for some coordinate \(t = t_i\) in the algebraic simplex, this is a non-zerodivisor, so a regular closed embedding of codimension \(1\).

The morphisms given in \ref{lemm:senh-struc} are application of \(\phi^!\), then purity \(\phi^!\simeq \phi^* \otimes \bb{T}^{\otimes r}\), for \(r\) the relative dimension of \(\phi\), then application of support exchange. The fundamental class \(\eta_\phi\) is given by the purity equivalence in the case of \(\phi\) smooth, and in the case of \(\phi\) a regular closed embedding by \ref{thm:fund-class}. The Cartesian square 
\[\begin{tikzcd} [row sep = large, column sep = large]
\tilde{W} \ar[d, "\tilde{i}", swap] \ar[r, "\tilde{\phi}", ""{name=h1}] & W \ar[d, "i"]\\
U' \times \Delta^{q'} \ar[r, "\phi", swap, ""{name=h2}] & U \times \Delta^q,
\ar[phantom, from=h1, to=h2, "\Delta"]
\end{tikzcd}\]
is tor-independent in all three cases of \(\phi\): when étale, so \(U' = V\), \(q' = q\), and when a codegeneracy, so \(U' = U\), \(q' = q + 1\), \(\phi\) is flat. When a coface, so \(U' = U\), \(q' = q - 1\), \(W\) being reduced and meeting faces properly implies that \(t\) is a regular element on \(\mathscr{O}_W\), so that 
\[\begin{tikzcd}
0 \ar[r] & \mathscr{O}_W \ar[r, "t \cdot"] & \mathscr{O}_W \ar[r] & \mathscr{O}_{V(t)} \ar[r] & 0
\end{tikzcd}\]
is a resolution for the structure sheaf of the \(\mathscr{O}_{V(t)}\) on \(U \times \Delta^q\), whose exactness is preserved under tensoring by \(\mathscr{O}_W\). Support exchange is the base change equivalence \(\phi^! i_* \simeq \tilde{i}_* \tilde{\phi}^!\), and composing the fundamental class \(\eta_\phi\) with this gives the refined Gysin map \(\Delta^* \eta_\phi\). As the square is tor-independent, by transverse base change, \(\Delta^* \eta_\phi = \eta_{\tilde{\phi}}\), the Gysin morphism \(\tilde{\phi}^!\).
\end{proof}

\begin{theorem}[Structure of the spectral replacement]
Writing \(\gamma_{U, q} \colon \pi_0 \big(F_n(U, q)\big) \to z_{n + 1}(U, q)\) for the isomorphism in \ref{lemm:EM}, this upgrades to an equivalence of functors 
\[\gamma \colon \pi_0\big(F_n( -, \bullet)\big) \simeq z_{n + 1}( -, \bullet) \colon X_\nis^\op \times \Delta^\op \to \mathrm{Ab}\]
\end{theorem}

\begin{proof}
\(z_{n + 1}(U, \bullet)\) is a presheaf of simplicial Abelian groups, so proving that the isomorphism \(\gamma_{U, q}\) is functorial in Nisnevich opens \(U\) and simplicial coordinates \(q\) proves both that \(\pi_0\big(F_n( -, \bullet)\big)\) is a presheaf of simplicial Abelian groups, and that it is equivalent to \(z_{n + 1}( -, \bullet)\).

The map \(\pi_0 F_U^W \to \pi_0 F_n(U, q)\) is the inclusion of direct summands, so is compatible with the structure morphisms of \ref{lemm:senh-struc}. Take some integral \(W \in \mathcal{S}_n(U, q)\). We must show that \(\gamma_{\tilde{W}}\big(\tilde{\phi}^![W]\big) = \tilde{\phi}^*[W]\), in which \(\tilde{\phi}^*\) is the structure morphism on \(z_{n + 1}(- , \bullet)\). \(\nu \colon \kgl \to H\bb{Z}\) is a morphism of motivic spectra, and the functoriality of \(\Kgl^\bm\) being compatible with that of \(\G\)-theory \cite[Corollarly 3.3.7]{Jin} gives a commutative diagram
\[\begin{tikzcd} [row sep = large, column sep = large]
G_0(W) \ar[d, "\tilde{\phi}^*", swap] \ar[r, "\sim", swap] \ar[rrr, "\nu_d", bend left = 12] & \Kgl_{2d, d}^\bm(W/\bb{Z}) \ar[d, "\tilde{\phi}^*"] \ar[r, "\sim", swap] & \kgl_{2d, d}^\bm(W/\bb{Z}) \ar[d] \ar[r, "\nu", swap] & H\bb{Z}_{2d, d}^\bm(W/\bb{Z}) \ar[d, "\tilde{\phi}^!"]\\
G_0(\tilde{W}) \ar[rrr, "\nu_{d + r}", bend right = 12, swap] \ar[r, "\sim"] & \Kgl_{2(d + r), d + r}^\bm \ar[r, "\sim"] & \kgl_{2(d + r), d + r}^\bm(\tilde{W}/\bb{Z}) \ar[r, "\nu"] & H\bb{Z}_{2(d + r), d + r}^\bm(\tilde{W}/\bb{Z})
\end{tikzcd}\]
in which \(\tilde{\phi}^*\) is contravariant functoriality of \(\G\)-theory. In all cases of \(\tilde{\phi}\), tor-independence implies that \(\tilde{\phi}^*[\mathscr{O}_W] = [\mathscr{O}_{\tilde{W}}]\). So 
\[\tilde{\phi}^![W] \stackrel{\eqref{eq:nu-int}} = \tilde{\phi}^!\nu_d[\mathscr{O}_W] = \nu_{d + r}\big(\tilde{\phi}^*[\mathscr{O}_W]\big) = \nu_{d + r}[\mathscr{O}_{\tilde{W}}],\]
and application of \ref{prop:cyc-coh} thus gives 
\[\gamma_{\tilde{W}}\big(\nu_{d + r}[\mathscr{O}_{\tilde{W}}]\big) = \sum_{y \in \tilde{W}_{(d + r)}} \ell(\mathscr{O}_{\tilde{W}, y}\big) \cdot \eta_y,\]
noting that the length of the stalk over itself is just its length. In the case of \(\tilde{\phi}\) smooth, all lengths are 1 since \(W\) is reduced. The structure morphism on the Bloch-Levine side is simply flat pullback, mapping \([C]\) to \([\tilde{\phi}^{- 1}(C)]\), and the identification \(F_V^{\tilde{W}} \cong z_{n + 1}(V, q + r)\) maps \(\eta_y\) to \([Z(y)]\), so that these are compatible.

In the case that \(\tilde{\phi}\) is a coface inclusion, then the map of the cycle group is intersecting \([W]\) with the effective Cartier divisor \(V(t)\) \cite[\S 2.3]{Ful}, which is a linear combination of \([Z(y)]\) with coefficients the lengths above. Thus the structure maps are compatible in all cases of \(\phi\), so upon taking colimits, we arrive at the desired equivalence.
\end{proof}

\begin{theorem}[Equivalence A]\label{thm:equiv-A}
\(F_n( -, \bullet)\) is a presheaf of simplicial spectra, and for any \(j \in \bb{Z}\), there is an equivalence 
\begin{equation}
R\Gamma\big(X_\nis, |F_n(-, \bullet)|\big) \simeq R\Gamma \big(X_\mathrm{Zar}, \mathcal{Z}_{n + 1}(X, *)\big)
\end{equation}
\end{theorem}

\begin{proof}
By \ref{lemm:pi-nought}, the equivalence of simplicial presheaves of Abelian groups \(\gamma \colon \pi_0 \big(F_n( -, \bullet)\big) \simeq z_{n + 1}( -, \bullet)\) lifts to an equivalence of simplicial presheaves of (Eilenberg-Mac Lane) spectra,
\begin{equation}\label{eq:EMFn}
F_n( -, \bullet) \simeq Hz_{n + 1}( - , \bullet)
\end{equation}
giving the first result.
Taking geometric realizations, there is an equivalence of presheaves of spectra 
\[|F_n( -, \bullet)| \simeq |Hz_{n + 1}( -, \bullet)|\]
The geometric realization of a presheaf is computed objectwise, so fix some \(U \in X_\nis\). Then  
\[|H z_{n + 1}(U, \bullet)| \simeq H N_*z_{n + 1}(U, \bullet) \simeq H \mathcal{Z}_{n + 1}(X, *)(U),\]
by Dold-Kan and the definition of \(\mathcal{Z}_{n + 1}(X, *)(U)\) as the alternating complex of \(z_{n + 1}(U , \bullet)\). So we have an equivalence of presheaves of spectra, 
\[|F_n( -, \bullet)| \simeq |Hz_{n + 1}(- , \bullet)| \simeq H \mathcal{Z}_{n + 1}(X, *)\]
Forming derived global sections gives an equivalence  
\[R\Gamma(X_\nis, |F_n( -, \bullet)|) \cong R\Gamma \big(X_\nis, \mathcal{Z}_{n + 1}(X, *)\big)\]
\(X\) is essentially of finite type over the spectrum of a Dedekind ring, so by Geisser \cite[Proposition 3.6]{Gei} there is an equivalence 
\[R\Gamma\big(X_\nis, \mathcal{Z}_{n + 1}(X, *)\big) \simeq R\Gamma \big(X_\mathrm{Zar}, \mathcal{Z}_{n + 1}(X, *)\big),\]
completing the proof.
\end{proof}

\section{Slice Computation B}\label{slice-comp-B}
This section finishes the proof of representability of higher Chow groups, by proving that derived global sections of the spectral replacement \(|F_n( -, \bullet)|\) are equivalent to those of the motivic Borel-Moore complex. This is done by filtering by niveau, constructing a morphism of filtered presheaves of spectra, and proving that it is a Nisnevich-local equivalence on each graded piece.

\subsection{Filtrations and the filtered morphisms in the homotopy category}

In this subsection, we construct filtrations on the spectral replacement and motivic Borel-Moore homology spectrum, and describe a morphism between them at the level of objects and morphisms. This morphism is upgraded coherently in the next subsection, but this explanation helps give an intuitive explanation of the geometric content of the construction.

\begin{construction}[Filtration on \(F_n(U, q)\)]
For \(q \in \bb{N}\), write \(d = n + q\) as before. Given \(U \in X_\nis\), \(q \geq 0\), define the poset category 
\[\mathcal{S}_{\leq m}^n(U, q) = \left\{W \in \mathcal{S}_n(U, q) \, \big| \, \delta\big(\overline{\pr_1(W)}\big) \leq m\right\},\]
This is a filtered category, as \(\pr_1(W_1 \cup W_2) = \pr_1(W_1) \cup \pr_1(W_2)\), and \(\delta\) takes the dimension of the largest irreducible component, preserving the dimension bound under union. Define 
\[\fil_m F_n(U, q) = \colim_{W \in \mathcal{S}_{\leq m}^n(U, q)} F_U^W\]
The inclusion functor \(\mathcal{S}_{\leq m - 1}^n(U, q) \to \mathcal{S}_{\leq m}^n(U, q)\) is fully faithful, so gives rise to a canonical morphism of spectra 
\[\fil_{m - 1} F_n(U, q) \to \fil_m F_n(U, q),\]
so that \(\fil_\star F_n(U, q) \in \fil \, \Sp\). Note that this is a bounded filtration, as \(X\) is of finite Krull dimension.
\end{construction}

\begin{lemma}
This upgrades to a filtration on the simplicial presheaf of spectra \(F_n( -, \bullet)\).
\end{lemma}

\begin{proof}
We know that taking colimits over \(\mathcal{S}_n(U, q)\) is compatible with étale morphisms in \(U \in X_\nis\) and face and degeneracy morphisms in \(q \in \Delta\), so we must demonstrate compatibility with the dimension bound on closures of projection to the first factor with étale, face, and degeneracy morphisms. Take \(W \in \mathcal{S}_{\leq m}^n(U, q)\).

\begin{enumerate}

\item Étale morphisms: Given an étale morphism \(g \colon V \to U\) over \(X\) inducing étale \(V \times \Delta^q \to U \times \Delta^q\), the following diagram
\[\begin{tikzcd} [row sep = large, column sep = large] g^{- 1}(W) \ar[d] \ar[r] & W \ar[d]\\
V \times \Delta^q \ar[d, "\pr_1", swap] \ar[r, ""] & U \times \Delta^q \ar[d, "\pr_1"]\\
V \ar[r, "g", swap] & U \end{tikzcd}\]
commutes. if \(\delta(\overline{\pr_1(g^{- 1}(W))}) > m\), then there is a point \(v \in \overline{\pr_1(g^{- 1}(W))}\) with \(\delta(v) > m\). Writing \(u = g(v)\), \(\trd \big(\kappa(v) : \kappa(u)\big) = 0\) by g étale implies that \(\delta(u) = \delta(v) > m\), contradiction to our assumption that \(W\) lies in the specified level of the filtration. Thus 
\[\delta \! \left(\overline{\pr_1(g^{- 1}(W))}\right) \leq m,\]
so that \(g^{- 1}(W) \in \mathcal{S}_{\leq m}^n(V, q)\).

\item Degeneracy morphisms:
Let \(p \colon U \times \Delta^{q + 1} \to U \times \Delta^q\) be a codegeneracy morphism, which gives \(p^{- 1}(W) \in \mathcal{S}_n(U, q + 1)\). We have a commutative diagram 
\[\begin{tikzcd} [row sep = large, column sep = large] p^{- 1}(W) \ar[d] \ar[r] & W \ar[d]\\
U \times \Delta^{q + 1} \ar[dr, "\pr_1", swap] \ar[r, "p"] & U \times \Delta^q \ar[d, "\pr_1"]\\
& U \end{tikzcd}\]
the lower triangle commuting due to \(U \times \Delta^\bullet\) having the structure of a cosimplicial scheme. Thus 
\[\delta \! \left(\overline{\pr_1(p^{- 1}(W))}\right) \leq \delta \! \left(\overline{\pr_1(W)}\right) \leq m,\]
so that \(p^{- 1}(W) \in \mathcal{S}_{\leq m}^n(U, q + 1)\).

\item Face morphisms: letting \(\iota \colon U \times \Delta^{q - 1} \to U \times \Delta^q\) be a coface inclusion, a nearly identical argument to that above demonstrates that 
\[\delta \! \left(\overline{\pr_1(\iota^{- 1} W)}\right) \leq \delta \! \left(\overline{\pr_1(W)}\right) \leq m,\]
and thus \(\iota^{- 1} W \in \mathcal{S}_{\leq m}^n(U, q - 1)\).

\end{enumerate}
Thus these three classes of morphism respect formation of colimits over \(\mathcal{S}_{\leq m}^n(U, q)\), meaning \(\fil_\star F_n( -, \bullet) \in \fil \, \psh(X_\nis \times \Delta, \Sp)\).
\end{proof}

\begin{definition}[\(\mathcal{S}_{\leq m}\)]
For \(U \in X_\nis, m \geq - 1\), define \(\mathcal{S}_{\leq m}(U)\) to be the category of reduced closed subschemes \(Z \subset U\) such that \(\delta(Z) \leq m\). This is filtered since \(\delta\) is preserved under unions. 
\end{definition}

\begin{construction}[Filtration on \(\bb{Z}(n)^\bm \lbrack -2n \rbrack\)]\label{con:fil-ZBM}
Consider 
\[\bb{Z}(n)^\bm( -)\lbrack -2n \rbrack = \map_{\sh( -)}\big(\mathbbm{1}_{( -)}, f_{( -)}^! s^{- n} \Kgl_\bb{Z}\big) \in \psh(X_\nis, \Sp)\]
to be a constant simplicial presheaf of spectra. We equip this with a similar niveau filtration obtained from \(\delta\). For \(U \in X_\nis\), let 
\[\fil_m \, \bb{Z}(n)^\bm(U)\lbrack -2n \rbrack = \colim_{Z \in \mathcal{S}_{\leq m}(U)} \map_{\sh(U)}\big(\mathbbm{1}_U, \Gamma_Z f_U^! s^{- n}\Kgl_\bb{Z}\big),\]
Formation of sections with support is functorial by \ref{lemm:supp-fun}, and respects the presheaf structure of \(\bb{Z}(n)^\bm[2n]\) by reasoning analogous to \ref{lemm:senh-struc} (iii), only using six-functors and support exchange, which are fully coherent. 
\end{construction}

\begin{remark}[Projections, pullbacks, and filtrations]\label{rem:proj-pull-fil}
For \(W \in \mathcal{S}_{\leq m}^n(U, q)\), denote \(Z_W = \overline{\pr_1(W)}\), a closed subscheme of \(U\). The filtration condition on \(W\) is equivalent to \(Z_W \in \mathcal{S}_{\leq m}(U)\). Let \(i \colon Z_W \hookrightarrow U\) and \(i_W \colon W \hookrightarrow U \times \Delta^q\), 
denote the closed embeddings, which via universal property, produce the commutative diagram 
\[\begin{tikzcd} [row sep = large, column sep = large] W \ar[ddr, "i_W", bend right = 15, swap] \ar[dr, "\iota"] \ar[drr, "\pr_1", bend left = 15]\\
& Z_W \times \Delta^q \ar[d, "\tilde{i}", swap] \ar[r, "\pr_1"] & Z_W \ar[d, "i"]\\
& U \times \Delta^q \ar[r, "\pr_1", swap] & U
\end{tikzcd}\]
in which the square is a pullback. \(\iota\) is a morphism \(W \to Z_W \times \Delta^q\), giving a natural transformation
\[\beta^\iota \colon \Gamma_W \Longrightarrow \Gamma_{Z_W \times \Delta^q}\]
By \ref{rem:tau-inv}, there is a natural equivalence\footnote{\(\tilde{\tau}^{Z_W}\) is too clunky, so we slightly abuse notation here} 
\[\tilde{\tau}^W \colon \Gamma_{Z_W \times \Delta^q} \, \pr_1^! \stackrel{\sim} \Longrightarrow \pr_1^! \Gamma_{Z_W} \colon \sh(U) \to \sh(U \times \Delta^q),\]
such that \(\pr_1^!\epsilon^i \circ \tilde{\tau}^W = \epsilon^{\tilde{i}} \pr_1^!\).
\end{remark}

\begin{construction}[The colimitand morphism \(\phi_m^q(W)\)]
Use the same notation as \ref{rem:proj-pull-fil}. The natural transformation \(\beta^\iota\) induces a morphism of spectra, 
\[F_U^W \to F_U^{Z_W \times \Delta^q}\footnote{Note that \(Z_W \times \Delta^q\) is almost never an object of \(\mathcal{S}_{\leq m}^n(U, q)\), but this mapping spectrum does not appear in the colimit defining the filtration on the spectral replacement. We merely use it as a go-between.}\]
As \(f_{U, q}^! = \pr_1^! \circ f_U^!\), support exchange \(\tau^{\tilde{W}}\) gives an equivalence 
\[\Gamma_{Z_W \times \Delta^q} f_{U, q}^! s^{- d} \Kgl_\bb{Z} \simeq \pr_1^! \Gamma_{Z_W} f_U^! s^{- d} \Kgl_\bb{Z}\]
Smooth purity gives \(\pr_1^! \simeq \pr_1^* \otimes \bb{T}^{\otimes q}\), and \(\bb{A}^1\)-invariance of \(\pr_1\) means \((\pr_1)_\sharp \mathbbm{1}_{U \times \Delta^q} \simeq \mathbbm{1}_U\). Recalling that \(- d + q = - n\), smashing with the Tate motive and adjointness give an equivalence 
\[F_U^{Z_W \times \Delta^q} \simeq \map_{\sh(U)}\big(\mathbbm{1}_U, \Gamma_{Z_W} f_U^! s^{- n} \Kgl_\bb{Z}\big)\]
This mapping spectrum is a colimitand of \(\fil_m \bb{Z}(n)^\bm[2n](U)\), since \(Z_W \in \mathcal{S}_{\leq m}(U)\) by assumption.  Denote by \(\phi_m^q(W)\) the composite of this map, which sends a colimitand of \(\fil_m F_n(U, q)\) to a colimitand of \(\fil_m \bb{Z}(n)^\bm[2n](U)\).
\end{construction}

To obtain a filtered morphism of the two simplicial presheaves of spectra, we want to take colimits, but doing so requires justifying compatibility.

\begin{lemma}\label{lemm:colim-leg-comm}
Given a morphism \(\begin{tikzcd} W_1 \ar[rr, "\iota_1", bend right = 17, swap] \ar[r, "\iota_{12}"] & W_2 \ar[r, "\iota_2"] & U \times \Delta^q \end{tikzcd}  \! \! \) in \(\mathcal{S}_{\leq m}^n(U, q)\), the diagram 
\begin{equation}\label{eq:fil-leg-comm}
\begin{tikzcd} [column sep = large, row sep = large] F_U^{W_1} \ar[d] \ar[r, "\Phi_m^q(W_1)"] & \map_{\sh(U)}\big(\mathbbm{1}_U, \Gamma_{Z_{W_1}} f_U^! s^{- n}\Kgl_\bb{Z}\big) \ar[d]\\
F_U^{W_2} \ar[r, "\Phi_m^q(W_2)", swap] & \map_{\sh(U)}\big(\mathbbm{1}_U, \Gamma_{Z_{W_2}} f_U^! s^{- n} \Kgl_\bb{Z}\big) \end{tikzcd}
\end{equation}
commutes, the right vertical morphism induced by the natural transformation\ \(\beta^{i_{12}} \colon \Gamma_{Z_{W_1}} \Longrightarrow \Gamma_{Z_{W_2}}\). 
\end{lemma}

\begin{proof}
We summarize the numerous closed embeddings and their relations via the following commutative diagram, 
\begin{equation}\label{eq:big-chungus-diagram}
\begin{tikzcd} W_1 \ar[dddr, "\iota_{W_1}", swap] \ar[dr, "\iota_{12}"] \ar[rr, "\iota_1"] & & Z_{W_1} \times \Delta^q \ar[dddl, "\tilde{i}_1"] \ar[dr, "\tilde{i}_{12}"] \ar[rr, "\pr_1"] & & Z_{W_1} \ar[dddl, "i_1"] \ar[dr, "i_{12}"]\\
& W_2 \ar[dd, "\iota_{W_2}"] \ar[rr, "\iota_2", crossing over] & & Z_{W_2} \ar[ddll, "\tilde{i}_2"] \times \Delta^q \ar[rr, "\pr_1", crossing over] & & Z_{W_2} \ar[ddll, "i_2"]\\
\\
& U \times \Delta^q \ar[rr, "\pr_1", swap] & & U\end{tikzcd}
\end{equation}
with top right squares of Cartesian. Use the abbreviated notation \(\tilde{\tau}^1 = \tilde{\tau}^{Z_{W_1}}, \tilde{\tau}^2 = \tilde{\tau}^{Z_{W_2}}\). Smashing by the Tate motive is an equivalence, as is the adjunction \((\pr_1)_\sharp \dashv \pr_1^*\), so this comes down to proving that \(\pr_1^! \beta^{i_{12}} \circ \tilde{\tau}^1 \circ \beta^{\iota_1}\pr_1^! \simeq \tilde{\tau}^2 \circ \beta^{\iota_2} \pr_1^! \circ \beta^{\iota_{12}}\pr_1^! \) in 
\[\Map_{\fun\big(\sh(U), \sh(U \times \Delta^q)\big)}\big(\Gamma_{W_1}\pr_1^!, \pr_1^! \Gamma_{Z_{W_2}}\big)\]
\(\pr_1^! \Gamma_{Z_{W_2}}\) has values in \(\sh_{Z_{W_2} \times \Delta^q}(U \times \Delta^q)\), with counit \(\pr_1^!\epsilon^{i_2}\).
\[\pr_1^! \epsilon^{i_2} \circ \left(\pr_1^! \beta^{i_{12}} \circ \tilde{\tau}^1 \circ \beta^{\iota_1}\pr_1^!\right) \simeq \pr_1^!\epsilon^{i_1} \circ \tilde{\tau}^1 \circ \beta^{\iota_1}\pr_1^! \stackrel{\eqref{eq:twid-tau-epsilon}} \simeq \epsilon^{\tilde{i}_1}\pr_1^! \circ \beta^{\iota_1}\pr_1^!\]
\[\simeq \epsilon^{\tilde{i}_2}\pr_1^! \circ \beta^{\tilde{i}_{12}}\pr_1^! \circ \beta^{\iota_1}\pr_1^! \simeq \epsilon^{\tilde{i}_2}\pr_1^! \circ \beta^{\iota_2 \circ \iota_{12}}\pr_1^! \stackrel{\eqref{eq:twid-tau-epsilon}} \simeq \pr_1^! \epsilon^{i_2} \circ \left(\tilde{\tau}^2 \circ \beta^{\iota_2}\pr_1^! \circ \beta^{\iota_{12}}\pr_1^!\right),\]
so by \ref{lemm:fun-fact}, these are equivalent functors.
\end{proof}

\begin{construction}[The colimit morphism \(\Phi_m^q(U)\)]\label{con:Phi-naive}
By \ref{lemm:colim-leg-comm}, \(\phi_m^q( - )\) takes a coconce of the functor on \(\mathcal{S}_{\leq m}^n(U, q) \to \Sp\) to a cocone of the functor \(\mathcal{S}_{\leq m}^n(U, q)\). As a result, we have morphism 
\[\Phi_m^q(U) \colon \colim_{W \in \mathcal{S}_{\leq m}^n(U, q)} F_U^W \to \colim_{Z \subset \mathcal{S}_{\leq m}(U)} \map_{\sh(U)}\big(\mathbbm{1}_U,\Gamma_Z f_U^! s^{- n} \Kgl_\bb{Z}\big),\]
that is, a morphism \(\Phi_m^q(U) \colon \fil_m F_n(U, q) \to \fil_m \bb{Z}(n)^\bm(U)\lbrack -2n\rbrack\), which is well-defined at the level of the homotopy category.
\end{construction}

\subsection{Extended comparison morphism}

To ensure that this is highly coherent, and not just compatible at the level of morphisms, we must extend \(\Phi_\star^q(U)\) to \(\Phi_\star \colon \fil_\star \, F_n( -, \bullet) \to \fil_\star \, \bb{Z}(n)^\bm\). We mimic the constructions of the previous section, but with the full machinery of \(\infty\)-functors.

\begin{definition}
Let \(\mathrm{Pair}_X\) be the the category of \((Y, W)\), for \(Y \in \sm_X\) and \(W \subset Y\) closed. A morphism \((Y_1, W_1) \to (Y_2, W_2)\) is a smooth morphism \(h \colon Y_1 \to Y_2\) over \(X\) such that \(h^{- 1}(W_2) \subset W_1\). 
\[\mathrm{Pair}_X \to \mathrm{Arr}(\sm_X), \quad (Y, W) \mapsto W - Y \subset Y,\]
is a functor to the arrow category of \(\sm_X\). Applying \(M_X\) and taking cofibres gives a functor 
\[M_X(-/-) \colon \mathrm{Pair}_X \to \sh(X), \quad (Y, W) \mapsto \mathrm{cofib}\big(M_X(Y - W) \to M_X(Y)\big)\]
We call \(M_X(Y/W)\) the \emph{motive of the pair \((Y, W)\)}.
\end{definition}

\begin{lemma}\label{lemm:mot-pair-eq}
For \(Y = U \times \Delta^q\), writing \(\mf{G}_X = f^! s^{- n} \Kgl_\bb{Z}\), then 
\begin{equation}\label{eq:mot-pair-map}
\map_{\sh(X)}(M_X(Y/W), \mf{G}_X) \simeq F_U^W
\end{equation}
In the case of \(q = 0\), with \(Z \subset U\), we thus have that 
\[\map_{\sh(X)}\big(M_X(U/Z), \mf{G}_X) \simeq \map_{\sh(U)}\big(\mathbbm{1}_U, \Gamma_Z f_U^! s^{- n} \Kgl_\bb{Z}\big)\]
\end{lemma}

\begin{proof}
Let \(f_{U, q} \colon Y = U \times \Delta^q \to \spec \bb{Z}\) denote the structure morphism, with \(i \colon W \hookrightarrow Y\) the closed embedding with complementary open \(j \colon Y - W \hookrightarrow Y\). Applying \(\map_{\sh(X)}( -, \mf{G}_X)\) is exact, and \(Y - W\) and \(Y\) are smooth over \(X\), so the mapping spectrum in question is the cofibre of 
\[\map_{\sh(X)}\big((f_{Y - W})_\sharp \mathbbm{1}_X, \mf{G}_X\big) \to \map_{\sh(X)}\big((f_{U, q})_\sharp \mathbbm{1}_X, \mf{G}_X\big),\]
by adjointness the cofibre of
\[\map_{\sh(U \times \Delta^q)}\big(\mathbbm{1}_{U \times \Delta^q}, j_*j^* f_{U, q}^* \mf{G}_X\big) \to \map_{\sh(U \times \Delta^q)}\big(\mathbbm{1}_{U \times \Delta^q}, f_{U, q}^* \mf{G}_X\big),\]
which the localization sequence is  
\[\map_{\sh(U \times \Delta^q)}\big(\mathbbm{1}_{U \times \Delta^q}, \Gamma_W f_{U, q}^* \mf{G}_X\big) \simeq \map_{\sh(U \times \Delta^q)}\big(\mathbbm{1}_{U \times \Delta^q}, \Gamma_W f_{U, q}^! s^{- d} \Kgl_\bb{Z}\big) = F_U^W,\]
the final equivalence by smooth purity.
\end{proof}

\begin{definition}[Filtration on \(\bb{Z}(n)^\bm \lbrack - 2n \rbrack\)]
Define 
\[\mf{G}_X^{\leq m} = \colim_{Z \in \mathcal{S}_{\leq m}(X)} \underline{\Hom}\big(M_X(X/Z), \mf{G}_X\big)\]
This object is functorial in the filtration index \(m\). Define 
\[\fil_m \bb{Z}(n)^\bm(U)\lbrack -2n \rbrack = \map_{\sh(X)}\big(M_X(U), \mf{G}_X^{\leq m}\big)\]
This assignment is an object of \(\fil \,\psh(X_\nis, \Sp)\), and satisfies 
\[\bb{Z}(n)^\bm(U)\lbrack -2n \rbrack \simeq \map_{\sh(X)}\big(M_X(U), \mf{G}_X\big),\]
so as \(\mf{G}_X^{\leq m} = \mf{G}_X\) for \(m \gg 0\), this gives a bounded filtration on (shifted) motivic Borel-Moore homology.
\end{definition}

\begin{lemma}
For each \(U \in X_\nis\), \(m \in \bb{Z}\), there is an equivalence
\[\fil_m \bb{Z}(n)^\bm(U) \lbrack - 2n \rbrack \simeq \colim_{Z \in \mathcal{S}_{\leq m}(U)} \map_{\sh(U)}\big(\mathbbm{1}_U, \Gamma_Z f_U^! s^{- n} \Kgl_\bb{Z}\big),\]
identifying this coherent filtration with the naive filtration of the previous section.
\end{lemma}

\begin{proof}
Compactness of \(M_X(U)\) and adjointness of internal hom gives 
\[\fil_m \bb{Z}(n)^\bm(U)\lbrack -2n \rbrack = \map_{\sh(X)}\! \left(M_X(U), \colim_{Z \in \mathcal{S}_{\leq m}(X)} \underline{\Hom}\big(M_X(X/Z), \mf{G}_X\big)\right)\]
\[\simeq \colim_{Z \in \mathcal{S}_{\leq m}(X)} \map_{\sh(X)}\big(M_X(U) \otimes M_X(X/Z), \mf{G}_X\big)\]
\(M_X\) is a symmetric monoidal functor, so the smash product is the cofibre of 
\[M_X(U \times_X (X - Z)) \to M_X(U \times_X X),\]
which is by definition \(M_X(U/U \cap Z)\). \(\mathcal{S}_{\leq m}(X) \to \mathcal{S}_{\leq m}(U)\), \(Z \mapsto Z \cap U\), is cofinal, so the limit is equivalent to 
\[\colim_{Z \in \mathcal{S}_{\leq m}(U)} \map_{\sh(X)}\big(M_X(U/Z), \mf{G}_X\big),\]
which is the desired colimit by \ref{lemm:mot-pair-eq}.
\end{proof}

\begin{definition}[Base categories]
Write \(\mathcal{B} = X_\nis^\op \times \Delta^\op \times (\bb{Z}, \leq)\), so that a functor \(\mathcal{B} \to \Sp\) is a filtered simplicial presheaf. A morphism in \(\mathcal{B}\) is a triple 
\[(g, \theta, \leq) \colon (U, q, m) \to (V, p, m'),\]
with \(g \colon V \to U\) étale over \(X\), \(\theta \colon p \to q\) in \(\Delta\), and \(m \leq m'\). We also denote by \(\theta\) the induced morphism of algebraic simplices, \(\theta \colon \Delta^p \to \Delta^q\), and write \(g_\theta = g \times \theta \colon V \times \Delta^p \to U \times \Delta^q\).

Similarly, let \(\mathcal{B}' = X_\nis^\op \times (\bb{Z}, \leq)\), with \(\psi \colon \mathcal{B} \to \mathcal{B}'\) the forgetful functor. \(\psi^*\) converts a filtered presheaf of spectra to a filtered presheaf of constant simplicial spectra.
\end{definition}

\begin{definition}[Fibre posets]\label{def:fib-poset}
Given \((U, q, m) \in \mathcal{B}\), let \(\mathcal{F}(U, q, m)\) be the poset of pairs \((W, Z)\) such that \(W \in \mathcal{S}_{\leq m}^n(U, q)\), \(Z \in \mathcal{S}_{\leq m}(U)\), and \((g_U \circ \pr_1)(W) \subset Z\), for \(g_U \colon U \to X\) the étale morphism encoding \(U \in X_\nis\). The assignment \((U, q, m) \mapsto \mathcal{F}(U, q, m)\), sending morphism \((g, \theta, \leq)\) to the morphism 
\[(W, Z) \mapsto \big(g_\theta^{- 1}(W), Z)\big)\]
is a functor from \(\mathcal{B}\) to \(\mathrm{Cat}\), since the scheme-theoretic inverse image of a closed subscheme strictly commutes with composition.
\end{definition}

\begin{construction}
Let \(\pi \colon \mathcal{T} \to \mathcal{B}\) be the Grothendieck construction of \(\mathcal{F} \colon \mathcal{B} \to \mathrm{Cat}\).
Writing \(B = (U, q, m) \in \mathcal{B}\), a morphism \((B, W, Z) \to (B', W', Z')\) in \(\mathcal{T}\) is a morphism \(B \to B'\) in \(\mathcal{B}\) such that \(g_\theta^{- 1}(W) \subset W'\) and \(Z \subset Z'\).
\end{construction}

\begin{lemma}\label{lemm:colim-prep}
\(\pi\) is coCartesian, with contractible fibres \(\mathcal{F}(B)\), and the forgetful functor \(\mathcal{F}(B) \to \mathcal{S}_{\leq m}^n(U, q)\) has left adjoint \(W \mapsto (W, Z_W)\), so is cofinal.
\end{lemma}

\begin{proof}
The fibre \(\mathcal{F}(B)\) contains \((\varnothing, \varnothing)\), and is closed under taking unions, so is contractible. And \((W, Z_W) \leq (W', Z')\) exactly when \(W \subset W'\), as \(Z_W \subset Z_{W'} \subset Z'\).
\end{proof}

\begin{construction}[Coherent filtration on spectral replacement]
Write \(F \colon \mathcal{T} \to \Sp\) as the composition 
\[\begin{tikzcd} [column sep = huge]
\mathcal{T} \ar[r] & \mathrm{Pair}_X^\op \ar[r, "{M_X(-/-)}"] & \sh(X)^\op \ar[r, "{\map(\!-,\mathfrak{G}_X)}"] & \Sp
\end{tikzcd}\]
the first functor mapping \((U, q, m, W, Z)\) to \((U \times \Delta^q, W)\). Define 
\[\fil_\star F_n( -, \bullet) = \lan_\pi F \in \fun(\mathcal{B}, \Sp) = \fil \, \psh(X_\nis \times \Delta, \Sp)\]
\end{construction}

\begin{lemma}
For \((U, q, m) \in \mathcal{B}\), there is an equivalence 
\[\lan_\pi F(U, q, m) \simeq \colim_{W \in \mathcal{S}_{\leq m}^n(U)} F_U^W,\]
identifying the left Kan extension as the filtration constructed in the previous section.
\end{lemma}

\begin{proof}
\ref{lemm:colim-prep} tells us that left Kan extensions along \(\pi\) are computed as colimits over fibres since it is a coCartesian fibration, and also that we may compute this colimit over \(\mathcal{S}_{\leq m}^n(U, q)\) by cofinality. 
\[F(B, W, Z) = \map_{\sh(U \times \Delta^q)}\big(M_X(Y/W), \mf{G}_X\big) \stackrel{\ref{eq:mot-pair-map}} \simeq F_U^W,\]
so we have the desired result.
\end{proof}

\begin{construction}\label{con:phi}
Given \((B, W, Z) \in \mathcal{T}\), write \(Y = U \times \Delta^q\) for \(B = (U, q, m)\) as before, and write \(\phi(B, W, Z)\) for the composite 
\[\begin{tikzcd} [row sep = large]
F(B, W, Z) = \map_{\sh(X)}\big(M_X(Y/W), \mf{G}_X\big) \ar[r, "(i)"] & \map_{\sh(X)}\big(M_X(Y/(g_U \circ \pr_1)^{- 1}Z), \mf{G}_X\big) \ar[dl, "(ii)", swap]\\ \map_{\sh(X)}\big(M_X(Y), \underline{\Hom}\big(M_X(X/Z), \mf{G}_X\big)\big) \ar[r, "(iii)"] & \map_{\sh(X)}\big(M_X(Y), \mf{G}_X^{\leq m}\big) \ar[dl, "(iv)"]\\
\map_{\sh(X)}\big(M_X(U), \mf{G}_X^{\leq m}\big) = \pi^* \psi^* \fil_m \bb{Z}(n)^\bm(U) \lbrack -2n \rbrack
\end{tikzcd}\]
The connecting morphisms are given as follows:
\begin{enumerate}[(i)]

\item Functoriality of \(M_X( -/-)\) applied to the morphism \((Y, W) \to \big(Y, (g_U \circ \pr_1)^{- 1}(Z)\big)\) in \(\mathrm{Pair}_X\) (which exists since \(Z_W \subset Z\)) in the source of the mapping spectrum.

\item Symmetric monoidality of \(M_X\) in the source of the mapping spectrum, then adjoint equivalence of smash product and inner hom.

\item The canonical inclusion into the colimit 
\[\underline{\Hom}\big(M_X(X/Z), \mf{G}_X\big) \to \colim_{Z \in \mathcal{S}_{\leq m}(X)} \underline{\Hom}\big(M_X(X/Z), \mf{G}_X\big) = \mf{G}_X^{\leq m}\]
in the target of the mapping spectrum

\item The inverse of the \(\bb{A}^1\)-invariance equivalence \(\pr_1 \colon M_X(U \times \Delta^q) \simeq M_X(U)\) applied to the source of the mapping spectrum.
\end{enumerate}
Each step is natural, relying on strict commutativity of \(\mathcal{T} \to \mathrm{Pair}_X\), functoriality of \(M_X\), adjointness, or universal properties of colimits. 
As a result, \(\phi\) is a natural transformation of functors 
\[F \Longrightarrow \pi^* \psi^* \bb{Z}(n)^\bm \lbrack -2n \rbrack \in \fil \, \psh(X_\nis \times \Delta)\]
\end{construction}

\begin{remark}
This \(\phi\) is the familiar morphism on the colimitands, refined to a natural transformation in a coherent setting. We will further refine the morphism on the colimits themselves through left Kan extensions.
\end{remark}

\begin{definition}[Filtered morphism \(\Phi_\star\)]
Define the morphism of filtered presheaves of simplicial spectra \(\Phi_\star \colon \fil_\star F_n( -, \bullet) \to \psi^* \bb{Z}(n)^\bm [-2n]\) as 
\[\Phi_\star = \lan_\pi \phi \colon \fil_\star F_n( -, \bullet) \simeq  \lan_\pi F \to \lan_\pi \pi^* \psi^* \bb{Z}(n)^\bm \lbrack -2n \rbrack \simeq \psi^* \bb{Z}(n)^\bm \lbrack -2n \rbrack,\]
the final morphism given by the counit of the adjunction \(\lan_\pi \dashv \pi^*\), which is an equivalence since \(\pi\) has weakly contractible fibres. Note that using the colimit description of left Kan extensions on objects, this is precisely the morphism built at the level of the homotopy category in \ref{con:Phi-naive}.
\end{definition}



\subsection{Filtrations and realizations}
To convert the morphism of filtrations in \(\fil \, \psh(X_\nis \times \Delta, \Sp)\) to one in \(\fil \, \psh(X_\nis, \Sp)\), we take geometric realizations.

\begin{definition}[The geometrically realized morphism \(|\Phi_\star|\)]
The morphism \(\psi \colon \mathcal{B} \to \mathcal{B}'\) is a coCartesian fibration with fibre \(\Delta^\op\), so left Kan extensions along \(\psi\) are computed objectwise. So 
\[\lan_\psi \fil_\star F_n( -, \bullet)(U, q, m) \simeq \colim_{q \in \Delta^\op} \fil_m F_n(U, q) = |\fil_m F_n(U, \bullet)|,\]
identifying this left Kan extension with the geometric realization of the spectral replacement.

\(\Delta^\op\) is a sifted category, thus weakly contractible, meaning once again the counit of \(\lan_\psi \dashv \psi\) is an equivalence, so that 
\[\lan_\psi \psi^*\fil_\star \bb{Z}(n)^\bm \lbrack -2n \rbrack \simeq \fil_\star \bb{Z}(n)^\bm \lbrack - 2n\rbrack\]
Write \(|\Phi|_\star = \lan_\psi \, \Phi_\star\). This is a morphism in \(\fil \, \psh(X_\nis, \Sp)\), 
\[|\Phi|_\star \colon \fil_\star |F_n( -, \bullet)| \to \fil_\star \, \bb{Z}(n)^\bm\lbrack -2n \rbrack\]
Left Kan extension is a left adjoint, thus preserves colimits and is exact when its source and target are stable \(\infty\)-categories. Thus this morphism has associated gradeds 
\[\gr_m |\Phi|_\star \simeq |\gr_m \Phi_\star| \colon \gr_m |F_n( -, \bullet)| \simeq |\gr_m F_n( -, \bullet)| \to \gr_m \bb{Z}(n)^\bm \lbrack -2n \rbrack\]
for each \(m \in \bb{Z}\).
\end{definition}

Proving that \(|F_n( -, \bullet)|\) is equivalent to \(\bb{Z}(n)^\bm\lbrack -2n \rbrack\) is, as far as the author is aware, reliant on moving lemmas in a far more general setting than available. However, proving Nisnevich-local equivalence is possible, which allows for an equivalence of derived global sections. We further convert the filtered morphism to one on the derived global sections of the relevant sheaves of spectra.

\begin{lemma}[Descent and hypersheaves]
A presheaf \(\mathcal{F} \in \psh(X_\nis, \Sp)\) is Nisnevich-local if and only if \(\mathcal{F}(\varnothing) = 0\) and \(\mathcal{F}\) maps an elementary distinguished square to a Cartesian square in \(\Sp\).
\end{lemma}

\begin{proof}
By \cite[Proposition 2.8.1 (1)]{BEM}, \(\shv(X_\nis, \Sp)\) are exactly \(\mathcal{F} \in \psh(X_\nis, \Sp)\) satisfying the conditions above. \(X\) is of finite Krull dimension, so by \cite[Theorem 1.7]{CM}, any sheaf of spectra is a hypersheaf, that is, local. And any hypersheaf is a sheaf, so we have the desired equivalence of conditions.
\end{proof}

\begin{corollary}\label{corr:loc-znbm}
\(\fil_m \bb{Z}(n)^\bm \lbrack -2n \rbrack\) is a Nisnevich sheaf of spectra.
\end{corollary}

\begin{proof}
\(M_X\) maps elementary distinguished squares of schemes to Cartesian squares of spectra, and \(\map_{\sh(X)}\big( -, \mf{G}_X^{\leq m}\big)\) maps these to coCartesian, equivalently by stability, Cartesian, squares. Likewise, \(M_X\) maps \(\varnothing\) to \(0\), so the mapping spectrum functor does as well.
\end{proof}

We use this to see that derived global sections of (shifted) motivic Borel-Moore homology are simply obtained by evaluation.

\begin{lemma}
There are filtrations \(\fil_\star R\Gamma(X_\nis, |F_n( -, \bullet)|), \fil_\star \, \bb{Z}(n)^\bm(X)\lbrack - 2n \rbrack\) with associated gradeds 
\[\gr_m R\Gamma(X_\nis, |F_n( -, \bullet)|) \simeq R\Gamma \big(X_\nis, \gr_m |F_n( -, \bullet)|\big) \simeq R\Gamma \big(X_\nis, \big|\gr_m \, F_n( -, \bullet)\big|\big)\]
\end{lemma}

\begin{proof}
The functor \(R\Gamma(X_\nis, -)\) factors as 
\[\begin{tikzcd}
\psh(X_\nis, \Sp) \ar[r, "L_\nis"] & \shv(X_\nis, \Sp) \ar[r, "\mathrm{ev}_X"] & \Sp
\end{tikzcd}\]
Nisnevich localization \(L_\nis\), is given by left Bousfield localization, hence is a left adjoint, preserving colimits. Its target and source are stable categories, so it is exact. Limits of sheaves are computed objectwise, so \(\mathrm{ev}_X\) preserves limits, and is thus exact as well.
\end{proof}

\begin{theorem}[Equivalence B]\label{thm:equiv-B}
If \(\gr_m |\Phi| \colon \gr_m |F_n( -, \bullet)| \to \gr_m \, \bb{Z}(n)^\bm \lbrack -2n \rbrack\) is a local equivalence for each \(m \in \bb{Z}\), then there is an equivalence 
\[R\Gamma \big(X_\nis, |F_n( -, \bullet)|\big) \simeq \bb{Z}(n)^\bm(X)\lbrack -2n \rbrack\]
\end{theorem}

\begin{proof}
We prove that each \(R\Gamma(|\Phi_m|)\) is an equivalence using induction on \(m\). For \(m \leq - 2\), \(R\Gamma(|\Phi_m|)\) is simply self-equivalence of \(0\), so the base case holds. Assuming that this is true for some \(m \geq - 2\), we have a commutative diagram of spectra
\[\begin{tikzcd} [row sep = huge, column sep = large]
\fil_{m - 1} R\Gamma\big(X_\nis, |F_n( -, \bullet)|\big) \ar[d, "R\Gamma(\lvert\Phi_{m - 1}\rvert)", swap] \ar[r] & \fil_m R\Gamma\big(X_\nis, |F_n( -, \bullet)|\big) \ar[d, "R\Gamma(\lvert\Phi_m\rvert)"] \ar[r] & \gr_m R\Gamma\big(X_\nis, F_n( -, \bullet)|\big) \ar[d, "R\Gamma(\gr_m \lvert\Phi_\star\rvert)"]\\
\fil_{m - 1} \, \bb{Z}(n)^\bm(X)\lbrack -2n \rbrack \ar[r] & \fil_m \, \bb{Z}(n)^\bm(X)\lbrack -2n \rbrack \ar[r] & \gr_m \, \bb{Z}(n)^\bm(X)\lbrack -2n \rbrack, \end{tikzcd}\]
where we have simplified the bottom row by virtue of \ref{corr:loc-znbm}.
The left vertical leg is an equivalence by the inductive hypothesis, and the right vertical leg is an equivalence since derived global sections map local equivalences to equivalences. As the rows are fibre sequences in a stable category, the middle vertical leg must be an equivalence as well, completing the induction. As both objects are bounded,
\[\fil_m R\Gamma\big(X_\nis, |F_n( -, \bullet)|\big) \simeq R\Gamma\big(X_\nis, |F_n( -, \bullet)|\big), \; \fil_m \, \bb{Z}(n)^\bm(X)\lbrack -2n \rbrack \simeq \bb{Z}(n)^\bm(X)\lbrack -2n \rbrack,\]
for \(m \gg 0\), giving the desired result.
\end{proof}

\subsection{Graded equivalence: shells}
In this section, we break down the spaces and morphisms between them into simpler pieces. This is the first half of proving that \(\gr_m |\Phi|\) is a Nisnevich-local equivalence. 

\begin{remark}\label{rem:cyc-in-closed}
Given \(U \in X_\nis\) and \(q \in \Delta\), if \(Z \subset U\) is closed, then 
\[\mathcal{S}_n(Z, q) = \{W \in \mathcal{S}_n(U, q) \, | \, W \subset Z \times \Delta^q\},\]
as \(W \cap (Z \times \Delta^{q - 1}) = W \cap (U \times \Delta^{q - 1})\), so that meeting faces of \(Z \times \Delta^q\) properly is equivalent to meeting faces of \(U \times \Delta^q\) properly. Thus there is an inclusion \(\mathcal{S}_n(Z, q) \to \mathcal{S}_n(U, q)\) of posets.
\end{remark}

\begin{lemma}[Cycles on closed subschemes]\label{lemm:Z-S-fil}
Let \(Z \subset U\) be closed.
\begin{enumerate}

\item Then the inclusion \(\mathcal{S}_n(Z, q) \to \mathcal{S}_n(U, q)\) induces a pushforward \(F_n(Z, q) \to F_n(U, q)\) compatible with simplicial structures, which on \(\pi_0\) is \(z_{n + 1}(Z, q) \subset z_{n + 1}(U, q)\) under the identification \ref{lemm:EM}.

\item For any \(m \geq - 1\), \(\fil_m F_n(U, q)\) is an Eilenberg-Mac Lane spectrum, with \(\pi_0 \fil_m F_n(U, q) \subset z_{n + 1}(U, q)\) generated by the classes \([C]\) such that \(\delta\big(\pr_1(\eta_C)\big) \leq m\).

\item Taking colimits over pushforwards gives an equivalence 
\[\begin{tikzcd}
\displaystyle \colim_{Z \in \mathcal{S}_{\leq m}(U)} F_n(Z, q) \ar[r, "\sim"] & \fil_m F_n(U, q)
\end{tikzcd}\]
\end{enumerate}
\end{lemma}

\begin{proof}
(1) follows from the fact that the colimitand \(F_U^W\) only depends on the scheme \(W\), and that the coface and codegeneracy morphisms of \(Z \times \Delta^\bullet\) are inherited from \(U \times \Delta^\bullet\). For (2), \(\fil_m F_n(U, q)\), a filtered colimit of Eilenberg-Mac Lane spectra, is an Eilenberg-Mac Lane spectrum, with \(\pi_0\) the colimit of \(\pi_0 F_U^W\) for \(W \in \mathcal{S}_{\leq m}^n(U, q)\), so the colimit over such \(W\) of 
\[\bigoplus_{x \in W_{(d)}} \bb{Z} \cdot \eta_x\]
\(x \in W\) has \(\delta = d\) exactly when \(Z(x)\) meets faces properly, and as \(\delta\big(\pr_1(W)\big) \leq m\), \(\delta\big(\pr_1(x)\big) \leq m\) as well. Integral closed subschemes of \(U \times \Delta^q\) have unique generic point, proving the result.

The final result concerns two Eilenberg-Mac Lane spectra, so we just need to show that the statement holds after taking \(\pi_0\). The morphism is the inclusion 
\[\bigcup_{Z \in \mathcal{S}_{\leq m}(U)} z_{n + 1}(Z, q) \hookrightarrow z_{n + 1}(U, q),\]
which (2) shows is an isomorphism.
\end{proof}

\begin{proposition}[Graded pieces of spectral replacement]\label{prop:gr-sp-rep}
For any \(U \in X_\nis\), \(m \in \bb{Z}\), \(q \in \Delta\), \(\gr_m F_n(U, q)\) is an Eilenberg-Mac Lane spectrum, and 
\begin{equation}
\pi_0 \! \left(\gr_m F_n(U, q)\right) \cong \bigoplus_{x \in U_{(m)}} z_{n + 1}(U, q)_x,
\end{equation}
in which \(z_{n + 1}(U, q)_x\) be the free Abelian subgroup of \(z_{n + 1}(U, q)\) generated by \([C]\) such that \(\pr_1(\eta_C) = x\). This preserves simplicial structure, thus expressing \(\gr_m F_n(U, \bullet)\) as a direct sum of simplicial spectra over \(x \in U_{(m)}\).
\end{proposition}

\begin{proof}
We know that \(\fil_{m - 1} F_n(U, q) \to \fil_m F_n(U, q)\) is the cofibre of Eilenberg-Mac Lane spectra, so Eilenberg-Mac Lane itself. On \(\pi_0\), it is the inclusion of \([C]\) such that \(\delta\big(\pr_1(\eta_C)\big) \leq m - 1\). \(\gr_m F_n(U, q)\), the cofibre, is thus generated at the level of \(\pi_0\) by the remaining classes, which are exactly the \([C]\) such that \(\delta(\pr_1(C)\big) = m\). Partitioning these into the points of \(U\), they project to give the desired formula.

Given \([C]\) such that \(\pr_1(\eta_C) = x\), the Gysin map induced by a coface \(\iota \colon U \times \Delta^{q - 1} \hookrightarrow U \times \Delta^q\) takes \([C]\) to a linear combination of \([C_i]\), in which \(\eta_C \leadsto \iota(\eta_{C_i})\), thus 
\[\pr_1(\eta_{C_i}) = \pr_1\big(\iota(\eta_{C_i})\big) \subset \overline{\pr_1(\eta_C)} = Z(x)\]
This means that \(\pr_1(\eta_{C_i}) = x\), or has \(\delta < m\), and the latter terms vanish in \(\gr_m\). So this equivalence respects cofaces. A codegeneracy \(p \colon U \times \Delta^{q + 1} \to U \times \Delta^q\) is a trivial \(\bb{A}^1\)-bundle over \(U \times \Delta^q\) mapping \(C\) to another integral closed subscheme, and since \(p\) respects \(\pr_1\), this integral closed subscheme also projects to \(x\). So this equivalence respects codegeneracies, and thus holds at the simplicial level.
\end{proof}

\begin{remark}\label{rem:cyc-res}
For any \(x \in U_{(m)}\), there is a canonical map 
\[\rho_x \colon z_{n + 1}(U, q)_x \to \colim_{V \subset Z(x)} z_{n + 1}(V, q)\]
in which \(V\) decreases over all opens: given integral closed \(C \subset U \times \Delta^q\) such that \(\pr_1(\eta_C) = x\), then \(C \subset Z(x) \times \Delta^q\), meaning this lies in some \(V \times \Delta^q\), for open \(V \ni x\), and \(V' \subset V\) gives the same restriction to \(V' \times \Delta^q\). This is injective, as a cycle which becomes zero in the colimit must be zero on some \(V \times \Delta^q\), and thus not intersect \(V \times \Delta^q\), so not contain \(x\), which must be the zero cycle in \(z_{n + 1}(U, q)_x\).
\end{remark}

\begin{definition}[Cycle filtration]
Define a filtration on cycle groups by 
\[\fil_m \, z_{n + 1}(U, q) = \bb{Z} \! \left\{[C] \, \big| \, \delta\big(\pr_1(C)\big) \leq m\right\}\] 
This means that 
\[\gr_m \, z_{n + 1}(U, q) = \bigoplus_{x \in U_{(m)}} z_{n + 1}(U, q)_x\]
Just as in the proof of \ref{prop:gr-sp-rep}, this filtration is compatible with cofaces and codegeneracies on \(U \times \Delta^\bullet\), so this may be considered a filtration of simplicial Abelian groups.
\end{definition}

\begin{example}
Let \(U = \spec \bb{Z}[x]\), writing the generic point as \(\eta\). Let \(n = - 1\) and \(q = 2\), so that \(d = 1\). Define \(C \subset U \times \Delta^2\) as the closed subscheme given by \(x = t_1 = t_2\) and \(t_0 = 1 - 2x\). This means \(C \cong U\), with \(\delta(C) = 1 = d\). On the open \(V = D(x)\), \(C|_V\) meets faces properly, not intersecting \(V(t_1)\) or \(V(t_2)\) or the vertices \(V(t_i, t_j)\), and intersecting \(V(t_0)\) at 
\[\spec \bb{Z}[x]/(2x - 1) \cong \spec \bb{Z}[1/2],\]
which has relative dimension \(\delta = 0 = d - 1\). Hence the class of \(C|_V\) gives an element of the colimit
\[\colim_{V \subset Z(\eta)} z_0(V, 2)\]
On \(U\), \(C\) does not meet faces properly, as \(C\) intersects the vertex \(V(t_0, t_1)\) at \(V(x) \cong \spec \bb{Z}\), which has \(\delta = 0 > d - 2\). A cycle in \(z_0(U, 2)_\eta\) mapping to the class of \(C|_V\) in the colimit would need to restrict to the same element on some open \(V' \subset V\), but this would mean \(C\) is a component of that cycle, contradiction to the cycle meeting faces properly. Thus we have that this map is not guaranteed to be surjective. 
\end{example}

We want to be able to express all of our filtrations as direct sums over generic points of colimits, but see that this is not possible at the group level. The failure of surjectivity in the classical setting of intersection theory is remedied by moving lemmas, and we will use moving lemmas of Levine to solve this problem in our context as well after taking geometric realizations of the entire simplicial object. We first reduce to smooth neighborhoods in the colimit.

\begin{lemma}[Smooth reduction]\label{lemm:sm-red}
Let \(p \in \spec \bb{Z}\) be the image of \(x \in U_{(m)}\). Then the nonempty opens \(V \subset Z(x)\) of the form 

\begin{enumerate}[i)]

\item \(V\) smooth over \(\kappa(p) \cong \bb{F}_p\) of relative dimension \(m + 1\) for \(p\) prime

\item \(V\) smooth over \(\spec \bb{Z}\) of relative dimension \(m\) when \(p = 0\)

\end{enumerate}
are cofinal among the opens of \(Z(x)\).
\end{lemma}

\begin{proof}
In the case of \((i)\), \(Z(x)\) is an integral scheme over \(\bb{F}_p\), with \(\kappa(x)\) separably generated over \(\bb{F}_p\), thus \(Z(x)\) is smooth at its generic point. In the case of \((ii)\), \(Z(x) \to \spec \bb{Z}\) is flat and has smooth generic fibre at \(x\), thus has an open smooth locus containing \(x\) as we saw in the proof of \ref{lemm:HZBM}. In the case of \(p\) prime, then \(\mathrm{ht} \, p = 1\), so \(\dim V = \delta(V) + 1 = m + 1\).
\end{proof}

\begin{definition}[Shells]
For the presheaf of spectra \(\mathcal{G} = \bb{Z}(n)^\bm \lbrack -2n \rbrack\) or \(|F_n( -, \bullet)|\), define the \emph{\(m^\mathrm{th}\) shell of \(\mathcal{G}\),} denoted \(\mathcal{G}^{\lceil m \rceil}\), to be 
\[\mathcal{G}^{\lceil m \rceil}(U) = \colim_{Z' \in \mathcal{S}_{\leq m - 1}(U)} \fil_m \, \mathcal{G}(U - Z')\]
This has transition maps given by restriction: \(Z_1' \subset Z_2'\) means \(U - Z_2' \hookrightarrow U - Z_1'\) is open, so étale, thus induces a morphism \(\mathcal{G}(U - Z_1')\) to \(\mathcal{G}(U - Z_2')\) since \(\mathcal{G}\) is a presheaf on \(X_\nis\). If \(g \colon U' \to U\) is a morphism in \(X_\nis\), then étale pullback gives a functor \(\mathcal{S}_{\leq m - 1}^n(U, q) \to \mathcal{S}_{\leq m - 1}^n(U', q)\), thus a morphism 
\[\mathcal{G}^{\lceil m \rceil}(U) \to \mathcal{G}^{\lceil m \rceil}(U')\]
As a result, \(\mathcal{G}^{\lceil m \rceil}\) is a presheaf of spectra on \(X_\nis\).
\end{definition}

\begin{lemma}[Shells and vanishing]
For \(\mathcal{G} = |F_n( -, \bullet)|\) or \(\bb{Z}(n)^\bm\lbrack -2n \rbrack\), 
\[\colim_{Z' \in \mathcal{S}_{\leq m - 1}(U)} \fil_{m - 1} \, \mathcal{G}(U - Z') \simeq 0\]
As a result, restriction induces a morphism  
\[\begin{tikzcd}
\gr_m \, \mathcal{G}(U) \ar[r, "\mf{s}_U"] & \colim_{Z' \in \mathcal{S}_{\leq m - 1}(U)} \gr_m \, \mathcal{G}(U - Z')
\end{tikzcd}\]
This is natural in étale morphisms over \(X\) and in morphisms of filtered preseheaves of spectra, so that there is a commutative diagram 
\begin{equation}\label{eq:niv-shell}
\begin{tikzcd} [row sep = large, column sep = large]
\gr_m |F_n( -, \bullet)| \ar[d, "\lvert \mf{s} \rvert", swap] \ar[r, "\gr_m |\Phi|"] & \gr_m \, \bb{Z}(n)^\bm\lbrack -2n \rbrack \ar[d, "\mf{s}"]\\
\lvert F_n( -, \bullet)\rvert^{\lceil m \rceil} \ar[r, "|\Phi|^{\lceil m \rceil}"'] & \left(\bb{Z}(n)^\bm\right)^{\lceil m \rceil}\lbrack -2n \rbrack
\end{tikzcd}
\end{equation}
in the category \(\psh(X_\nis, \Sp)\).
\end{lemma}

\begin{proof}
Geometric realization commutes with colimits, meaning if we prove the result for \(q \in \Delta\), it holds for \(|F_n( -, \bullet)|\). 
\[\pi_0 \! \left(\colim_{Z' \in \mathcal{S}_{\leq m - 1}(U)} \fil_{m - 1} F_n(U - Z', q)\right) \simeq \colim_{Z' \in \mathcal{S}_{\leq m - 1}(U)} \pi_0 \! \left(\fil_{m - 1} F_n(U - Z', q)\right),\]
so if each cycle in \(\pi_0\) vanishes in the colimit, we have the first statement. Given class \([C]\), \(C \subset (U - Z') \times \Delta^q\), with \(\delta(Z_C) \leq m - 1\), so \(Z_C \in \mathcal{S}_{\leq m - 1}(U)\). Then taking \(Z'' = Z' \cup Z_C) \in \mathcal{S}_{\leq m - 1}(U)\) means that \([C]\) restricts to \(0\) under \(U - Z' \subset U - Z''\), so that the colimit is indeed the zero spectrum.

In the case of \(\mathcal{G} = \bb{Z}(n)^\bm\), up to adjoint equivalence we have 
\[\fil_{m - 1} \bb{Z}(n)^\bm(U - Z')\lbrack -2n \rbrack \simeq \colim_{Z \in \mathcal{S}_{\leq m - 1}(U - Z')} \bb{Z}(n)^\bm(Z) \lbrack - 2n \rbrack,\]
So
\[\colim_{Z' \in \mathcal{S}_{\leq m - 1}(U)} \fil_{m - 1} \, \bb{Z}(n)^\bm(U - Z')\lbrack -2n \rbrack = \colim_{Z' \in \mathcal{S}_{\leq m - 1}(U)} \! \left(\colim_{Z \in \mathcal{S}_{\leq m - 1}(U - Z')} \bb{Z}(n)^\bm(Z) \lbrack - 2n \rbrack\right)\]
Letting \(Z'' = Z' \cup \overline{Z}\), then the transition map induced by \(U - Z'' \subset U - Z'\) maps the \(\bb{Z}(n)^\bm(Z) \lbrack - 2n \rbrack\) term to \(0\). Both colimits are filtered, so this ends up as the zero spectrum as well. 

Given \(U \in X_\nis\), \(Z' \in \mathcal{S}_{\leq m - 1}(U)\), the diagra 
\[\begin{tikzcd} [row sep = huge, column sep = large]
\fil_{m - 1} \, |F_n(U - Z', \bullet)| \ar[d, "|\Phi|_{m - 1}", swap] \ar[r] & \fil_m \, |F_n(U - Z', \bullet)| \ar[d, "|\Phi|_m"] \ar[r] & \gr_m \, |F_n(U - Z', \bullet)| \ar[d, "\gr_m |\Phi|"]\\
\fil_{m - 1} \, \bb{Z}(n)^\bm(U - Z')\lbrack -2n \rbrack \ar[r] & \fil_m \, \bb{Z}(n)^\bm(U - Z')\lbrack -2n \rbrack \ar[r] & \gr_m \, \bb{Z}(n)^\bm(U - Z')\lbrack -2n \rbrack
\end{tikzcd}\]
commutes by construction.
Taking colimits is exact, so preserves the cofibre sequences in both rows and gives a commutative diagram 
\[\begin{tikzcd} [row sep = large]
 0 \ar[d] \ar[r] & \lvert F_n(-, \bullet)\rvert^{\lceil m \rceil}(U) \ar[d] \ar[r] & \displaystyle \colim_{Z' \in \mathcal{S}_{\leq m - 1}(U)} \gr_m |F_n(U - Z', \bullet)| \ar[d]\\
0 \ar[r] & \left(\bb{Z}(n)^\bm\right)^{\lceil m \rceil}(U)\lbrack -2n \rbrack \ar[r] & \displaystyle \colim_{Z' \in \mathcal{S}_{\leq m - 1}(U)} \gr_m \, \bb{Z}(n)^\bm(U - Z')\lbrack -2n \rbrack
\end{tikzcd}\]
so that the final two terms in each row are equivalent. The restriction induced by \(U \subset U - Z'\) gives a morphism \(\gr_m |F_n(U, \bullet)| \to \gr_m |F_n(U - Z', \bullet)|\) compatible with restrictions \(U - Z_1' \subset U - Z_2'\), so that \(\gr_m |F_n(U, \bullet)|\) maps to the colimit, thus to \(|F_n( -, \bullet)|^{\lceil m \rceil}(U - Z')\) under equivalence. Likewise, there is a morphism from \(\bb{Z}(n)^\bm(U) \lbrack -2n \rbrack\) to the colimit, thus to \(\big(\bb{Z}(n)^\bm\big)^{\lceil m \rceil}(U)\lbrack -2n \rbrack\). \(\gr_m|\Phi|(U)\) sends \(\gr_m|F_n(U, \bullet)|\) to \(\gr_m \, \bb{Z}(n)^\bm(U) \lbrack -2n \rbrack\), compatibly through each restriction, so factors through the colimit, giving the commutative diagram \ref{eq:niv-shell} at \(U\). The fact that this holds at the level of presheaves is a result of \(\mathcal{G}^{\lceil m \rceil}\) being a presheaf, as well as naturality of diagrams of cofibre sequences.
\end{proof}

\begin{definition}[Decomposable pairs]
A pair \((Z, Z')\) such that \(Z' \in \mathcal{S}_{\leq m - 1}(U)\), and \(Z \subset U - Z'\), in which \(Z\) has the form 
\[Z = \coprod_{i = 1}^s V_i,\]
where \(V_i\) equals \(Z(x_i) - Z'\) for distinct \(x_i \in U_{(m)}\) and is of the form \ref{lemm:sm-red} is called \emph{decomposable}. 
\end{definition}

\begin{lemma}
Partially ordering pairs \((Z, Z')\) such that \(Z' \in \mathcal{S}_{\leq m - 1}(U)\), \(Z \in \mathcal{S}_{\leq m}(U - Z')\) by 
\[(Z_1, Z_1') \leq (Z_2, Z_2') \text{ if } Z_1' \subset Z_2', Z_1 - Z_2' \subset Z_2,\]
then decomposable pairs are cofinal among such pairs.
\end{lemma}

\begin{proof}
Take a pair \((Z, Z')\). Let \(x_1, \dots , x_s\) be the generic points of the irreducible components of \(\delta = m\) of \(Z\), with closure taken in \(U\). Define \(Z''\) to be the union of \(Z'\) with the irreducible components of \(Z\) of \(\delta < m\), the pairwise intersections \(Z(x_i) \cap Z(x_j)\), for \(i \neq j\), and the complement of the smooth locus of each \(Z(x_i)\).

\(Z''\) is a finite union of closed subschemes of \(\delta < m\), thus is of this form itself. Removing \(Z''\) from \(Z(x_i)\) excises all intersections of top-dimensional components and the non-smooth locus. In addition, it guarantees that the \(V_i\) do not meet. So \(Z\) decomposes into a disjoint union as desired. As \(Z''\) contains \(Z'\), we have \((Z, Z') \leq (Z, Z'')\), giving cofinality.
\end{proof}

\begin{corollary}\label{corr:shell-decomp}
Pushing forward from \(V_i\) induces an equivalence 
\[\mathcal{G}^{\lceil m \rceil}(U) \simeq \bigoplus_{x \in U_{(m)}} \colim_{V \subset Z(x)} \mathcal{G}(V),\]
in which \(\mathcal{G} = |F_n( -, \bullet)|\) or \(\bb{Z}(n)^\bm \lbrack -2n \rbrack\), giving a commutative diagram
\begin{equation}\label{eq:shell-decomp}
\begin{tikzcd} [row sep = large, column sep = large]
\displaystyle \bigoplus_{x \in U_{(m)}} \colim_{V \subset Z(x)} |F_n(V, \bullet)| \ar[d] \ar[r] & \displaystyle \bigoplus_{x \in U_{(m)}} \colim_{V \subset Z(x)} \bb{Z}(n)^\bm(V) \lbrack -2n \rbrack \ar[d]\\
\lvert F_n( -, \bullet)\rvert^{\lceil m \rceil}(U) \ar[r, "|\Phi|^{\lceil m \rceil}(U)", swap] & \left(\bb{Z}(n)^\bm\right)^{\lceil m \rceil}(U)\lbrack -2n \rbrack
\end{tikzcd}
\end{equation}
\end{corollary}

\begin{proof}
By definition of the filtration on \(\mathcal{G} = \bb{Z}(n)^\bm \lbrack -2n \rbrack\), and by \ref{lemm:Z-S-fil} (3) on \(\mathcal{G} = |F_n( -, \bullet)|\), we have 
\[\mathcal{G}^{\lceil m \rceil}(U) = \colim_{Z' \in \mathcal{S}_{\leq m - 1}(U)} \colim_{Z \in \mathcal{S}_{\leq m}(U - Z')} \mathcal{G}(Z)\]
By cofinality, this is equivalent to  
\[\colim_{Z' \in \mathcal{S}_{\leq m}(U)} \colim_{V_i \in \mathcal{S}_{\leq m}(U - Z')} \bigoplus_{i = 1}^s \mathcal{G}(V_i),\]
Given \(Z_1' \subset Z_2'\) in \(\mathcal{S}_{\leq m - 1}(U)\), the transition map restricts the affine smooth open \(Z(x_i) - Z_1''\) to affine smooth open \(Z(x_i) - Z_2''\). There are no restrictions to \(Z(x_j) - Z_2''\) for \(i \neq j\), as all intersections of \(Z(x_i)\) and \(Z(x_j)\) lie in \(Z_2''\) as well as \(Z_1''\). Thus restrictions of the \(V_i\) respect the direct sum. Increasing \(Z\) induces inclusion into a direct sum of more points of \(U_{(m)}\), so we are done.
\end{proof}

We'll use this diagram to show that \(\gr_m |\Phi|\) is a local equivalence. To do so, we first need to show that the top two vertical legs, the \(\mf{s}_U\), are local equivalences.

\begin{proposition}[Borel-Moore shell equivalence]\label{prop-BM-shell}
For any \(U \in X_\nis\), the canonical morphism \(\mf{s}_U \colon \gr_m \, \bb{Z}(n)^\bm(U) \lbrack -2n \rbrack \to \big(\bb{Z}(n)^\bm)^{\lceil m \rceil}(U)\lbrack -2n \rbrack\) is an equivalence.
\end{proposition}

\begin{proof}
Take \(Z \in \mathcal{S}_{\leq m}(U)\), and \(Z' \in \mathcal{S}_{\leq m - 1}(U)\). Write \(i_Z \colon Z \hookrightarrow U\) and \(i \colon Z' \cap Z\hookrightarrow Z\) for the closed embeddings, and \(j \colon Z' - Z \hookrightarrow Z\) for the complementary open embedding of the second morphism. Applying the exact functor \(\map_{\sh(Z)}(\mathbbm{1}_Z, -)\) to the cofibre sequence 
\[\begin{tikzcd}
i_* i^! i_Z^! f_U^! s^{- n} \Kgl_\bb{Z} \ar[r] & i_Z^! f_U^! s^{- n} \Kgl_\bb{Z} \ar[r] & j_*j^! i_Z^! f_U^! s^{- n} \Kgl_\bb{Z}
\end{tikzcd}\]
in \(\sh(Z)\) gives, up to adjoint equivalence, the cofibre sequence 
\[\begin{tikzcd} [column sep = small]
\map_{\sh(U)}\big(\mathbbm{1}_U, \Gamma_{Z \cap Z'} f_U^! s^{- n} \Kgl_\bb{Z}\big) \ar[r] & \map_{\sh(U)}\big(\mathbbm{1}_U, \Gamma_Z f_U^! s^{- n} \Kgl_\bb{Z}\big) \ar[r] & \map_{\sh(U)}\big(\mathbbm{1}_U, \Gamma_{Z - Z'} f_U^! s^{- n} \Kgl_\bb{Z}\big)
\end{tikzcd}\]
in \(\Sp\), which is natural in the pair. Taking colimits over \(\mathcal{S}_{\leq m}(U) \times \mathcal{S}_{\leq m - 1}(U)\) preserves cofibres. The functor \(\mathcal{S}_{\leq m}(U) \times \mathcal{S}_{\leq m - 1}(U) \to \mathcal{S}_{\leq m - 1}(U)\), \((Z, Z') \mapsto (Z \cap Z')\) is cofinal, so we may take the first colimit over this latter poset. Likewise, \(\mathcal{S}_{\leq m}(U) \times \mathcal{S}_{\leq m - 1}(U) \to \mathcal{S}_{\leq m}(U)\), \((Z, Z') \mapsto Z\) and \(\mathcal{S}_{\leq m - 1} \times \mathcal{S}_{\leq m}(U) \mapsto \mathcal{S}_{\leq m}(U - Z')\), \(Z \mapsto Z - Z'\) are cofinal, so we may take the colimit over these posets in the second and third terms of the cofibre sequence, respectively. Taking those colimits gives the cofibre sequence 
\[\begin{tikzcd} 
\fil_{m - 1} \, \bb{Z}(n)^\bm(U) \lbrack -2n \rbrack \ar[r] & \fil_m \, \bb{Z}(n)^\bm(U) \lbrack -2n \rbrack \ar[r] & \left(\bb{Z}(n)^\bm\right)^{\lceil m \rceil}(U)\lbrack -2n \rbrack
\end{tikzcd}\]
The second map is given componentwise by restriction, so the induced morphism, which is precisely \(\mf{s}_U \colon \gr_m \, \bb{Z}(n)^\bm(U) \lbrack -2n \rbrack \to \big(\bb{Z}(n)^\bm\big)^{\lceil m \rceil}(U)\lbrack -2n \rbrack\), is an equivalence.
\end{proof}

To prove an analogous statement for \(\mathcal{G} = |F_n( -, \bullet)|\), we'll import some localization results from the Bloch-Levine complex.

\begin{definition}
For a finite type separated scheme \(U\), and \(m \geq -1\), define 
\[z_{n + 1}(U, \bullet)^{\lceil m \rceil} = \colim_{Z' \in \mathcal{S}_{\leq m - 1}(U)} \fil_m z_{n + 1}(U - Z', \bullet)\]
By the proof of \ref{eq:niv-shell} for \(|F_n( -, \bullet)|\), there is a natural morphism \(\rho_U^{\lceil m \rceil} \colon \gr_m  z_{n + 1}(U, \bullet) \to z_{n + 1}(U, \bullet)^{\lceil m \rceil}\).

Using the same cofinality argument proof of \ref{prop-BM-shell}, the shell of a cycle may equivalently be defined as 
\[z_{n + 1}(U, \bullet)^{\lceil m \rceil} = \bigoplus_{x \in U_{(m)}} \colim_{V \subset Z(x)} z_{n + 1}(V, \bullet)\]
\end{definition}

\begin{remark}\label{rem:shell-cyc-spc}
The proof of \ref{eq:niv-shell} (3) shows that the shells and filtrations for the cycle groups are defined exactly so that \(\pi_0\) of shells of the unrealized \(F_n(U, \bullet)\) are shells of \(z_{n + 1}(U, \bullet)\), and likewise for the filtrations. Consequently, \(\pi_0\) maps the unrealized \(\mf{s}_U\) to \(\rho_U\).
\end{remark}

\begin{theorem}[Cycle shell equivalence]\label{thm:cyc-shell-gr}
Let \(Y\) be a finite type separated over a field or semi-local Dedekind domain. Then \(\rho_Y^{\lceil m \rceil}\) is an equivalence.
\end{theorem}

\begin{proof}
Let \(Y\) be a scheme with the properties above. For any pair \(Z \in \mathcal{S}_{\leq m}(Y)\), \(\mathcal{S}_{\leq m - 1}(Y)\), we have an open-closed pair \(Z \cap Z' \hookrightarrow Z\) and \(Z - Z \cap Z' \hookrightarrow Z\), all schemes finite type and separated over a field or semi-local Dedekind domain. Levine \cite[Theorem 1.7]{LevTechLoc} proves that 
\[\begin{tikzcd} 
z_{n + 1}(Z \cap Z', \bullet) \ar[r] & z_{n + 1}(Z, \bullet) \ar[r] & z_{n + 1}(Z - Z \cap Z', \bullet)
\end{tikzcd}\]
is a distinguished triangle of simplicial Abelian groups. At this point, simply run the proof of \ref{prop-BM-shell} again.
\end{proof}

\begin{lemma}[Continuity of cycles]
For finite type separated \(U\) and prime \(p\), write 
\[U_{(p)} = \colim_{p \nmid N} U[1/N], \quad U_\bb{Q} = \colim_N U[1/N]\]
Then restriction induces isomorphisms of simplicial Abelian groups, 
\begin{equation}\label{eq:cyc-cont}
\colim_{p \nmid N}  z_{n + 1} \big(U[1/N], \bullet \big) \simeq z_{n + 1}\big(U_{(p)}, \bullet\big), \quad \colim_N z_{n + 1} \big(U[1/N], \bullet \big) \simeq z_{n + 1}\big(U_\bb{Q}, \bullet\big)
\end{equation}
The equivalences are compatible with \(\fil_\star\) and \(( - )^{\lceil m \rceil}\).
\end{lemma}

\begin{proof}
\(U_{(p)} \cong U \times \bb{Z}_{(p)}\), so a point in \(U_{(p)}\) is a point in \(U\) which lies over \((p)\) or \(0\), and has the same residue field and \(\delta\) as that point. As a result, an integral closed subscheme \(C \subset U_{(p)} \times \Delta^q\) is the restriction of its closure \(\overline{C}\) in \(U \times \Delta^q\).

If \(C\) meets all faces properly, we must show that a sufficiently small restriction of \(\overline{C}\) does as well. If \(C_0\) is an irreducible component of \(\overline{C}\) that doesn't meet all faces \(U \times \Delta^{q - 1}\) properly, then its generic point must lie over a prime \(p_0 \neq p\). This means that \(C_0\) lies over \(p_0\), and only finitely such irreducible components can exist. Hence for large enough \(N\), \(\overline{C}|_{U[1/N]}\) meets all faces \(U \times \Delta^{q - 1}\) properly. Restricting to a further open \(U[1/N']\) preserves \(\delta\) and the condition of meeting faces properly, thus the isomorphism holds at \(q\).

To see that this holds at the simplicial level, note that as \(U \times \Delta^\bullet\) is a cosimplicial scheme, the projection morphism commutes with taking cofaces or codegeneracies, so that if a generic point in \(U \times \Delta^q\) does not lie over \((p)\) or \(0\), then the generic point obtained by smoothly pulling back to \(U \times \Delta^{q + 1}\) will not, nor will any of the generic points obtained by pulling back along the closed embedding to \(U \times \Delta^{q - 1}\). This compatibility paired with \(z_{n +1}(U_{(p)}, \bullet)\) being a simplicial Abelian group ensures that the restriction isomorphism respects the simplicial structure.

Compatibility with the filtration follows immediately from commutativity of direct sums with colimits. To prove shell compatibility, we use the direct sum definition. Mapping \(C \subset U_{(p)} \times \Delta^q\) to its closure in \(U \times \Delta^q\) is compatible with \(\pr_1\), and preserves \(\delta\), so is compatible with the filtration on the cycle groups. And if \(x \in U_{(m)}\) lies over a prime other than \(p\), then \(z_{n + 1}(U[1/N], q) = 0\) for sufficiently large \(N\), meaning \(z_{n + 1}(U_{(p)}, q)\), which is a colimit of these, is also zero. Thus 
\[\colim_{p \nmid N} \! \left(z_{n + 1}(U[1/N], q)^{\lceil m \rceil}\right) = \colim_{p \nmid N} \bigoplus_{x \in U[1/N]_{(m)}} \colim_{V \subset Z(x)} z_{n + 1}(V, q)\]
removes all points of \(U\) not lying over \((p)\) or \(0\), reducing the direct sum to being over \(U_{(p)}\). Opens in the reduced closure of \(x\) in \(U\), for \(x\) lying over \((p)\) or \(0\), restrict cofinally to opens over the reduced closure of \(x\) in \(U_{(p)}\), which we write as \(Z(x)_{(p)}\). So this evaluates to 
\[\bigoplus_{x \in {U_{(p)}}_{(m)}} \colim_{V \subset Z(x)_{(p)}} z_{n + 1}(V, q) = \big(z_{n + 1}(U_{(p)}, q)\big)^{\lceil m \rceil}\]
giving compatibility of taking shells with localizing at \(p\).

The case of \(U_\bb{Q}\) is the exact same, except we invert every nonzero prime, so only deal with primes of \(U\) lying over \(0\).
\end{proof}

Finally, we need one technical lemma regarding simplicial spectra.

\begin{lemma}\label{lemm:dsc-sim-sp} Let \(S_\bullet, T_\bullet\) be simplicial spectra which are Eilenberg-Mac Lane at each \(q \in \Delta\). If \(\phi \colon S_\bullet \to T_\bullet\) is a morphism such that \(\pi_0\phi\) is a weak equivalence of simplicial Abelian groups, then \(|\phi| \colon |S_\bullet| \to |T_\bullet|\) is an equivalence of spectra.
\end{lemma}

\begin{proof}
We first prove that there is an isomorphism 
\[\pi_i|S_\bullet| \cong H_i\big(N(\pi_0 S_\bullet)\big),\]
natural in \(S_\bullet\). Equip \(|S_\bullet|\) with the skeletal filtration, which yields a first quadrant spectral sequence \(E_{s, t}^2 = H_s\big(N(\pi_t S_\bullet)\big) \Longrightarrow \pi_{s + t} |S_\bullet|\) \cite[Remark 1.2.4.4]{HA}.  In this case, \(\pi_t S_\bullet = 0\) for \(t \neq 0\), so \(E^2\) is concentrated in the row \(t = 0\), and the spectral sequence collapses. Applying this and Whitehead's theorem thus proves that \(|\phi|\) is an equivalence.
\end{proof}

\begin{proposition}
The morphism \(|\mf{s}| \colon \gr_m |F_n( -, \bullet)| \to |F_n( -, \bullet)|^{\lceil m \rceil}\) is a Zariski-local equivalence of presheaves \(\psh(X_\nis, \Sp)\).
\end{proposition}

\begin{proof}
Fix \(U \in X_\nis\), and take some point \(x \in U\). Its image in \(\spec \bb{Z}\) is the ideal \((p)\), for \(p\) prime or 0. Write \(B = \bb{Z}_{(p)}\) or \(\bb{Q}\) for the base, respectively. Opens \(U'[1/N]\) for \(U'\) an affine neighborhood \(U'\) of \(x\) such that \(N \not \in (p)\) are cofinal among the open neighborhoods of \(x\).

Associated gradeds and colimits commute with geometric realization, so \(|\mf{s}|_U\) is the geometric realization of a morphism between simplicial spectra which are Eilenberg-Mac Lane at each \(q \in \Delta\). This holds for the stalk of \(|\mf{s}|_U\) at \(x\), \(|\mf{s}|_x\), as well, since it is computed as a colimit. If \(\pi_0 \mf{s}_x\) is a weak equivalence of simplicial Abelian groups for each \(q\), then by \ref{lemm:dsc-sim-sp}, \(|\mf{s}|_x\) is an equivalence. Compute: 
\[\pi_0 \mf{s}_x \colon \colim_{N \not \in (p)} \colim_{x \in U'} \pi_0 \! \left(\gr_m F_n(U'[1/N], q)\right) \to \colim_{N \not \in (p)} \colim_{x \in U'} \pi_0 \! \left(F_n(U'[1/N], q)^{\lceil m \rceil}\right)\]
Interchanging the colimits and commuting the filtered colimit over \(N\) with \(\pi_0\) and \(\gr_m\), this is the right leg of 
\[\begin{tikzcd} [row sep = large, column sep = large]
\displaystyle \colim_{x \in U'}  \pi_0 \! \left(\gr_m \displaystyle \colim_{N \not \in (p)} F_n(U'[1/N], q)\right) \ar[d] \ar[r, "\stackrel{\ref{eq:cyc-cont}} \sim"] & \displaystyle \colim_{x \in U'} \pi_0 \! \left(\gr_m F_n(U_B', q)\right) \ar[d]\\
\displaystyle \colim_{x \in U'}  \pi_0 \! \left(\displaystyle \colim_{N \not \in (p)} F_n(U'[1/N], q)^{\lceil m \rceil} \right) \ar[r, "\stackrel{\ref{eq:cyc-cont}} \sim"] & \displaystyle \colim_{x \in U'} \pi_0 \! \left(F_n(U_B', q)^{\lceil m \rceil} \right)
\end{tikzcd}\] 
By \ref{rem:shell-cyc-spc}, that map is exactly the filtered colimit of 
\[\rho_{U_B'} \colon \gr_m z_{n + 1}(U_B', q) \to z_{n + 1}(U_B', q)^{\lceil m \rceil},\]
thus an equivalence by \ref{thm:cyc-shell-gr}. Thus \(\pi_0 \mf{s}_x\) is a filtered colimit of weak equivalences, so a weak equivalence itself, meaning \(|\mf{s}|_x\) is an equivalence, so that \(|\mf{s}|\) is a Zariski-local, hence Nisnevich-local, equivalence.
\end{proof}

\subsection{Graded equivalence: the homotopy coniveau tower}
At this point, we have pared down the argument to showing that the lower morphism in \ref{eq:shell-decomp} is an equivalence. To finish this off, we utilize some results of Levine on homotopy coniveau towers \cite{LevHC}, whose great influence on this paper should be apparent.

\begin{construction}[The homotopy coniveau tower]
This construction and the following remark are laid out in \cite[\S 2]{LevHC}. Fix a base scheme \(B\), and let \(E \in \psh(\sm_B, \Sp)\). For smooth \(X \to \spec B\), and closed subscheme \(W \subset X\), Levine writes 
\[E^W(X) = \mathrm{fib}\big(E(X) \to E(X - W)\big),\]
and \(\mathcal{S}_X^{(p)}(r)\) for the filtered poset of closed subschemes \(W \subset X \times \Delta^r\) such that 
\[\codim_{X \times F} \big(W \cap (X \times F)\big) \geq p\]
for each face. Define 
\[E^{(p)}(X, r) = \colim_{W \in \mathcal{S}_X^{(p)}(r)} E^W(X \times \Delta^r)\]
\(E^{(p)}(X, \bullet)\) is a simplicial spectrum, and there is a natural tower of simplicial spectra, 
\[\begin{tikzcd} 
\cdots \ar[r] & E^{(p + 1)}(X, \bullet) \ar[r] & E^{(p)}(X, \bullet) \ar[r] & \cdots \ar[r] & E^{(0)}(X, \bullet)
\end{tikzcd}\]
This is called the \emph{homotopy coniveau tower of \(E\)}. Varying \(X \in \sm_B\) and taking geometric realizations, we may view \(|E^{(p)}|\) as a presheaf of spectra on \(\sm_B\).
\end{construction}

\begin{remark}[Forgetting support]\label{rem:supp-gone}
As \(\mathcal{S}_X^{(0)}(r)\) has terminal object \(X \times \Delta^r\) and \(E(0) = *\), \(E^{(0)}(X, r) \simeq E^W(X \times \Delta^q) \simeq E(X \times \Delta^q)\). Thus when \(E\) is \(\bb{A}^1\)-invariant, the morphism \(E^{(p)} \to E^{(0)}\) may be viewed as \(E^{(p)} \to E\). This morphism is given by taking the colimit over \(W \in \mathcal{S}_X^{(p)}(r)\) of the composition 
\[\begin{tikzcd}
E^W(X \times \Delta^q) \ar[r] & E(X \times \Delta^q) \ar[r, "\sim"] & E(X),
\end{tikzcd}\]
then taking geometric realizations gives a morphism \(|E^{(p)}| \to E\) in \(\psh(\sm_B, \Sp)\). 
\end{remark}

\begin{remark}[\(\bb{P}^1\)-delooping]
There is \cite[\S 3.2]{BEM} a functor \(\omega^\infty \colon \sh(X) \to \shv_{\nis, \bb{A}^1}(\sm_X, \Sp)\), called \emph{\(\bb{P}^1\)-delooping}, taking a motivic spectrum \(E \in \sh(X)\) to an \(\bb{A}^1\)-invariant Nisnevich sheaf of spectra by the rule 
\[\omega^\infty E(Y) = \map_{\sh(X)}\big(M_X(Y), E\big)\]
Sheaves of this form are one model for the category of \(S^1\)-motivic spectra, \(\sh_{S^1}(X)\). This category only appears a single time, used in a theorem of Levine, so we do not expound upon it any further.
\end{remark}

The utility of the functor \(\omega^\infty\) is that it allows us to compare homotopy coniveau with the spectral replacement. The ``level zero'' manifestation of this is a comparison between motivic Borel-Moore homology and the \(\bb{P}^1\)-delooping:

\begin{lemma}[Coniveau and spectral replacement]\label{lemm:LvT}
Let \(V\) be as in \ref{lemm:sm-red}. Write \(c = m - n\) for the relative codimension, and \(B = \bb{F}_p\) or \(\bb{Z}[1/N]\) for the base over which \(V\) is smooth. For each \(q \in \Delta\), the following hold:
\begin{enumerate}

\item \(f_{V, q}^! s^{- d} \Kgl_\bb{Z} \simeq s^c \Kgl_{V \times \Delta^q}\), meaning 
\[\bb{Z}(n)^\bm(V) \lbrack -2n \rbrack \simeq \map_{\sh(V)}\big(\mathbbm{1}_V, s^c \Kgl_V\big) \simeq \bb{Z}(c)^\bb{A}(V)[2c]\]

\item \(\mathcal{S}_n(V, q) = \mathcal{S}_V^{(c)}(q)\).

\item Let \(E = \omega^\infty s^c \Kgl_B\). Then there are equivalences of spectra 
\[\bb{Z}(n)^\bm(V) \lbrack -2n \rbrack \simeq E(V)\]

\end{enumerate}
\end{lemma}

\begin{proof}
\(f_{V, q}\) is smooth, and factors as \(f_V \circ \pr_1\). The fact that \(\pr_1\) is smooth of relative dimension \(q\) paired with smooth and possibly arithmetic purity gives \(f_{V, q}^! \simeq f_V^! \otimes \bb{T}^{\otimes m}\), so commuting , so that 
\[f_{V, q}^! \simeq \pr_1^! \circ f_V^! \simeq \bb{T}^{\otimes (m + q)} \otimes \pr_1^* \circ f_V^* \simeq \bb{T}^{\otimes (m + q)} \otimes f_{V, q}^*,\]
so Bott periodicity and commutativity of slices with smooth pullback gives the first part of (1). The second part follows by adjointness under \((\pr_1)_\sharp \dashv \pr_1^*\).

\(W \in \mathcal{S}_n(V, q)\) is a closed subscheme of \(V \times \Delta^q\) with \(\delta(W) = d\), and \(V \times F\) is itself integral, with \(\delta(V \times F) = m + q - \codim_{V \times \Delta^q} \, F\) for any face \(F\) of \(V \times \Delta^q\). Hence the condition that \(\delta\big(W \cap (V \times F)\big) \leq d - \codim_{V \times \Delta^q} \, F\) is equivalent to the condition that \(\codim_{V \times F}\big(W \cap (V \times \Delta^q)\big) \geq m + q - d = c\), giving (2).

\(E(V)\) is the mapping spectrum 
\[\map_{\sh(B)}\big(M_B(V), s^c \Kgl_B\big) \simeq \map_{\sh(B)}\big((f_V)_\sharp \mathbbm{1}_V, s^c \Kgl_B\big) \simeq \map_{\sh(V)}\big(\mathbbm{1}_V, f_V^* s^c \Kgl_B\big)
\]
As \(f_V^*\) is smooth, it commutes with formation of slices \cite[Remark 3.41]{BEM}, pulling back \(\Kgl_B\) to \(\Kgl_V\). Paired with (1) at \(q = 0\), this is 
\[\map_{\sh(V)}\big(\mathbbm{1}_V, f_V^! s^{- n} \Kgl_\bb{Z}\big) = \bb{Z}(n)^\bm(V) \lbrack -2n \rbrack,\]
proving (3).
\end{proof}

\begin{remark}[Setup]
Fix \(U \in X_\nis\), and \(i \colon V \hookrightarrow U\) closed such that \(V\) is of the form \ref{lemm:sm-red} and \(\delta(V) = m\). Write \(Y = U \times \Delta^q\), and \(Z_V = \overline{g_{U}(V)}\), for \(g_{U} \colon U \to X\) the structure morphism. As this is étale, \(Z_V \in \mathcal{S}_{\leq m}(X)\).
\end{remark}

\begin{construction}
Let \(\xi \colon \mathcal{T}_V \to \Delta^\op\) be the Grothendieck construction of the functor \(q \mapsto \mathcal{S}_n(V, q)\), the poset defined as in \ref{rem:cyc-in-closed} for (étale)-locally closed \(V\). We obtain a commutative diagram of functors 
\begin{equation}\label{eq:intw-fib}
\begin{tikzcd} [row sep = large, column sep = large]
\mathcal{T}_V \ar[d, "\xi", swap] \ar[r, "j_V"] & \mathcal{T} \ar[d, "\pi"]\\
\Delta^\op \ar[r, "\jmath", swap] & \mathcal{B}
\end{tikzcd}
\end{equation}
in which \(j_V(q, W) = (U, q, m, W, Z_V)\) and \(\jmath(q) = (U, q, m)\). Define 
\[F_n(V, \bullet) = \lan_\xi (F \circ j_V),\]
which thus has value \(F_n(V, q)\) (in the sense of \ref{con:spec-rep}) at \(q\).
\end{construction}

\begin{definition}[The morphism \(\Xi\)]
Let \(\mf{G}_V(\bullet) \colon \Delta^\op \to \Sp\) be the functor 
\[\mf{G}_V(q) = \map_{\sh(X)}\big(M_X(Y/V \times \Delta^q), \mf{G}_X\big),\]
which is equivalent to 
\[q \mapsto \bb{Z}(n)^\bm(V \times \Delta^q)\lbrack -2n \rbrack\]
The morphism of pairs \((Y, V \times \Delta^q) \to (Y, W)\) induces a natural transformation \(\tilde{\alpha} \colon F \circ j_V \Longrightarrow \xi^* \mf{G}_V\), which corresponds to an adjoint natural transformation 
\[\alpha \colon \lan_\xi (F \circ j_V) \Longrightarrow \mf{G}_V(\bullet)\]
Likewise, the morphism of pairs \((Y, V \times \Delta^q) \to (U, V)\) induced by \(\pr_1\), gives a morphism 
\[\psi^*\map_{\sh(X)}\big(M_X(U/V), \mf{G}_X) \to \mf{G}_V,\]
an equivalence at each \(q\) by \(\bb{A}^1\)-invariance, so denote its inverse as \(\beta\). Define \(\Xi_V\) to be the composite 
\[\begin{tikzcd}
\lvert F_n(V, \bullet) \rvert \ar[r, "\lvert \alpha \rvert"] & \lvert \mf{G}_V(\bullet) \rvert \ar[r, "\lvert \beta \rvert"] & \map_{\sh(X)}\big(M_X(U/V), \mf{G}_X\big)
\end{tikzcd}\]
The final object is equivalent to \(\bb{Z}(n)^\bm(V)\lbrack -2n \rbrack\) by the localization sequence.
\end{definition}

While very similar to our previous construction, this step is necessary to compare Levine's functoriality, which is only defined on smooth schemes over perfect fields, with the functoriality of the spectral replacement, which is defined in a broader setting.

\begin{definition}\label{def:i-push}
Define \(i_* \colon \bb{Z}(n)^\bm(V) \lbrack -2n \rbrack \to \fil_m \bb{Z}(n)^\bm(U) \lbrack -2n \rbrack\) to be as the composite 
\[\begin{tikzcd} [row sep = large]
\map_{\sh(X)}\big(M_X(U/V), \mf{G}_X\big) \ar[r, "h"] & \map_{\sh(X)}\big(M_X(U/g^{- 1} Z_V), \mf{G}_X) \ar[dl, "(ii)"]\\
\map_{\sh(X)} \Big(M_X(U), \underline{\Hom}\big(M_X(X/Z_V), \mf{G}_X\big)\Big) \ar[r, "(iii)", swap] & \map_{\sh(X)}\big(M_X(U), \mf{G}_X^{\leq m}\big)
\end{tikzcd}\]
for \((ii)\) and \((iii)\) the same as \ref{con:phi}. This is a coherent upgrade of inclusion of a colimitand of the filtration on (shifted) motivic Borel-Moore homology into the colimit.
\end{definition}

\begin{proposition}[Intertwining functorialities]\label{prop:int-fun}
If we write \(|j_V|_* \colon |F_n(V, \bullet)| \to \fil_m |F_n(U, \bullet)|\) for the geometric realization of the Beck-Chevalley transformation \(\lan_\xi j_V^* \Longrightarrow \jmath^* \lan_\pi\) induced by the Cartesian square \eqref{eq:intw-fib}, then 
\[|\Phi|_m(U) \circ |j_V|_* \simeq i_* \circ \Xi_V\]
\end{proposition}


\begin{proof}
\(g_U(V) \subset Z_V\) means \(W \subset V \times \Delta^q \subset (g_U \circ \pr_1)^{- 1}(Z_V)\), giving a factorization of the morphism of pairs, 
\[\begin{tikzcd}
\big(Y, (g_U \circ \pr_1)^{- 1}(Z_V)\big) \ar[r, "\id"] & (Y, V \times \Delta^q) \ar[r, "\id"] & (Y, W)
\end{tikzcd}\]
Contravariant functoriality of \(\map_{\sh(X)}\big(M_X( -/-), \mf{G}_X\big)\) 
on this factors the map \((i)\) in \ref{con:phi} as \(\tilde{\alpha}\) composed with \(\map_{\sh(X)}\big(M_X( -/-), \mf{G}_X\big)\) 
applied to the upper horizontal morphism of the diagram 
\[\begin{tikzcd} [row sep = large, column sep = large]
\big(Y, (g_U \circ \pr_1)^{- 1} Z_V\big) \ar[d, "\pr_1", swap] \ar[r, "\id"] & \big(Y, V \times \Delta^q\big) \ar[d, "\pr_1"]\\
\big(U, g^{- 1}Z_V\big) \ar[r, "\id", swap] & \big(U, V\big),
\end{tikzcd}\]
Strict commutativity of that diagram means this upper horizontal morphism is \(\pr_1^* \circ h \circ \beta\), that is, \((i) \simeq \pr_1^* \circ h \circ \beta \circ \tilde{\alpha}\). Naturality of the second and third natural transformations \((ii)\) and \((iii)\) in \ref{con:phi}, and the fact that the fourth morphism \((iv)\) is exactly the inverse of \(\pr_1^*\) means that 
\[\phi \simeq (\pr_1^*)^{- 1} \circ (iii) \circ (ii) \circ \pr_1^* \circ h \circ \beta \circ \tilde{\alpha} \simeq (iii) \circ (ii) \circ h \circ \beta \circ \tilde{\alpha} \simeq i_* \circ \beta \circ \tilde{\alpha},\]
the final equivalence by \ref{def:i-push}. This gives an equivalence of functors 
\[j_V^*\phi \simeq C_{i_*} \circ \xi^*\beta \circ \tilde{\alpha}\]
on the value \((q, W)\), and every step is natural, so we have an equivalence of are natural transformations \(j_V^*F \Longrightarrow \xi^* C_\mathcal{M} \simeq C_\mathcal{M}\), for 
\[\mathcal{M} = \fil_m \bb{Z}(n)^\bm(U)\lbrack -2n \rbrack = \map_{\sh(X)}\big(M_X(U), \mf{G}_X^{\leq m}\big)\]
 on \(\mathcal{T}_V\). thus their adjoints under \(\lan_\xi \dashv \xi^*\), morphisms in \(\fun(\Delta^\op, \Sp)\), are equivalent as well. The left adjoint of \(C_{i_*} \circ \xi^*\beta \circ \tilde{\alpha}\) is \(C_{i_*} \circ \beta \circ \alpha\). The Beck-Chevalley map is the left adjoint to the map \(\jmath^*\) of the unit \(\eta_F \colon F \to \pi^* \lan_\pi F\), meaning \(\jmath^* \Phi\) composed with the Beck-Chevalley map is adjoint to \(j_V^*(\pi^*\Phi \circ \eta_F) = j_V^* \phi\). Thus taking geometric realizations gives the desired equivalence.
\end{proof}

\begin{proposition}\label{prop:TL}
Let \(V\), \(B\), \(c\), and \(E\) be as in \ref{lemm:LvT}. There is a commutative diagram of spectra,
\[\begin{tikzcd} [row sep = large, column sep = large]
\lvert F_n(V, \bullet) \rvert \ar[d] \ar[r, "\Xi_V"] & \bb{Z}(n)^\bm(V)\lbrack -2n \rbrack \ar[d]\\
\lvert E^{(c)}(V, \bullet)\rvert \ar[r] & E(V)
\end{tikzcd}\]
with lower horizontal Levine's homotopy coniveau tower morphism and vertical morphisms equivalences. This is natural in restriction given by \(V' \subset V\).
\end{proposition}

\begin{proof}
Let \(\mathcal{P}_{U|V}\) be the category of pairs \((Y, W)\) such that \(Y \in \sm_U\) and \(W\) is closed in \(Y \times_U V\), with the same morphisms as in \(\mathrm{Pair}_X\). We build a natural transformation of functors \(\mathcal{P}_{U|V} \to \Sp\), 
\[\map_{\sh(X)}\big(M_X(Y/W), \mf{G}_X\big) \simeq \map_{\sh(B)}\big(M_B(Y \times_U V/W), s^c \Kgl_B\big),\]
with a succession of maps: \(M_X \simeq g_\sharp M_U\), so adjointness gives an equivalence 
\[\map_{\sh(X)}\big(g_\sharp M_U(Y/W), \mf{G}_X\big) \simeq \map_{\sh(U)}\big(M_U(Y/W), g^* \mf{G}_X\big)\]
\(W \subset Y \times_U V\), meaning \(M_U(Y/W) \in \sh_V(U)\) by the localization sequence, and so the unit 
\[M_U(Y/W) \to i_*i^! M_U(Y/W)\]
is a natural equivalence. Applying this unit in the source and using adjointness again equates this to 
\[\map_{\sh(V)}\big(i^*M_U(Y/W), i^! \circ g^* \mf{G}_X\big)\]
\(i^* M_U( -) \simeq M_V( - \times_U V)\), so \(W \subset Y \times_U V\) implies that 
\[i^* M_U\big(Y/W\big) \simeq M_V\big(Y \times_U V/W\big),\]
and 
\[i^! g^*\mf{G}_X = i^! g^* f^! s^{- n} \Kgl_\bb{Z} \simeq (f \circ g \circ i)^! s^{- n} \Kgl_\bb{Z} \simeq s^c \Kgl_V,\]
the final equivalence by \ref{lemm:LvT} (1) with \(q = 0\). And \(f_V^* s^c \Kgl_B \simeq s^c \Kgl_V\), so using the adjunction \((f_V)_\sharp \dashv f_V^*\) and \((f_V)_\sharp M_V \simeq M_B\). Replacing the source and target with these equivalent objects gives the desired expression. Mapping from a cofibre gives a fibre, so 
\[\map_{\sh(B)}\big(M_B(Y/W), s^c \Kgl_B\big) \simeq E^W(Y)\]
So we have a natural equivalence between the colimitands of \(F_n(V, \bullet)\) and \(E^{(c)}(V, \bullet)\), and taking colimits over \(W\) gives the left vertical equivalence. Restriction along the functor \(\mathcal{T}_V \to \mathcal{P}_{U|V}\), \((q, W) \mapsto (Y, W)\) gives \(F_n(V, q)\) on the spectral replacement side and \(E^W(V \times \Delta^q)\) on the Levine side. Likewise, restriction along the functor \(\Delta^\op \to \mathcal{P}_{U|V}\), \(q \mapsto (Y, V \times \Delta^q)\) gives \(\mf{G}_V(q)\) and \(E(V \times \Delta^q)\), and on the object \((U, V)\) gives \(\bb{Z}(n)^\bm(V)\lbrack -2n \rbrack\) and \(E(V)\), respectively. The morphism of pairs 
\[(Y, W) \to (Y, V \times \Delta^q) \to (U, W)\]
is precisely the morphism inducing 
\[F_n(V, q) \to \bb{Z}(n)^\bm(V \times \Delta^q)\lbrack -2n \rbrack \to \bb{Z}(n)^\bm(V)\lbrack -2n \rbrack,\]
and the composition on the Levine side is 
\[E^W(V) \to E(V \times \Delta^q) \to E(V),\]
Naturality of the equivalence on \(\mathcal{P}_{U|V}\) gives compatibility with restriction along functors into \(\mathcal{P}_{U|V}\), so taking colimits and geometric realizations gives the diagram. Naturality of restricting to open \(V' \subset V\) is inherited by naturality in morphisms of pairs.
\end{proof}

\begin{theorem}
For any \(m \in \bb{Z}\), the morphism 
\[\gr_m |\Phi| \colon \gr_m |F_n( -, \bullet)| \to \gr_m \bb{Z}(n)^\bm \lbrack -2n \rbrack\]
is a Nisnevich-local equivalence.
\end{theorem}

\begin{proof}
Take a point \(x \in U_{(m)}\), and a decomposable pair \((Z, Z')\), with \(Z\) the disjoint union of \(V_i\) of the form \ref{lemm:sm-red}. Write \(U' = U - Z'\). \eqref{prop:int-fun} gives a commutative diagram 
\[\begin{tikzcd} [row sep = large, column sep = huge]
\lvert F_n(V_i), \bullet \rvert \ar[d, "\lvert j_{V_i} \rvert_*", swap] \ar[r, "\Xi_{V_i}"] & \bb{Z}(n)^\bm(V_i)\lbrack -2n \rbrack \ar[d, "i_*"]\\
\fil_m \lvert F_n(U', \bullet) \rvert \ar[r, "\lvert \Phi \rvert^m(U')", swap] & \fil_m \bb{Z}(n)^\bm(U') \lbrack - 2n \rbrack
\end{tikzcd}\]
Taking colimits and using cofinality of decomposable pairs gives the square in \ref{eq:shell-decomp}, identifying the upper morphism in that square as 
\[\bigoplus_{x \in U_{(m)}} \colim_{V \subset Z(x)} \Xi_V\]
So we must show that this morphism is an equivalence of spectra. If \(c \leq 0\), then for each \(q \geq 0\), every closed \(W \subset V \times \Delta^q\) lies in \(\mathcal{S}_V^{(c)}(q)\), so \(|\epsilon| \colon |E^{(c)}(V \times \Delta^\bullet)| \to |E( V \times \Delta^\bullet)|\) is an equivalence.

So assume that \(c > 0\). If \(x\) lies over the prime \(p\), then \(V\) is smooth over \(B = \spec \bb{F}_p\). \(s^c \Kgl_{\bb{F}_p}\) is \(c\)-effective in \(\sh(\bb{F}_p)\), so \(E = \omega^\infty s^c \Kgl_{\bb{F}_p}\) is \(c\)-effective in \(\sh_{S^1}(\bb{F}_p)\) by \cite[Theorem 9.0.3]{LevHC}, meaning \(|E^{(c)}(V, \bullet)| \to E(V)\) is an equivalence by \cite[Lemma 7.3.1]{LevHC}.\footnote{\(E\) is, by definition, the \(\bb{P}^1\)-delooping of a motivic spectrum, so satisfies the requisite axiom A3 for the finite field \(\bb{F}_p\) by \cite[Remark 2.1.1]{LevHC}} This means that \(\Xi_V\) is an equivalence by \ref{prop:TL}.

If \(x\) lies over \(0\), then \(V\) is smooth over \(B = \spec \bb{Z}[1/N]\), for some \(N \in \bb{N}\). Consider the commutative diagram with vertical equivalences 
\[\begin{tikzcd} [row sep = large, column sep = large]
\displaystyle \colim_{V \subset Z(x)} \left \lvert F_n(V, \bullet) \right \rvert \ar[d] \ar[r, "\colim \Xi_V"] & \displaystyle \colim_{V \subset Z(x)} \bb{Z}(n)^\bm(V) \lbrack -2n \rbrack \ar[d]\\
\displaystyle \colim_{V \subset Z(x)} \left \lvert E^{(c)}(V, \bullet) \right \rvert \ar[r] & \displaystyle \colim_{V \subset Z(x)} E(V)
\end{tikzcd}\]
The colimit along \(V \subset Z(x)\) can be computed as the double colimit over \(N\) and \(V\) of \(V[1/N]\). Take a smooth scheme \(f_S \colon S \to V[1/N]\). Then 
\[\omega^\infty s^c \Kgl_{V[1/N]}(S) = \map_{\sh(V[1/N])}(M_{V[1/N]}(S), s^c \Kgl_{V[1/N]}) \simeq \map_{\sh(V[1/N])}\big(\mathbbm{1}_{V[1/N]}, f_S^* s^c \Kgl_{V[1/N]}\big),\]
so once again exploiting continuity of \(\sh\) \cite[1.3.9]{BonDeg}, 
\[\colim_V \colim_N \omega^\infty s^c \Kgl_{V[1/N]}(S) \simeq \colim_V s^c \Kgl_{V_\bb{Q}}(S),\]
meaning 
\[\colim_{V \subset Z(x)} E(V) \simeq \colim_V E_\bb{Q}(V_\bb{Q}),\]
in which \(E_\bb{Q} = \omega^\infty s^c \Kgl_{V_\bb{Q}}\). Similarly, by compatibility of colimits with geometric realization and with other colimits, we have that 
\[\colim_{V \subset Z(x)} |E^{(c)}(V, \bullet)| \simeq \colim_V \colim_N |E^{(c)}(V, \bullet)| \simeq \colim_V \left|\colim_N E^{(c)}(V, \bullet) \right| \simeq \colim_V \left|E_\bb{Q}^{(c)}(V_\bb{Q}, \bullet)\right|\]
The morphism \(|E_\bb{Q}^{(c)}(V_\bb{Q}, \bullet)| \to E_\bb{Q}(V_\bb{Q})\) is an equivalence by Levine \cite[Lemma 7.3.1]{LevHC} again, this time for \(k = \bb{Q}\).

In either case of \(x \in U_{(m)}\), \(\colim_V \Xi_V\) is an equivalence, so that the upper horizontal morphism of \ref{eq:shell-decomp} is an equivalence. The vertical legs of this map are euqivalences, meaning \(|\Phi|^{\lceil m \rceil}(U)\) is an equivalence for each \(U \in X_\nis\). The vertical legs in \eqref{eq:niv-shell} are at least a Nisnevich-local equivalences, meaning the top, \(\gr_m |\Phi|\), is Nisnevich-local.
\end{proof}

So combining equivalences \(A\) and \(B\), \ref{thm:equiv-A} and \ref{thm:equiv-B}, we obtain our crown equivalence, \ref{main-intro} (5).

\bibliographystyle{plain}
\bibliography{references}
\end{document}